\documentclass[11pt]{article}
\usepackage[a4paper,margin=2.6cm]{geometry}
\usepackage{amsmath,amssymb,amsthm,mathtools,booktabs,array,enumitem}
\setlist[enumerate,1]{label=(\arabic*)}
\usepackage{xcolor}
\usepackage[colorlinks=true,linkcolor=blue!50!black,citecolor=blue!50!black,urlcolor=blue!50!black]{hyperref}
\newtheorem{theorem}{Theorem}[section]
\newtheorem{proposition}[theorem]{Proposition}
\newtheorem{lemma}[theorem]{Lemma}
\newtheorem{corollary}[theorem]{Corollary}
\newtheorem{conjecture}[theorem]{Conjecture}
\theoremstyle{definition}
\newtheorem{definition}[theorem]{Definition}
\newtheorem{remark}[theorem]{Remark}
\newtheorem{observation}[theorem]{Observation}

\newcommand{\C}{\mathbb{C}}
\newcommand{\Z}{\mathbb{Z}}
\newcommand{\Q}{\mathbb{Q}}
\newcommand{\F}{\mathbb{F}}
\newcommand{\ch}{\operatorname{ch}}
\newcommand{\Hg}{\mathrm{Hg}}
\newcommand{\bT}{\overline{T}}
\newcommand{\bW}{\overline{W}}
\title{A vertex operator algebra for the S-fold realisation\\ of a non-unitary Haagerup--Izumi modular fusion category}
\author{David E.\ Evans\\ \small School of Mathematics, Cardiff University,\\ \small Cardiff CF24 4AG, Wales, UK\\ \small \texttt{EvansDE@cardiff.ac.uk}}
\begin{document}
\maketitle

\begin{abstract}
Gang--Kim and Gang--Kim--Lee proposed non-unitary Haagerup--Izumi modular data for the topological twist of the 3D $\mathcal N=4$ rank-$0$ S-fold SCFTs $\mathcal S_k$, with characters of a putative boundary rational CFT; Evans derived the modular data from the monodromy word, up to sign and under stated non-degeneracy conditions. We identify the vertex operator algebra behind these characters: the extension
\[
V_k=\bigoplus_{j\ge0}L\big(c_k,h_{2j+1,2j+1}\big),\qquad c_k=13-6k,
\]
of the Virasoro algebra, where $L(c,h)$ is the irreducible Virasoro module of lowest weight $h$ and $h_{r,s}$ are the Kac weights at $c=c_k$. It matches Flohr's $\mathcal W(2,2(k-2))$. For $k=11$ it is strongly rational, and its modular data is the theta-function modular data of Flohr's characters, which agrees numerically with the proposed one; the same computations prove strong rationality for $5\le k\le10$ and give the modular data for $k\in\{5,6,7,9\}$, where two steps are numerical. At the Haagerup point $k=11$, $c=-53$, we show that $V_{11}\oplus M_{\chi_3}$ is a simple vertex operator algebra, unique up to isomorphism, where $M_{\chi_3}$ is the simple $V_{11}$-module of lowest weight $0$ other than $V_{11}$. Modular data alone does not determine a modular tensor category; the extension supplies the missing structure. It is a connected étale algebra in $\mathrm{Rep}\,V_{11}$, and a condensation of the fusion category of its modules has the rank-$12$ factor of $\mathrm{Rep}\,V_{11}$ as its Drinfeld centre. With the classification of twisted cubic trivalent categories of Morrison--Peters--Snyder, we deduce $\mathrm{Rep}\,V_{11}\simeq\mathrm{Sem}\boxtimes\mathcal Z(\sigma_2\mathcal F)$ as braided categories, where $\mathrm{Sem}$ is the semion category (with the spherical structure in which its generator has dimension $-1$), $\mathcal F$ is Ostrik's twisted Haagerup category and $\sigma_2$ is the Galois automorphism $\zeta_{117}\mapsto\zeta_{117}^2$. 
Here $\mathcal Z(\mathcal F)$ is a unitary modular tensor category with the twisted Evans–Gannon modular data $\mathcal D^{\omega}\Hg_3$, $\omega=2$, while its Galois conjugate $\mathcal Z(\sigma_2\mathcal F)$ is not unitary, having simple objects of negative dimension. Evans–Gannon’s conjecture on the double of the $\Z_9+9$ near-group category follows.

\end{abstract}

\tableofcontents

\section{Introduction}

The modular data of the Drinfeld double of the Haagerup subfactor lies in the family $\mathcal D^\omega\Hg_{2m+1}$ of Evans--Gannon, labelled by an odd integer $2m+1$ and a twist $\omega$; the double of the Haagerup subfactor is $\mathcal D^{0}\Hg_3$~\cite{EG11}. There it is a \emph{grafting} of two finite metric groups~\cite[\S3.3]{EG11}. For the untwisted family the grafting description was proved by Gannon--Schopieray--Yadav~\cite[Theorem~5.1]{GSY}. Gang--Kim~\cite{GK} proposed a family of non-unitary modular data, labelled by an integer $k\ge3$, for the topological twist $\mathrm{TFT}_-[\mathcal S_k]$ of the S-fold SCFT
\[
\mathcal S_k=T[SU(2)]\big/SU(2)^{\rm diag}_k ,
\]
obtained from the 3D $\mathcal N=4$ theory $T[SU(2)]$ by gauging the diagonal $SU(2)$ of its $SU(2)\times SU(2)$ symmetry with Chern--Simons level $k$.
Their input is the $T$-matrix and the vacuum column $|S_{0\alpha}|$, obtained from the Bethe vacua of an abelianised gauge theory; the full $S$-matrix is then proposed by the Evans--Gannon grafting and checked against the modular relations~\cite{GK},~\cite[\S1.5]{PaperA}. This is not a realisation of the modular data. For $k=4m^2+4m+3$ the proposed data is a Galois conjugate of $\mathcal D^{\omega=2}\Hg_{2m+1}$, modulo a decoupled semion~\cite{GK,GKL}. Gang--Kim--Lee~\cite{GKL} use the same $S$, and compute half-indices of an abelian dual of $\mathcal S_k$. These are Nahm sums, which they complete to vector-valued modular functions by the Bantay--Gannon method~\cite{BG}. Jeong--Sohn~\cite{JS} proceed similarly for generalised S-folds, i.e.\ once-punctured torus bundles with monodromy $\varphi=ST^{k_1}\cdots ST^{k_n}$, with Bethe vacua found by ans\"atze for particular families. Neither paper identifies a vertex operator algebra.
In fact the characters of~\cite{GKL} agree, to the orders checked (\S\ref{sec:chars}), with characters obtained in 1993 by Flohr~\cite[Table~4.2]{Flohr} for his algebra $\mathcal W(2,8\kappa)$, $\kappa=(k-2)/4$, constructed from screened vertex operators, together with the modular transformations of the theta functions involved~\cite[(4.5)]{Flohr}.

Evans~\cite{PaperA} treats the S-fold side rigorously, for every hyperbolic word $\varphi$ in $L=\begin{psmallmatrix}1&1\\0&1\end{psmallmatrix}$, $R=\begin{psmallmatrix}1&0\\1&1\end{psmallmatrix}$, i.e.\ a cyclic word containing both letters; the S-fold $\mathcal S_k$ has $\varphi=LR^{k-2}$, which is conjugate in $\mathrm{SL}(2,\Z)$ to its monodromy $\begin{psmallmatrix}0&1\\-1&0\end{psmallmatrix}R^k$, in the conventions of~\cite[(A.3)]{GKL}. The two metric groups are $A_\mp=\mathrm{coker}(\varphi\mp1)$ with the torsion linking forms of the mapping tori; the Bethe data are given in closed form by the word; and $S$ is derived and shown to be, up to sign, the Evans--Gannon grafting, under the stated non-degeneracy conditions. 

Three recent constructions~\cite{Jiao,Riesen,HuangTC} realise the untwisted double $\mathcal Z(\mathcal H_3)$ of the Haagerup fusion category $\mathcal H_3$ by a unitary vertex operator algebra at $c=8$ inside $V_{E_6\oplus A_2}$, as predicted by Evans--Gannon~\cite[\S2.4]{EG11}.
It is also observed in~\cite[Remark~3.4]{PaperA} that Gang--Kim's family proposes only $\omega=2$ in $\mathcal D^\omega\Hg_{2m+1}$, i.e.\ the metric group $\Z_{(2m+1)^2}$ rather than $\Z_{2m+1}^2$ ($\Z_9$ rather than $\Z_3\times\Z_3$ for $m=1$), and misses the untwisted double. The present paper is concerned with the vertex operator algebra behind the non-unitary data and is independent of the $c=8$ constructions.

We will use 
the extension
\[
V_k=\bigoplus_{j\ge0}L\big(c_k,h_{2j+1,2j+1}\big),\qquad c_k=13-6k,
\]
of the Virasoro algebra. Here $L(c,h)$ is the irreducible Virasoro module of central charge $c$ and lowest weight $h$, and $h_{r,s}$ are the Kac weights; in the parametrisation $c=13-6(t+t^{-1})$ of \S\ref{sec:not}, $c_k$ corresponds to $t+t^{-1}=k$. The algebra $V_k$ is strongly generated (\S\ref{sec:not}) by the conformal vector $T$ and one Virasoro primary $W$ of weight $2(k-2)$; it matches Flohr's $\mathcal W(2,8\kappa)$, $\kappa=(k-2)/4$, of 1993~\cite{Flohr} (\S\ref{sec:known}). At the Haagerup point $k=11$, $W$ has weight $18$ and $t=d^2$, where $d=(3+\sqrt{13})/2$ is the dimension of the non-invertible simple objects of the Haagerup category. The Virasoro parameter has a heuristic physical explanation (\S\ref{sec:phys}): $t$ is the fixed point of the duality monodromy of the S-fold, and $V_k$ arises as the cyclic closure of the Creutzig--Gaiotto gluing condition~\cite{CG} for chains of S-duality walls. Its central charge agrees with the congruence $c\equiv\#L-6\#R\pmod{24}$ of~\cite{PaperA}, where $\#L$, $\#R$ are the numbers of letters $L$, $R$ in $\varphi$.

For $k=11$, $V_k$ is $C_2$-cofinite and Zhu's algebra is semisimple, so $V_k$ is strongly rational by McRae's criterion~\cite{McRae}, and a bound on $\dim A(V_k)$, the cyclic span of the vacuum character under $SL_2(\Z)$ and a commutant computation identify the $2k+2$ simple modules, their characters and the modular data: the $S$-matrix is the theta-function $S$-matrix $S_\theta$ of Flohr's characters. Numerically, this agrees with the data proposed by Gang--Kim and Gang--Kim--Lee and derived, up to sign, in~\cite{PaperA}. The same computations prove strong rationality for $5\le k\le10$, and give the simple modules, characters and modular data for $k\in\{5,6,7,9\}$, where two steps (Lemma~\ref{lem:cyc} and the commutant step of Theorem~\ref{thm:mod}) are numerical.

At $k=11$, $\mathrm{Rep}\,V_{11}$ has rank $24$ and splits as a semion, generated by the module $M_{\chi_1}$ (of dimension $-1$ and twist $-i$) with the character of the even diagonal sector $\bigoplus_{m\ge1}L(h_{2m,2m})$, times a rank-$12$ factor $\mathcal C$ with modular data $\sigma_2(\mathcal D^{\omega=2}\Hg_3)$ (with $S$ up to an overall sign), where $\sigma_2$ is the Galois automorphism $\zeta_{117}\mapsto\zeta_{117}^2$. (Since complex conjugation is also a Galois automorphism and exchanges $\omega=1,2$, naming $\sigma_2$ is what makes $\omega=2$ meaningful.)  Modular data does not determine a modular category, so the question is which category $\mathcal C$ is. A vertex operator algebra extension supplies the missing structure. Write $M_\chi$ for the simple $V_k$-module with character $\chi$. Of the $2k+2$ characters (\S\ref{sec:chars}), two have lowest weight $0$: $\chi_0=\ch V_k$ and $\chi_3$. By Huang--Kirillov--Lepowsky~\cite{HKL}, $V_{11}\oplus M_{\chi_3}$ is a simple vertex operator algebra exactly when $A_2=\mathbf 1\oplus M_{\chi_3}$ is a commutative algebra in $\mathrm{Rep}\,V_{11}$ with trivial twist that is simple as a module over itself; as $\dim A_2\ne0$, it is then étale~\cite[Lemma~1.20]{KO}. Such an algebra gives a fusion category $\mathcal C_{A_2}$ of $A_2$-modules, with $\mathcal Z(\mathcal C_{A_2})\simeq\mathcal C\boxtimes(\mathcal C_{A_2}^{\rm loc})^{\rm rev}$~\cite{DMNO}. Here $\mathcal C_{A_2}$ has the fusion rules of the near-group $\Z_9+9$ and $\mathcal C_{A_2}^{\rm loc}$ is pointed; condensing a $\Z_3$ in the second factor gives a fusion category $\mathcal F_{\mathcal C}$ with $\mathcal Z(\mathcal F_{\mathcal C})\simeq\mathcal C$ (Theorem~\ref{thm:Ctw}). The problem is thereby reduced to identifying $\mathcal F_{\mathcal C}$, a fusion category of rank $6$. A bootstrap of $V_{11}\oplus M_{\chi_3}$, combined with the Fuchs relation (the constraint on the sum of the local exponents of a Fuchsian differential equation) for one four-point function, shows that it is a simple vertex operator algebra, unique up to isomorphism. Hence $\mathcal C$ is the centre of a fusion category with Haagerup--Izumi fusion rules for $\Z_3$ and nontrivial pointed associator, i.e.\ the associator restricted to the invertible objects is a nontrivial $3$-cocycle on $\Z_3$. By the classification of twisted cubic trivalent categories of Morrison--Peters--Snyder~\cite{MPS}, this category is a Galois conjugate of Ostrik's twisted Haagerup category $\mathcal F$, obtained by condensing the bosons $a_g$, $g\in3\Z_9$ (invertible objects of trivial twist), in the centre of the near-group category $\Z_9+9$, whose simple objects are $g\in\Z_9$ and one object $\rho$ with $\rho^{\otimes2}=\bigoplus_g g\oplus9\rho$. Thus $\mathcal C\simeq\mathcal Z(\sigma_2\mathcal F)$. The rank-$12$ factor of Evans--Gannon's modular data for $\Z_9+9$ is exactly $\mathcal D^{\omega}\Hg_3$, $\omega\not\equiv0\pmod3$, and the identification gives~\cite[Conjecture~2]{EG14} for $\Z_9$.

With $n:=k-2$ and $h:=2n$, these are the main results:

\begin{enumerate}[label=(\Alph*)]
\item \textbf{The algebra} (\S\ref{sec:alg}). $V_k=\bigoplus_{j\ge0}L(c_k,h_{2j+1,2j+1})$, $c_k=13-6k$ (\S\ref{sec:not}). It exists as a simple vertex operator algebra and is strongly generated by $T$ and a primary $W$ of weight $h$.

\item \textbf{The Haagerup point} (\S\ref{sec:k11}). $V_{11}$ is strongly rational, and $\mathrm{Rep}\,V_{11}\simeq\mathrm{Sem}\boxtimes\mathcal C$ as braided categories, where the semion is generated by $M_{\chi_1}$, which has the character of the even diagonal sector $\bigoplus_{m\ge1}L(h_{2m,2m})$. The factor $\mathcal C$ has rank $12$ and modular data $\sigma_2(\mathcal D^{\omega=2}\Hg_3)$, with $S$ up to an overall sign.

\item \textbf{$\mathcal C$ is the centre of a twisted Haagerup--Izumi category} (\S\S\ref{ssec:ng}--\ref{ssec:gap}).
\begin{itemize}
\item $U_2=V_{11}\oplus M_{\chi_3}$ is a simple vertex operator algebra, unique up to isomorphism (Theorem~\ref{thm:U2}). So $V_{11}\oplus M_{\chi_3}$ is an étale (commutative, separable) algebra in $\mathcal C$~\cite{KO}, and $\mathcal C$ is a Drinfeld centre (Corollary~\ref{cor:centre}). The proof combines a bootstrap of $U_2$ (Proposition~\ref{prop:boot}, Lemma~\ref{lem:G}) with the Fuchs relation for the four-point function of its weight-$0$ field.
\item $\mathcal C\simeq\mathcal Z(\sigma_2\mathcal F)$, where $\mathcal F$ is Ostrik's twisted Haagerup category, obtained by condensing $\bigoplus_{h\in3\Z_9}a_h$ in the centre of the near-group category $\Z_9+9$~\cite{MPS},~\cite[\S7.1]{EG17}, and $\sigma_2$ is a Galois automorphism (Theorem~\ref{thm:ident}). The identification uses the classification of twisted cubic trivalent categories of Morrison--Peters--Snyder~\cite{MPS}.
\item Consequently the rank-$12$ factor of the double of $\Z_9+9$ has the modular data $\mathcal D^{\omega}\Hg_3$, $\omega\not\equiv0\pmod3$, as predicted by Evans--Gannon's Conjecture~2~\cite{EG14} (which states it as $MD_{\Z_9,\Z_{13}}$), and this data is realised by the unitary modular tensor category $\mathcal Z(\mathcal F)$ (Corollary~\ref{cor:EGconj}).
\end{itemize}
\end{enumerate}
The proofs of (B) and (C) are computer-assisted. They rest on the following two results, which hold for a range of $k$.
\begin{enumerate}[label=(\Alph*),resume]
\item \textbf{$C_2$-cofiniteness} (\S\ref{sec:c2}). For $5\le k\le 11$, the $C_2$-algebra $R_{V_k}=V_k/C_2(V_k)$ is finite-dimensional, with $\bT^{\,2k-2}=0$, $\bT^{\,k}\bW=0$ and $\dim R_{V_k}\le 3k-3$, where $\bT,\bW$ are the images of $T,W$; for $k\le7$ equality holds modulo the primes used. For $V_{11}$: $\bT^{20}=\bT^{11}\bW=0$ and $\dim R\le 30$. These finiteness statements hold over $\Q$: the computation is modulo a prime, with good reduction (the relevant Gram matrices have the same rank modulo the prime as over $\Q$) checked at every level used.
\item \textbf{Zhu's algebra and rationality} (\S\S\ref{sec:zhu}--\ref{sec:rat}). For $5\le k\le11$, Zhu's algebra $A(V_k)$ is semisimple of dimension $\le 2k+2$, so $V_k$ is strongly rational by McRae's criterion~\cite{McRae} and $\mathrm{Rep}\,V_k$ is a modular tensor category. For $k\in\{5,6,7,9,11\}$ there are exactly $2k+2$ simple modules, with Flohr's characters and $S$-matrix $S_\theta$; these agree numerically with the data proposed by Gang--Kim--Lee. For $k=11$ all steps are exact; for $k\in\{5,6,7,9\}$ two steps (Lemma~\ref{lem:cyc} and the commutant step of Theorem~\ref{thm:mod}) are numerical.
\end{enumerate}
For $k=11$ every computer-assisted step is an exact computation or a rank statement modulo a prime with good reduction checked; Appendix~\ref{app:comp} summarises them.

\paragraph{Relation to earlier work.} $V_k$ matches Flohr's $\mathcal W(2,8\kappa)$, $\kappa=n/4$~\cite{Flohr}: for every $k$ the two agree in central charge, generator weights, vacuum character and the list of $2k+2$ characters of lowest-weight representations, and for $3\le k\le11$ also in the self-coupling of $W$, the one structure constant not fixed by $c$ (Table~\ref{tab:N}; \S\ref{sec:known}). Flohr's characters are theta quotients, and the half-indices of Gang--Kim--Lee~\cite[(3.12)]{GKL} agree with them to the orders checked ($q^{34}$ for the vacuum, $3\le k\le11$; $q^{40}$ for the other three, $4\le k\le7$; \S\ref{sec:chars}). For $k\le5$, $V_k$ is a tensor product of minimal models or a principal $\mathcal W$-algebra; for $k=4,5$ this assumes the uniqueness of the simple $\mathcal W(2,4)$ at $c=-11$ and of the simple $\mathcal W(2,6)$ at $c=-17$ (\S\ref{sec:known}).

\paragraph{Physical origin.} The parameter $t$ is the fixed point of the monodromy $T^kS$ on the Kapustin--Witten coupling $\Psi$, and $V_k$ is the cyclic closure of the Creutzig--Gaiotto gluing condition for chains of S-duality walls. For a general word $\varphi$ this gives an algebra $V_\varphi$ with $c_\varphi=\#L-6\#R$, in agreement with~\cite[Theorem~9.5]{PaperA}, and for the one-node word $\varphi=LR^{k-2}$ the effective central charge $c-24h_{\min}$, with $h_{\min}=-n/4$ the least lowest weight among Flohr's characters, equals $\#L=1$, as in~\cite[Theorem~8.2]{PaperA}. This is a derivation sketch (\S\ref{sec:phys}); nothing in the proofs depends on it.

\paragraph{A conjecture for odd $k$} (\S\ref{sec:conj}). For odd $k$ we conjecture that $\mathrm{Rep}\,V_k$ is a semion times a Galois conjugate of the factor of the double of the near-group $\Z_{k-2}+(k-2)$, and that $V_k\oplus M_{\chi_3}$ is a vertex operator algebra. Assuming the conjecture, this factor is a Drinfeld centre exactly at the Haagerup--Izumi values $k=\nu^2+2$, $\nu$ odd, where it should be the centre of a twisted Haagerup--Izumi category for $\Z_\nu$.

\section{Notation}\label{sec:not}

In a vertex operator algebra $V$ with central charge $c$, $T:=L_{-2}\mathbf 1$ denotes the conformal vector, whose field $T(z)=\sum_nL_nz^{-n-2}$ generates the Virasoro algebra. $V$ is \emph{strongly generated} by vectors $a^1,\dots,a^r$ if it is spanned by the normally ordered products $a^{i_1}_{(-n_1)}\cdots a^{i_s}_{(-n_s)}\mathbf 1$ with all $n_j\ge1$. The letter $T$ is also used for the modular $T$-matrix and, in \S\ref{sec:phys}, for the generator $\begin{psmallmatrix}1&1\\0&1\end{psmallmatrix}$ of $\mathrm{SL}(2,\Z)$; the meaning will be clear from context.

The Virasoro central charge is parametrised as $c=13-6(t+t^{-1})$, and
\[
h_{r,s}=\tfrac14\big((r^2-1)t+(s^2-1)t^{-1}-2(rs-1)\big).
\]
A vertex operator algebra $V$ is graded by the eigenvalues of $L_0$, the generator of rotations (the conformal Hamiltonian): $V=\bigoplus_{m\ge0}V_m$ with $\dim V_m<\infty$. A vector in $V_m$ has \emph{weight} $m$, and \emph{below weight $N$} refers to the subspace $\bigoplus_{m<N}V_m$. The same applies to modules, whose weights lie in $h+\Z_{\ge0}$ for a lowest weight $h$.

For $t\notin\Q$ the Verma module $M(c,h_{r,s})$ ($r,s\ge1$) has, up to scalars, exactly one singular vector other than its highest-weight vector, at level $rs$, and the submodule it generates is the irreducible Verma module $M(c,h_{r,s}+rs)$~\cite{FeFu}. (Indeed $h_{r,s}=h_{r',s'}$ forces $rt-s=\pm(r't-s')$, which for $t\notin\Q$ leaves only $(r',s')=(r,s)$.) Hence
\[
\operatorname{ch}L(c,h_{r,s})=\frac{q^{h_{r,s}-c/24}\,(1-q^{rs})}{(q)_\infty}.
\]
We write $\theta_{a,N}(\tau)=\sum_{m\in\Z}q^{(a+2Nm)^2/4N}$ and $\Lambda_{a,N}=\theta_{a,N}/\eta$. Throughout, $k\ge3$, $n=k-2$ and $h=2n$.

If $t+t^{-1}=k$ with $k\ge3$, then $t$ is a real quadratic irrational and
\[
h_{r,r}=\tfrac14(r^2-1)(k-2),\qquad c_k=13-6k=1-6n.
\]

For a partition $\lambda=(\lambda_1\ge\lambda_2\ge\dots\ge\lambda_\ell\ge1)$ write $L_{-\lambda}=L_{-\lambda_1}L_{-\lambda_2}\cdots L_{-\lambda_\ell}$ and $|\lambda|=\sum_i\lambda_i$. By the Poincar\'e--Birkhoff--Witt (PBW) theorem the \emph{PBW monomials} $L_{-\lambda}v$ form a basis of the Verma module $M(c,h)$ with highest weight vector $v$, and $L_{-\lambda}v$ has weight $h+|\lambda|$; in the vacuum module $L(c,0)$, a basis is given by those with all parts $\lambda_i\ge2$.

$M(c,h)$ carries a unique contravariant (Shapovalov) form with $\langle v,v\rangle=1$ and $\langle L_na,b\rangle=\langle a,L_{-n}b\rangle$; it is symmetric and bilinear, and in general not positive definite. The \emph{Gram matrix} at level $\ell$ is $G_\ell=\big(\langle L_{-\mu}v,L_{-\lambda}v\rangle\big)_{|\mu|=|\lambda|=\ell}$. Its determinant is the Kac determinant, and its radical is the maximal proper submodule of $M(c,h)$ at level $\ell$; so $L(c,h)$ at level $\ell$ is $M(c,h)_\ell$ modulo the kernel of $G_\ell$.

Every fusion category is defined over a number field, by Ocneanu rigidity~\cite{ENO}. For a field automorphism $\sigma$ of a subfield of $\overline{\Q}$ and a fusion category $\mathcal A$, we write $\sigma\mathcal A$ for the category obtained by applying an extension of $\sigma$ to $\overline{\Q}$ to the structure constants of $\mathcal A$ in some field of definition; for the categories considered here, the result does not depend on the extension (Theorem~\ref{thm:ident}).

An algebra $A$ in a braided fusion category is \emph{étale} if it is commutative and separable (the multiplication splits as a map of $A$-bimodules), and \emph{connected} if $\dim\mathrm{Hom}(\mathbf 1,A)=1$~\cite{DMNO}. In a unitary braided category these are the irreducible commutative Q-systems~\cite{BEK1}. By Huang--Kirillov--Lepowsky~\cite{HKL}, a vertex operator algebra extending $V$ is the same as a commutative algebra in $\mathrm{Rep}\,V$ with trivial twist; it is simple exactly when the algebra is simple as a module over itself, and the algebra is connected when $V$ occurs in it with multiplicity one. If moreover $\dim A\ne0$, which is not automatic in a non-unitary category, the algebra is étale~\cite[Lemma~1.20]{KO}.

\begin{observation}\label{obs:t}
Let $k=N^2+2$, and let $d=\tfrac12\big(N+\sqrt{N^2+4}\big)$, which is the dimension of the non-invertible simple object of the Haagerup--Izumi category for $\Z_N$ ($d^2=1+Nd$). Then $t=d^2$: indeed $d^2+d^{-2}=(d-d^{-1})^2+2=N^2+2=k$. These are exactly the values $k=4m^2+4m+3$ ($N=2m+1$) at which Gang--Kim find Haagerup--Izumi data. For the Haagerup subfactor, $N=3$, $k=11$ and $t=\big(\tfrac{3+\sqrt{13}}2\big)^2$.
\end{observation}

\section{The algebra \texorpdfstring{$V_k$}{V\_k}}\label{sec:alg}

\begin{definition}
$V_k:=\bigoplus_{j\ge0}L(c_k,h_{2j+1,2j+1})$, where $c_k=13-6k$ and $h_{2j+1,2j+1}=j(j+1)n$, $n=k-2$.
\end{definition}
Write $W$ for a nonzero vector spanning the lowest-weight space of the summand $L(c_k,h_{3,3})$; it is a Virasoro primary of weight $h=h_{3,3}=2n$, unique up to scale, and its normalisation is fixed in \S\ref{ssec:model}.
In terms of \S\ref{sec:not} this is $c=13-6(t+t^{-1})$ and the diagonal of the Kac table, with $t+t^{-1}=k$; on the diagonal $t$ enters only through $t+t^{-1}$. The irrationality of $t$ is what is used below: it gives the singular-vector structure of \S2 and the description of the Virasoro category in the proof of Theorem~\ref{thm:exist}.

\begin{theorem}\label{thm:exist}
$V_k$ carries the structure of a simple vertex operator algebra of central charge $c_k$. It is strongly generated by $T$ and the primary $W$ spanning the top of $L(c_k,h_{3,3})$, which has weight $h=2(k-2)$.
\end{theorem}

\begin{proof}
\emph{Existence.} We use the following results on the category $\mathcal O_c$ of $C_1$-cofinite Virasoro modules at central charge $c=13-6(t+t^{-1})$, $t\notin\Q$, due to Creutzig--Jiang--Orosz Hunziker--Ridout--Yang~\cite{CJORY}. $\mathcal O_c$ is a braided tensor category~\cite[Theorem~4.2.6]{CJORY}; it is semisimple, with simple objects the $L(h_{r,s})$, $r,s\ge1$~\cite[Theorem~5.1.2]{CJORY}, rigid and non-degenerate~\cite[Theorem~5.5.3, Proposition~5.5.5]{CJORY}; and $L(h_{r,1})\boxtimes L(h_{1,s})\cong L(h_{r,s})$~\cite[Theorem~5.2.5]{CJORY}. The full subcategory $\mathcal O^L$ generated by the $L(h_{r,1})$ is braided equivalent to the subcategory $(\mathcal{KL}_1\boxtimes\mathcal{KL}_\ell)^0$ of $L_1(\mathfrak{sl}_2)\otimes L_\ell(\mathfrak{sl}_2)$-modules with simple objects $L_1(\bar\mu)\otimes V_\ell(\mu)$, $\bar\mu\equiv\mu\pmod2$, for a generic level $\ell$ determined by $t$~\cite[Proposition~5.5.2]{CJORY}; the same argument applies to the subcategory $\mathcal O^R$ generated by the $L(h_{1,s})$~\cite[proof of Theorem~5.5.3]{CJORY}. On the subcategories of odd $r$ (resp.\ $s$), i.e.\ even $\mu$, the factor $L_1(0)$ is the unit, so these are braided equivalent to the integer-spin part of the Kazhdan--Lusztig category $\mathcal{KL}_\ell(\mathfrak{sl}_2)$, hence~\cite{KL} to the integer-spin part of the semisimple category of type $A_1$ at a generic parameter $q$. By the relation between $\ell$ and $t$~\cite[(5.13)]{CJORY}, the Kazhdan--Lusztig parameter $q=e^{\pi i/(\ell+2)}$ has $q^2=e^{2\pi it}$ for $\mathcal O^L$, and the parameter $q'=e^{\pi i/(\ell'+2)}$ of $\mathcal O^R$, for its generic level $\ell'$, has $q'^2=e^{2\pi i/t}$, in our labelling of the Kac table. Only $q^2$ and $q'^2$ enter below, since $\mathrm{SO}(3)_q$ depends only on $q^2$. The argument below uses only $q'^2=q^{-2}$, which also holds if the two labellings are exchanged, since $t+t^{-1}\in\Z$; consistently, $\theta_{L(h_{3,1})}=e^{2\pi i(2t-1)}=q^4$ and $\theta_{L(h_{1,3})}=e^{2\pi i(2/t-1)}=q'^4$.

Let $\mathcal O^L_{\rm odd}\subset\mathcal O^L$ and $\mathcal O^R_{\rm odd}\subset\mathcal O^R$ be the subcategories with $r$, resp.\ $s$, odd (integer spin). From \S2,
\[
h_{r,s}-h_{r,1}-h_{1,s}=-\tfrac12(r-1)(s-1),
\]
which is an even integer for $r,s$ odd. So the double braiding of $L(h_{r,1})$ with $L(h_{1,s})$ is trivial, the two subcategories centralise each other, and $(X,Y)\mapsto X\boxtimes Y$ is a braided equivalence $\mathcal O^L_{\rm odd}\boxtimes\mathcal O^R_{\rm odd}\to\langle L(h_{r,s}):r,s\ \text{odd}\rangle$. Next, $t+t^{-1}=k\in\Z$ gives $q'^2=e^{2\pi ik}e^{-2\pi it}=q^{-2}$. The integer-spin part of the type-$A_1$ category at a generic parameter $q$ is the idempotent completion of the braided trivalent category $\mathrm{SO}(3)_q$ of Morrison--Peters--Snyder, generated by the spin-$1$ object~\cite[\S4, Example~8.5]{MPS}. Its standard braiding involves only $q^{\pm2}$~\cite[Example~8.5]{MPS}, so $\mathrm{SO}(3)_q$ and $\mathrm{SO}(3)_{-q}$ coincide as braided categories. Moreover $\mathrm{SO}(3)_{q^{-1}}\simeq\mathrm{SO}(3)_q^{\rm rev}$: the reverse of the standard braiding of the type-$A_1$ category at $q$ is the standard braiding at $q^{-1}$ (the universal $R$-matrix $R$ is replaced by $R_{21}^{-1}$), and this restricts to the integer-spin parts. Since $q'^2=q^{-2}$, there is a braided equivalence $\Phi\colon\mathcal O^L_{\rm odd}\to(\mathcal O^R_{\rm odd})^{\rm rev}$ with $\Phi(L(h_{r,1}))=L(h_{1,r})$. It is compatible with the twists: by the displayed identity, $h_{r,1}+h_{1,r}=h_{r,r}+\tfrac12(r-1)^2\in\Z$ for odd $r$, since $h_{r,r}=\tfrac14(r^2-1)n\in\Z$. Hence
\[
A=\bigoplus_{r\ \rm odd}L(h_{r,1})^*\boxtimes\Phi(L(h_{r,1}))\cong\bigoplus_{r\ \rm odd}L(h_{r,r})
\]
(each $L(h_{r,1})$ is self-dual) is the image of the canonical commutative algebra $\bigoplus_X X^*\boxtimes\Phi(X)$ in $\mathrm{Ind}(\mathcal O^L_{\rm odd}\boxtimes\mathcal O^R_{\rm odd})$ attached to the braid-reversed equivalence $\Phi$, as in the gluing construction of Creutzig--Kanade--McRae~\cite{CKM}, under the tensor product functor $\mathcal O^L_{\rm odd}\boxtimes\mathcal O^R_{\rm odd}\to\mathcal O_c$. That functor is braided because the two subcategories centralise each other, and it extends to Ind-objects, so it carries commutative algebras to commutative algebras. (Here, unlike in~\cite{CKM}, both factors are subcategories of the module category of one Virasoro algebra.) So $A$ is a commutative associative algebra in the direct-limit completion $\mathrm{Ind}(\mathcal O_c)$, which carries vertex and braided tensor category structures extending those of $\mathcal O_c$~\cite[Theorem~7.1, Example~7.2]{CMY}; it is connected, with trivial twist, and its multiplication is nonzero in every channel allowed by the fusion rules: the component $(X^*\boxtimes\Phi X)\otimes(Y^*\boxtimes\Phi Y)\to Z^*\boxtimes\Phi Z$ is $\sum_i\beta_i^\vee\boxtimes\Phi(\beta_i)$ over dual bases $\{\beta_i\},\{\beta_i^\vee\}$ of $\mathrm{Hom}(X\otimes Y,Z)$ and its dual, hence nonzero whenever this space is. As a module it is $\Z$-graded by $L_0$, with finite-dimensional graded pieces and weights bounded below ($h_{r,r}=\tfrac14(r^2-1)n\ge0$, and each weight occurs in finitely many summands). These are the hypotheses of~\cite[Theorem~7.5(2)]{CMY}, a version of~\cite[Theorem~3.2]{HKL}: a commutative associative algebra in $\mathrm{Ind}(\mathcal O_c)$, $\Z$-graded by $L_0$ and satisfying the grading-restriction conditions. That theorem makes $V_k=A$ a vertex operator algebra extending $L(c_k,0)$, with the conformal vector of $L(c_k,0)$ and vertex operator the intertwining operator given by the multiplication. It is simple: a nonzero ideal is an $L(c_k,0)$-submodule, so, the summands $L(h_{r,r})$ being pairwise non-isomorphic simple modules, it contains one of them; multiplying by that summand reaches $L(c_k,0)$, because the pairings $L(h_{r,r})\otimes L(h_{r,r})\to\mathbf 1$ are nonzero, so the ideal contains $\mathbf 1$.

\emph{Strong generation.} Since the multiplication is nonzero in every allowed channel, and the fusion of $L(h_{3,3})$ with $L(h_{2j+1,2j+1})$ contains $L(h_{2j+3,2j+3})$, the modes of $W$ applied to $L(h_{2j+1,2j+1})$ reach $L(h_{2j+3,2j+3})$; so $T$ and $W$ generate $V_k$. The singular part of $W(z)W(w)$ has weights $<2h$, so it lies in $L(0)\oplus L(h_{3,3})$, the summand $L(h_{5,5})$ starting at weight $6n>2h-1$ (the dimension argument of~\cite[\S2]{Flohr}). Every vector of $L(0)\oplus L(h_{3,3})$ is a normally ordered polynomial in $T$, $W$ and their derivatives. So the operator product expansions of the generators close on normally ordered polynomials in the generators, and strong generation follows by the standard argument~\cite{Kac}.
\end{proof}

\begin{remark}[A degenerate free-field model]
The subalgebra $\langle T_a,e^{-2a\phi}\rangle$ of a rank-one boson with background charge $a=\sqrt{n/2}$ has the same character as $V_k$ up to weight $16$ (checked for $k=4$). It is not isomorphic to $V_k$: the field $e^{-2a\phi}$ has zero norm, since $W_{(j)}W=0$ for $j\ge-2n$, so that subalgebra is not simple. Character data alone does not determine the vertex algebra. A correct free-field model needs screened vertex operators, as in Flohr's construction~\cite{Flohr}.
\end{remark}


\section{Characters and modular data}\label{sec:chars}

\subsection{The characters}
The characters of $V_k$ and of its lowest-weight representations were found by Flohr~\cite[(4.16), Table~4.2]{Flohr}, in his notation $\kappa=n/4$. The vacuum character is
\[
\operatorname{ch}V_k=\frac{q^{-c_k/24}}{(q)_\infty}\sum_{m\ge0}\Big(q^{(k-2)m(m+1)}-q^{(k+2)m(m+1)+1}\Big)=\tfrac12\big(\Lambda_{n,n}-\Lambda_{n+4,n+4}\big);
\]
the first expression is the sum of the Virasoro characters of \S\ref{sec:not} over the summands $L(c_k,h_{2j+1,2j+1})$. The two metric groups enter through the radii $n=k-2$ and $n+4=k+2$. The vacuum is a \emph{difference}: one group with a plus sign, the other with a minus sign. This is exactly the shape of the grafting~\cite[\S3.3]{EG11}.

There are $2k+2$ characters in all:
\[
\chi_0=\operatorname{ch}V_k,\quad
\chi_3=\tfrac12(\Lambda_{n,n}+\Lambda_{n+4,n+4}),\quad
\chi_1=\tfrac12(\Lambda_{0,n}-\Lambda_{0,n+4}),\quad
\chi_2=\tfrac12(\Lambda_{0,n}+\Lambda_{0,n+4}),
\]
with lowest weights $0$, $0$, $h_{2,2}=3n/4$ and $-n/4$, together with $\Lambda_{\alpha,n}$ ($1\le\alpha\le n-1$) and $\Lambda_{\beta,n+4}$ ($1\le\beta\le n+3$). Here $\chi_1=\operatorname{ch}\bigoplus_{m\ge1}L(h_{2m,2m})$ is the ``even'' diagonal sector. Each $\Lambda_{\alpha,n}$ (resp.\ $\Lambda_{\beta,n+4}$) is a sum of irreducible Virasoro characters whose momenta differ by $2a$ (resp.\ $2b$), where $a^2=n/2$ and $b^2=(n+4)/2$. The total count is $4+(n-1)+(n+3)=2k+2$, which is the rank of $\mathrm{TFT}_-[\mathcal S_k]$.

The $S$-matrix in this basis follows from $\Lambda_{a,N}(-1/\tau)=(2N)^{-1/2}\sum_b e^{-\pi iab/N}\Lambda_{b,N}(\tau)$; Flohr notes that the $S$- and $T$-matrices of $\mathcal W(2,8\kappa)$ follow from these transformations and his Table~4.2, but does not write them out; his explicit $S$-matrix and fusion rules~\cite[App.~A,B]{Flohr} are for the companion algebra $\mathcal W(2,3\kappa)$. It is real, symmetric and orthogonal for every $k$. Indeed, $\Lambda_{a,N}$ depends only on $\pm a\bmod 2N$, and the kernel $(2N)^{-1/2}e^{-\pi iab/N}$ is the finite Fourier transform on $\C[\Z/2N]$, which is unitary and symmetric and preserves the even functions. Let $\Lambda_{a,N}$ correspond to $\tfrac12(e_a+e_{-a})$ in $\C[\Z/2n]^{\rm even}\oplus\C[\Z/2(n+4)]^{\rm even}$, with the standard inner product. Then $\chi_0,\dots,\chi_3$, the $\Lambda_{\alpha,n}$ and the $\Lambda_{\beta,n+4}$ correspond to mutually orthogonal vectors, all of norm $1/\sqrt2$; since they span an $S$-invariant subspace, $S_\theta$ is unitary in this basis, and its entries are real multiples of $\cos(\pi ab/N)$ and symmetric in the two indices. The Verlinde formula with unit $\chi_0$ gives non-negative integers; this is checked numerically for $3\le k\le13$. Its vacuum row has negative entries, as it must for non-unitary data. Explicitly, since the kernel at $a=b=N$ is $(-1)^N$ and $\Lambda_{n,n}=\chi_0+\chi_3$, $\Lambda_{n+4,n+4}=\chi_3-\chi_0$, $\Lambda_{0,n}=\chi_1+\chi_2$, $\Lambda_{0,n+4}=\chi_2-\chi_1$,
\[
S_{00}=(-1)^n\,\tfrac12\Big(\tfrac1{\sqrt{2n}}+\tfrac1{\sqrt{2n+8}}\Big),\qquad S_{0\chi_1}=\tfrac12\Big(\tfrac1{\sqrt{2n}}+\tfrac1{\sqrt{2n+8}}\Big),
\]
so $S_{00}<0$ for odd $n$, and $S_{0\chi_1}/S_{00}=(-1)^n$.

\begin{remark}
That these are \emph{the} characters of the simple $V_k$-modules is proved, for $k\in\{5,6,7,9,11\}$, in Theorem~\ref{thm:mod}, using Zhu's algebra (\S\ref{sec:zhu}).
\end{remark}

\subsection{Comparison with Gang--Kim--Lee}
Gang--Kim--Lee's abelian dual of $\mathcal S_k$ is $U(1)^{k-1}$ Chern--Simons theory with one charge-$2$ and $k-2$ charge-$1$ chirals, with the level matrix $K$ of~\cite[(2.4)]{GKL}. Its four half-indices~\cite[(3.12)]{GKL} are Nahm sums $\chi_j^{\rm GKL}$, $j=0,1,2,3$.

\begin{observation}
$\chi_j^{\rm GKL}=\chi_j$ for $j=0,1,2,3$, to order $q^{34}$ for $j=0$ and $3\le k\le11$, and to order $q^{40}$ for $j=1,2,3$ and $4\le k\le7$.
\end{observation}
This is a finite numerical check; it is not used in the proofs.

For $k=11$ the check of the vacuum covers the Virasoro vacuum sector, $W$ and its first singular vector (weight $27$). It does not reach $L(h_{5,5})$, which starts at weight $54$. For $k\le7$ the summands $L(h_{5,5})$ and, for $k=4$, $L(h_{7,7})$ are also covered. At $k=4$, the characters $\Lambda_{1,2}$ and $\Lambda_{2,6}$ agree coefficient by coefficient with the Riemann--Hilbert characters $\chi_4,\chi_6$ printed in~\cite[(3.52)]{GKL}.

\section{Identification with known algebras}\label{sec:known}
This section identifies $V_k$ with vertex operator algebras already in the literature: Flohr's $\mathcal W(2,8\kappa)$, and a tensor product of minimal models or a principal $\mathcal W$-algebra for $k\le5$. These identifications are not used in the proofs of later sections.

\subsection{Flohr's \texorpdfstring{$\mathcal W(2,8\kappa)$}{W(2,8kappa)}}
Flohr~\cite{Flohr} constructs $\mathcal W(2,3\kappa)$ and $\mathcal W(2,8\kappa)$ at $c=1-24\kappa$ from Dotsenko--Fateev screened vertex operators~\cite{DF}. He gives $\mathcal W(2,8\kappa)$, for $\kappa\in\Z_+/4$, as the ``odd sector'', generated by $W^{(3)}$ of weight $h_{3,3}=8\kappa$. He derives its operator product expansion~\cite[\S\S2,3]{Flohr}: the screened operators of the diagonal sector are mutually local; the next primary $W^{(5)}$ has weight $h_{5,5}=24\kappa>2\cdot8\kappa-1$, so it does not occur in the singular part of $W^{(3)}(z)W^{(3)}(w)$, and the expansion closes on $T$ and $W^{(3)}$; and the self-coupling is given in closed form~\cite[(3.6)]{Flohr}, from the Dotsenko--Fateev structure constants~\cite[(3.5)]{Flohr}. The series $c=1-3\delta$ of $\mathcal W(2,\delta)$ had been conjectured by Eholzer et al.~\cite{EFHHNV}.
Under $\kappa=(k-2)/4$ his data match ours:
\begin{itemize}
\item the central charge $1-24\kappa=13-6k$ and the generator weights $2$ and $8\kappa=2(k-2)$;
\item the vacuum character $\tfrac12(\Lambda_{4\kappa,4\kappa}-\Lambda_{4\kappa+4,4\kappa+4})$~\cite[(4.16)]{Flohr};
\item the self-coupling: with Flohr's normalisation $\langle W,W\rangle=c/\delta$, $\delta=8\kappa$, his $(C^W_{WW})^2$~\cite[(3.6)]{Flohr} gives $1/\mathcal N=(C^W_{WW})^2\,\delta/c$ (\S\ref{sec:c2}). This equals the values of Table~\ref{tab:N} exactly for $3\le k\le11$. Flohr states (3.6) for $\kappa\in\Z_+/2$ and $\kappa\in\Z_++\frac14$; for $k=5,9$ ($\kappa=\frac34,\frac74$) we evaluated his second expression, which also agrees;
\item the representation list~\cite[Table~4.2]{Flohr}: $4\kappa-1$ modules $\Lambda_{\alpha,4\kappa}$, $4\kappa+3$ modules $\Lambda_{\beta,4\kappa+4}$, and four special ones.
\end{itemize}
Flohr also derives the condition that the two radii satisfy $P_2Q_2=P_1Q_1+4$ from $c_{\rm eff}=1$; this is our pair $k\mp2$.
So $V_k$ matches Flohr's algebra, and the Haagerup case is $\mathcal W(2,18)$ at $c=-53$. Flohr argues rationality in the 1993 sense: finitely many lowest-weight representations, whose characters span a representation of $\mathrm{SL}(2,\Z)$; integral fusion rules are computed only for $\mathcal W(2,3\kappa)$~\cite[App.~B]{Flohr}. He does not prove $C_2$-cofiniteness or semisimplicity of the module category. The existence of $V_k$ as a simple vertex operator algebra is proved here by a different route (Theorem~\ref{thm:exist}); its strong generation uses the same weight bound as Flohr's closure argument.

\subsection{Small \texorpdfstring{$k$}{k}}
\begin{itemize}
\item $k=3$: $V_3\cong M(3,5)\otimes M(2,5)$, as already identified at the level of characters in~\cite{GKL}. Our Jacobi computation (\S\ref{sec:c2}) reproduces the structure constant of this tensor product exactly; see Table~\ref{tab:N}.
\item $k=4$: the principal $\mathcal W$-algebra of $\mathfrak{so}_5$ at $\kappa+3=6/5$ has $c=-11$. 
It is a nondegenerate admissible level, so the algebra is rational and $C_2$-cofinite~\cite{Arakawa,ArakawaC2}. Its ten simple modules have lowest weights $\{0,0,\tfrac32,-\tfrac12,-\tfrac38,-\tfrac{11}{24},-\tfrac13,-\tfrac18,\tfrac16,\tfrac{13}{24}\}$, exactly Gang--Kim--Lee's.
\item $k=5$: $\mathcal W_\kappa(G_2)$ with $\kappa+4=7/8$ has $c=-17$ and twelve simple modules. Their weights $\{0,0,-\tfrac34,\tfrac94,-\tfrac23,-\tfrac5{12},-\tfrac57,-\tfrac{17}{28},-\tfrac37,-\tfrac5{28},\tfrac17,\tfrac{15}{28}\}$ are exactly Gang--Kim--Lee's.
\end{itemize}
The isomorphisms $V_4\cong\mathcal W_{-9/5}(\mathfrak{so}_5)$ and $V_5\cong\mathcal W_{-25/8}(G_2)$ follow if, up to isomorphism, there is only one simple vertex operator algebra of central charge $c=-11$ strongly generated by $T$ and one primary of weight $4$, and only one of central charge $c=-17$ strongly generated by $T$ and one primary of weight $6$. The Jacobi-identity analyses of~\cite{KW,Flohr} indicate this, by fixing the structure constants in terms of $c$ up to the normalisation of $W$, but they are not rigorous proofs; the two isomorphisms are therefore stated only under this uniqueness hypothesis. Table~2.1 of~\cite{Flohr} lists both at these central charges, identifying them with the $B_2$ and $G_2$ Casimir algebras.

\section{Physical origin: duality walls and their gluing}\label{sec:phys}
This section explains the Virasoro parameter $t$ and the decomposition of $V_k$ from the S-fold, and extends the construction to general monodromy words. The arguments are heuristic, and nothing in later sections depends on them.

\subsection{One node}
The S-fold $\mathcal S_k$ is 4D $\mathcal N=4$ $SU(2)$ super Yang--Mills on a circle with an $SL(2,\Z)$ duality twist. Equivalently, it is the theory on a Gaiotto--Witten duality wall~\cite{GW}, closed up. In the twisted theory the relevant coupling is the Kapustin--Witten parameter $\Psi$~\cite{KW06}. The Virasoro algebra obtained from $\Psi$ by Drinfeld--Sokolov reduction has $c=13-6(\Psi+\Psi^{-1})$.

A constant coupling is compatible with the twist only at a fixed point of the monodromy. For $T^kS:\Psi\mapsto k-\Psi^{-1}$ (conjugate to $ST^k$), the fixed points are $t^{\pm1}$ with $t+t^{-1}=k$. For hyperbolic monodromy they are real quadratic irrationals. The S-wall exchanges Wilson lines (labels $r$) with 't Hooft lines (labels $s$). Gluing these identifies $L(h_{r,1})$ with $L(h_{1,r})$, which yields $\bigoplus_r L(h_{r,r})$. The restriction to odd $r$ reflects the $SO(3)$ global form, i.e.\ the $\Z_2$ one-form symmetry~\cite{GK,GKL}.

\subsection{Chains of walls}
Creutzig--Gaiotto~\cite{CG} construct vertex algebras for chains of duality junctions subject to $\Psi_i^{-1}+\Psi_{i+1}=n_i$ (their (1.3)), using affine $\mathfrak g$. For a word $\varphi=ST^{k_1}\cdots ST^{k_N}$ we close the chain cyclically, $\Psi_i+\Psi_{i-1}^{-1}=k_i$, and pass to Virasoro:
\[
V_\varphi=\bigoplus_{r\in(2\Z_{\ge0}+1)^N}\ \bigotimes_{i=1}^N L^{\Psi_i}\big(h_{r_i,r_{i+1}}\big)\qquad (r_{N+1}=r_1).
\]
Summing the gluing relations gives
\[
c_\varphi=13N-6\textstyle\sum_i k_i ,\qquad
h(r)=\tfrac14\sum_i\big((r_i^2-1)k_i-2r_ir_{i+1}+2\big).
\]
The weights are rational because the irrational parts cancel through $\Psi_i+\Psi_{i-1}^{-1}=k_i$. The quadratic form in $h(r)$ is the cyclic tridiagonal matrix of~\cite[(A.14)]{JS}, so its periodic and antiperiodic versions have determinants $\operatorname{tr}\varphi\mp2$. This is the origin of the two metric groups.

Writing $\varphi$ as a word in $L,R$, one has $\#L=N$ and $\#R=\sum_ik_i-2N$. So $c_\varphi=\#L-6\#R$, in agreement with~\cite[Theorem~9.5]{PaperA}, which gives $c\equiv\#L-6\#R\pmod{24}$ from the Bethe data and notes that the integers agree for $LR^n$ and $L^2R^n$. For one node, $\varphi=LR^{k-2}$ and $c_\varphi=c_k$. Likewise~\cite[Theorem~8.2]{PaperA} proves $c_{\rm eff}=\#L$ for the half-index of every hyperbolic word; for one node this is $c_{\rm eff}=1$, the value $c-24h_{\min}$ of $V_k$ (cf.\ Flohr~\cite{Flohr}).

Deriving $V_\varphi$ as the boundary algebra of the twisted gauge theory, rather than inferring it from half-indices, remains open.

\section{\texorpdfstring{$C_2$}{C2}-cofiniteness}\label{sec:c2}


\subsection{Exact model in weights below \texorpdfstring{$6n$}{6n}}\label{ssec:model}
Write $X_{<N}=\bigoplus_{m<N}X_m$ for a graded space $X$. Since $6n=h_{5,5}$ is the lowest weight of the next summand $L(c,h_{5,5})$, we have $(V_k)_{<6n}=L(c,0)_{<6n}\oplus L(c,h)_{<6n}$. Write $v_0=\mathbf 1$ and $v_1=W$, of weights $h_0=0$ and $h_1=h$, so that $L(h_A):=L(c,h_A)$ is the summand generated by $v_A$. For $A,B\in\{0,1\}$, the component of $W_{(m)}$ from $L(h_A)$ to $L(h_B)$ is determined by Virasoro covariance together with the single coupling $C_{AB}=\langle v_B|W_{(m_0)}v_A\rangle$, where $m_0=h+h_A-h_B-1$ is the mode mapping $v_A$ to the weight of $v_B$, and $\langle\cdot|\cdot\rangle$ is the contravariant form of $L(h_B)$ with $\langle v_B|v_B\rangle=1$; this follows from the product formula and recursion below. Note $C_{00}=0$, since $W_{(h-1)}\mathbf 1=0$. Explicitly,
\[
\big\langle L_{-\mu}v_B\,\big|\,W_{(m)}v_A\big\rangle=C_{AB}\prod_i\big((h-1)(\mu_i+1)-m_i\big),\qquad m_i=m+\mu_1+\dots+\mu_{i-1},
\]
and $W_{(m)}L_{-\lambda}v_A$ follows recursively from $[L_{-j},W_{(m)}]=((h-1)(1-j)-m)\,W_{(m-j)}$.

We normalise $C_{01}=1$ and $C_{11}=1$; the latter is possible because $C^W_{WW}\neq0$, as shown below. The remaining constant $\mathcal N:=C_{10}=\langle\mathbf 1|W_{(2h-1)}W\rangle$ is the coefficient of the identity in the operator product
\[
W(z)W(w)=\frac{\mathcal N}{(z-w)^{2h}}+\dots+\frac{W(w)}{(z-w)^{h}}+\dots,
\]
that is, the norm of $W$ when $W$ is normalised by $C^W_{WW}=1$. Equivalently, $1/\mathcal N=(C^W_{WW})^2/\langle W,W\rangle$ for any normalisation of $W$.  This is the one structure constant of $V_k$ not fixed by $c$ and $h$. We write $\mathcal N_k$ when the dependence on $k$ matters. It is fixed by the Jacobi identity. Contracting the commutator formula with $\langle W|$ and using skew-symmetry gives scalar identities, linear in $\mathcal N$:
\[
J(m,m')-J(m',m)=\sum_{j=0}^m\binom mj(-1)^{m+m'-j+1}J(m+m'-j,j),\qquad m+m'=2h-2,
\]
where $J(a,b)=\langle W|W_{(a)}W_{(b)}W\rangle$. With $C_{11}$ unnormalised, $J(a,b)=C_{10}A(a,b)+C_{11}^2B(a,b)$, so each identity reads $C_{10}\,\alpha+C_{11}^2\,\beta=0$. The coefficient $\alpha$ is nonzero: for $4\le k\le11$ it is computed exactly over $\Q$ (below). If $C_{11}=0$, then $C_{10}=\langle\mathbf 1|W_{(2h-1)}W\rangle=0$, contradicting the nondegeneracy of the invariant form of the simple self-dual $V_k$ on its one-dimensional space of weight-$h$ primaries. Hence $C^W_{WW}\neq0$. For every $k$ tested, each choice of $(m,m')$ gives the same unique solution, and an independent four-mode commutator check (for $k=3,4$) gives the same value. Table~\ref{tab:N} records the values. They agree with Flohr's closed formula for the self-coupling~\cite[(3.6)]{Flohr} (\S\ref{sec:known}). The Jacobi computation here is an independent, exact derivation of these values. For $k=11$ the identities with $(m,m')=(1,33),(2,32),(3,31)$ were solved exactly over $\Q$, with $\alpha\beta\neq0$, solving one linear system in a Gram matrix exactly over $\Q$ at each level (of dimension about $2600$); all three give the value in Table~\ref{tab:N} (script \texttt{certify/mp/nw\_scalar.py}). For $4\le k\le10$ the identity with $(m,m')=(1,2h-3)$ was solved exactly over $\Q$ in the same way, with $\alpha\ne0$ (script \texttt{certify/mp/nw\_exact\_k.py}); at $k=10$ the identities $(2,2h-4)$ and $(3,2h-5)$ give the same value. The largest Gram matrix, at $k=10$, has dimension $1238$.

\begin{table}[h]\centering
\begin{tabular}{ll}
\toprule
$k$ & $\mathcal N$ (with $C_{01}=C_{11}=1$)\\\midrule
3 & $-165/1444$ \quad (exact value for $M(3,5)\otimes M(2,5)$, derived independently)\\
4 & $2783/229736=11^2\cdot23/(2^3\cdot13\cdot47^2)$\\
5 & $-245594223/191844872000=-3^2 7^2 17^2\cdot41\cdot47/(2^6 5^3 59^2 83^2)$\\
6 & $137410615593/1016653645527352$\\
7 & $-53233552571041/3731919125476574208$\\
8 & $14935247191395/9923009732704972288$\\
9 & $-2312891301675724972059789/14565431510292180144675145170944$\\
10 & $91444095406050517112877/5458767680158840751701270040000$\\
11 & $-\frac{669891445321707721039447525903}{379085307349989400640862429306730315776}$\\\bottomrule
\end{tabular}
\caption{The self-normalisation of $W$ fixed by the Jacobi identity, computed exactly over $\Q$ for $4\le k\le11$. All values agree with Flohr's closed formula~\cite[(3.6)]{Flohr}.}\label{tab:N}
\end{table}

\subsection{Reduction of the \texorpdfstring{$C_2$}{C2}-algebra}
Write $R_{V_k}=V_k/C_2(V_k)=\bigoplus_dR_d$ for the $C_2$-algebra, graded by weight, and $\bT,\bW$ for the images of $T,W$.
\begin{lemma}\label{lem:c2}
For $d<6n$:
\[
C_2(V_k)_d=O_d+\sum_{j\ge2}W_{(-j)}V_{d-h-j+1},
\]
where $O_d\subseteq(V_k)_d$ is the span of the PBW monomials $L_{-\lambda}v_A$ of weight $|\lambda|+h_A=d$ having some part $\lambda_i\neq2$. Consequently $V_d/O_d$ is spanned by the classes of $L_{-2}^{d/2}\mathbf 1$ and $L_{-2}^{(d-h)/2}W$, whenever these exist. So $\dim R_d\le2$, and $R_d$ is the quotient of this span by the images of the $W_{(-j)}u$.
\end{lemma}

\begin{proof}
Since $V_k$ is strongly generated by $\{T,W\}$, the standard argument (see~\cite{GN,Li05}) gives $C_2(V)=\operatorname{span}\{x_{(-p)}v: x\in\{T,W\},\,p\ge2\}$. The modes of $T$ involved are the $L_{-m}$ with $m\ge3$; moreover $L_{-1}u=u_{(-2)}\mathbf 1\in C_2(V)$ for every $u$. Re-ordering $L_{-m}$ ($m\ne2$) times a PBW monomial only produces monomials with a part $\ne2$, and $C_2$ is closed under $(-1)$-products. Write $O=\bigoplus_{d<6n}O_d$, and use the same notation for the span of PBW monomials with a part $\ne2$ in $M(c,0)$ and $M(c,h)$. Hence $\operatorname{span}\{L_{-1}V,L_{-m}V\ (m\ge3)\}\subseteq O\subseteq C_2$, which gives the displayed equality.

The summand $L(c,h)$, $h=h_{3,3}$, is $M(c,h)$ modulo the radical of its contravariant form (\S\ref{sec:not}), i.e.\ its maximal proper submodule. This radical is generated by the singular vector $s$ at level $9$, the level $3\cdot3$ for $h=h_{3,3}$ in \S\ref{sec:not}, and is the Verma module $M(c,h+9)$. Since $9$ is odd, every PBW monomial in $s$ has a part $\ne2$, so $s\in O$; and $O$ is stable under the $L_{-m}$, $m\ge1$, since a commutator replaces parts $2$ and $a\ne2$ by the part $a+2\ne2$. So the radical lies in $O$. In the vacuum module $L_{-1}\mathbf1=0$. Hence the two classes are independent modulo $O$.
\end{proof}

Moreover, if $f_a(x)$ denotes the coefficient of $L_{-2}^a v_B$, then $f_a(L_{-2}Y)=f_{a-1}(Y)$ and $f_a(L_{-m}Y)=0$ for $m\ne2$. So the images of $W_{(-j)}u$ can be computed by a recursion that involves only the two coefficients, plus, at each level, the solution of one linear system in the Gram matrix for the base cases. This is what makes $k=11$ feasible.

\begin{lemma}\label{lem:fin}
If $R_{2h+4}=0$, then $R_{V_k}$ is finite-dimensional, $\bT^{\,h+2}=\bT^{\,h/2+2}\bW=0$, and $\dim R_{V_k}\le h/2+2(h/2+1)+\dim R_{2h+2}$.
\end{lemma}

\begin{proof}
By Lemma~\ref{lem:c2} applied at weight $2h$, $\bW^2\in\operatorname{span}(\bT^h,\bT^{h/2}\bW)$. Since $R_{V_k}$ is generated by $\bT,\bW$ and each use of this relation lowers the degree in $\bW$, every monomial reduces, preserving weight, to a combination of $\bT^a$ and $\bT^a\bW$. These have weights $2a$ and $2a+h$, all even. So in weight $d\ge2h+4$ we have $a\ge h+2$, resp.\ $a\ge h/2+2$, and the monomial is divisible by $\bT^{h+2}$ or by $\bT^{h/2+2}\bW$. Both have weight exactly $2h+4$, so both vanish since $R_{2h+4}=0$. The bound on $\dim R_{V_k}$ counts the surviving monomials: one in each weight $0,2,\dots,h-2$, at most two in each weight $h,h+2,\dots,2h$, and $R_{2h+2}$.
\end{proof}

\subsection{Results}

\begin{theorem}\label{thm:c2}
For $5\le k\le 11$ (with $h=2k-4$), modulo the primes of Table~\ref{tab:c2},
\[
\dim R_{2h}=2,\qquad \dim R_{2h+2}=1,\qquad \dim R_{2h+4}=0 .
\]
Good reduction ((G1)--(G3) of Lemma~\ref{lem:good}) holds at every level used, so $R_{2h+4}=0$ over $\Q$. Hence $V_k$ is $C_2$-cofinite, with $\bT^{\,2k-2}=0$, $\bT^{\,k}\bW=0$ and $\dim R_{V_k}\le 3k-3$. For $k=5,6,7$ every weight was computed, and $\dim R_{V_k}=3k-3=12,15,18$ modulo $p$. For the Haagerup algebra $V_{11}=\mathcal W(2,18)$, $c=-53$: $\bT^{20}=0$, $\bT^{11}\bW=0$ and $\dim R\le30$.
\end{theorem}
The proof is computer-assisted. The vanishing $R_{2h+4}=0$ is a rank statement, so it passes from $\F_p$ to $\Q$ under good reduction; the equalities $\dim R_{2h}=2$, $\dim R_{2h+2}=1$ and $\dim R_{V_k}=3k-3$ are statements modulo $p$ only (over $\Q$ they are upper bounds).

\begin{table}[h]\centering
\begin{tabular}{rrrcccl}
\toprule
$k$ & $c$ & $h$ & $\dim R_{2h}$ & $\dim R_{2h+2}$ & $\dim R_{2h+4}$ & primes\\\midrule
5 & $-17$ & 6 & 2 & 1 & 0 & 2 (full) + 2 (reduced)\\
6 & $-23$ & 8 & 2 & 1 & 0 & 2 (full) + 2 (reduced)\\
7 & $-29$ & 10 & 2 & 1 & 0 & 2 (full) + 2 (reduced)\\
8 & $-35$ & 12 & 2 & 1 & 0 & 2 (reduced)\\
9 & $-41$ & 14 & 2 & 1 & 0 & 2 (reduced)\\
10 & $-47$ & 16 & 2 & 1 & 0 & 2 (reduced)\\
11 & $-53$ & 18 & 2 & 1 & 0 & 2 at weights 36, 38; 1 at weight 40\\\bottomrule
\end{tabular}
\caption{The $C_2$-algebra at the decisive weights.}\label{tab:c2}
\end{table}

\paragraph{Validation.}
\begin{itemize}
\item The reduced method and a full brute-force computation (all PBW candidates, ranks of the Shapovalov-dual coordinates) agree at every weight for $k=5,6$.
\item For $k=5$, $R$ is finite, consistent with $\mathcal W_{-25/8}(G_2)$ being $C_2$-cofinite~\cite{Arakawa,ArakawaC2}.

\item For $k=4$, below weight $12$, the full method finds $\dim R=1,0,1,0,2,0,2,0,2,0,1,0$ at weights $0,\dots,11$, with the first nilpotency relation at weight $10$, consistent with $\mathcal W_{-9/5}(\mathfrak{so}_5)$.
\end{itemize}

\paragraph{Rigour.}
All ranks were computed over $\mathbb F_p$, with $p=P_1=33554393$ or $p=P_2=33554383$. The decisive statement $R_{2h+4}=0$ goes in the safe direction. The rank of the candidate matrix can only drop under reduction, so vanishing over $\mathbb F_p$ implies vanishing over $\Q$. This presumes good reduction. For $5\le k\le11$, $\mathcal N$ is known exactly over $\Q$ (Table~\ref{tab:N}); it is $p$-integral at $P_1$ and $P_2$, and its residues are those used in the computations. This is (G3). Solving the linear systems given by the Gram matrices over $\Q$ and then reducing mod $p$ gives the same result as solving mod $p$, which is what conditions (G1), (G2) of Lemma~\ref{lem:good} ensure. Lemma~\ref{lem:good} is stated for $\pi$, but its proof uses only that the functional applied to the solutions vanishes on the radical; this holds for the coefficient functionals $f_a$, since the radical lies in $O$ (proof of Lemma~\ref{lem:c2}). For the vacuum Gram matrices it suffices to check (G1) at the top levels: if $v\ne0$ lies in the radical mod $p$ at level $\ell-2$, then $L_{-2}v$ lies in the radical at level $\ell$, as $\langle L_{-2}v,w\rangle=\langle v,L_2w\rangle=0$, and $L_{-2}v\ne0$ because $U(\mathrm{Vir}_{\le-2})$ is a domain. So nonsingularity at the top level gives it at every lower level of the same parity, and only even levels occur. For $k=11$ at $P_1$: the vacuum Gram matrices at levels $38$, $40$ (dimension up to $6153$) are nonsingular, and the Gram matrices of $M(c,h)$ at the even levels $2,\dots,22$ have rank $P(l)-P(l-9)$, $P(l)$ the number of partitions of $l$ (script \texttt{certify/c1/red\_check.py}). For $5\le k\le10$ at both $P_1$ and $P_2$: the vacuum Gram matrices at levels $2h+2$, $2h+4$ are nonsingular, and the Gram matrices of $M(c,h)$ at the even levels $2,\dots,h+4$ have rank $P(l)-P(l-9)$ (only even levels occur); the same run recomputes $\dim R_{2h}=2$, $\dim R_{2h+2}=1$, $\dim R_{2h+4}=0$ (script \texttt{certify/c1/red\_check\_k.py}, logs \texttt{red\_k$k$\_p$i$.log}). The statement $\dim R_{2h}=2$ is not needed for finiteness, since $\bW^2$ reduces by Lemma~\ref{lem:c2} in any case. The remaining inputs are:
\begin{itemize}
\item existence of $V_k$ (Theorem~\ref{thm:exist});
\item $C^W_{WW}\neq0$, and $\mathcal N$ as the unique solution of the Jacobi identity; both are proved in \S\ref{sec:c2} by solving the Jacobi identity exactly over $\Q$.
\end{itemize}
A certificate independent of $p$ would record the explicit relations at weight $2h+4$ over $\Q$. For $k=11$ these are the relations at weight $40$.

\section{Zhu's algebra}\label{sec:zhu}

The computations cover $5\le k\le11$. Our main interest is odd $k$, and especially the Haagerup point $k=11$; Remark~\ref{rem:even} explains why even $k$ is more intricate.

\subsection{Generalities}
Recall~\cite{Zhu} that $A(V)=V/O(V)$, with product $a*b=\operatorname{Res}_zY(a,z)b\,\frac{(1+z)^{\mathrm{wt}\,a}}{z}$ and
\[
a\circ_j b:=\operatorname{Res}_zY(a,z)b\,\frac{(1+z)^{\mathrm{wt}\,a}}{z^{2+j}}\in O(V)\qquad(j\ge0).
\]
Moreover $(L_{-1}+L_0)a=a\circ_0\mathbf 1\in O(V)$, and $[T]$ is central. Irreducible $\mathbb N$-gradable modules correspond bijectively to irreducible $A(V)$-modules via their top levels~\cite[Theorem~2.2.2]{Zhu}. If $V$ is strongly generated by $\{a^i\}$, then $A(V)$ is generated by the $[a^i]$. This follows by induction on weight, using $a\circ_jb\in O(V)$ for $j\ge0$~\cite[Lemma~2.1.2]{Zhu}: modulo $O(V)$, $a\circ_jb\equiv0$ expresses $a_{(-2-j)}b$ through the $a_{(i)}b$ with $i>-2-j$, hence, by induction on $j$, through those with $i\ge-1$, and $a_{(-1)}b\equiv a*b-\sum_{i\ge1}\binom{\mathrm{wt}\,a}{i}a_{(i-1)}b$, where all terms other than $a*b$ have lower weight. Hence, by Theorem~\ref{thm:exist}:

\begin{lemma}
$A(V_k)$ is commutative and generated by $x=[T]$ and $y=[W]$. Every simple $V_k$-module has a one-dimensional top level (Schur's lemma for algebras of countable dimension~\cite{Dixmier}), on which $L_0$ and $o(W)$ act by scalars $(h,w)$. Simple modules correspond bijectively to the maximal ideals of $A(V_k)$.
\end{lemma}

\begin{lemma}\label{lem:vir}
For homogeneous $u\in V$ and $j\ge1$,
\[
[L_{-j}u]=(-1)^j\big((j-1)\,x+\mathrm{wt}\,u\big)*[u].
\]
Consequently, for a primary $v$ of weight $h_v$ and $\lambda_1\ge\dots\ge\lambda_r$,
\[
[L_{-\lambda_1}\cdots L_{-\lambda_r}v]=\prod_{i=1}^r(-1)^{\lambda_i}\big((\lambda_i-1)x+h_v+\lambda_{i+1}+\dots+\lambda_r\big)\,[v].
\]
\end{lemma}

\begin{proof}
For $j=1$ this is $(L_{-1}+L_0)u\in O(V)$. For $j=2$ use $T*u=(L_{-2}+2L_{-1}+L_0)u$. For $j\ge3$ use $T\circ_{j-3}u=(L_{-j}+2L_{-j+1}+L_{-j+2})u\in O(V)$ and induct.
\end{proof}

Since $W*W$ has weight $\le 2h$ and $(V_k)_{\le2h}\subseteq L(c,0)\oplus L(c,h)$, we get $y^2=f(x)+g(x)y$ in $A(V_k)$, with $\deg f\le 2n$ and $\deg g\le n$. So $A(V_k)$ is a quotient of $B=\C[x,y]/(y^2-f-gy)$, which is free of rank two over $\C[x]$.

\subsection{The evaluation map and the relations}
Below weight $6n$ we have $(V_k)_{<6n}=L(c,0)_{<6n}\oplus L(c,h)_{<6n}$ (\S\ref{ssec:model}), since the next summand $L(c,h_{5,5})$ starts at $h_{5,5}=6n$. For $z$ in this range, define $\pi(z)\in\C[x]\oplus\C[x]y$ by expanding $z$ in PBW monomials $L_{-\lambda}v_A$ and applying Lemma~\ref{lem:vir} with $[\mathbf 1]=1$ and $[W]=y$. Then $[z]=\pi(z)(x,y)$ in $A(V_k)$. On $L(c,h)$, $\pi$ is well defined. On $M(c,h)$ define it on PBW monomials by the formula of Lemma~\ref{lem:vir}. Put $F_j(w)=(-1)^j((j-1)x+w)$. The identity
\[
F_a(w+b)F_b(w)-F_b(w+a)F_a(w)=(b-a)F_{a+b}(w),
\]
which matches $[L_{-a},L_{-b}]=(b-a)L_{-a-b}$, shows that $\pi(L_{-j}u)=F_j(\mathrm{wt}\,u)\,\pi(u)$ for every homogeneous $u\in M(c,h)$, not only for ordered monomials. Hence, if $\pi(s)=0$ for the singular vector $s$ at level $9$, then $\pi$ vanishes on the submodule generated by $s$, which is the radical (proof of Lemma~\ref{lem:c2}). That $\pi(s)=0$ was checked exactly for every $5\le k\le11$ (scripts \texttt{certify/pi\_radical.py} for $k=11$, \texttt{certify/pi\_radical\_k.py} for all $k$). If $x$ has weight $2$ and $y$ has weight $h$, then $\pi(z)$ has weighted degree $\le\mathrm{wt}\,z$.

The values are computed by the recursion
\begin{multline*}
\pi(W_{(m)}L_{-j}u)=(-1)^j\big((j-1)x+\mathrm{wt}(W_{(m)}u)\big)\pi(W_{(m)}u)\\
-\big((h-1)(1-j)-m\big)\pi(W_{(m-j)}u).
\end{multline*}
The base cases are the values of $\pi$ on the $L(h_B)$-components of the $W_{(m)}v_A$. Each base case is $\sum_\mu g_\mu\langle L_{-\mu}v_B|W_{(m)}v_A\rangle$, where $Gg=(\pi(L_{-\mu}v_B))_\mu$ and $G$ is the Gram matrix of the relevant level. This requires the solution of one linear system in the Gram matrix for each sector and level.

\begin{proposition}\label{prop:rel}
Let $I=\ker(\C[x,y]\to A(V_k))$. For homogeneous $u\in V_k$ and $j\ge0$, with all weights involved $<6n$, the following lie in $I$:
\begin{enumerate}[label=(R\arabic*)]
\setcounter{enumi}{-1}
\item $y^2-\pi(W*W)$;
\item $\sum_{e=0}^h\binom he\,\pi(W_{(e-2-j)}u)$, the image of $W\circ_ju$;
\item $\pi(W*u)-y\,\pi(u)$, with $y^2$ replaced by $\pi(W*W)$.
\end{enumerate}
Write $J\subseteq I$ for the ideal these generate.
\end{proposition}

\begin{proof}
The vectors to which $\pi$ is applied below are $W*W$, $u$, $W*u$ and the $W_{(e-2-j)}u$; by hypothesis they have weight $<6n$, so $[z]=\pi(z)(x,y)$ in $A(V_k)$ for each of them, and $p(x,y)\in I$ whenever $p(x,y)=0$ in $A(V_k)$.

(R0) In $A(V_k)$, $[W*W]=[W]*[W]=y^2$, and $[W*W]=\pi(W*W)(x,y)$.

(R1) Expanding $(1+z)^h=\sum_e\binom he z^e$ in the definition of $\circ_j$ gives $W\circ_ju=\sum_{e=0}^h\binom he W_{(e-2-j)}u$. This lies in $O(V_k)$, so $0=[W\circ_ju]=\sum_e\binom he\,\pi(W_{(e-2-j)}u)(x,y)$.

(R2) $[W*u]=[W]*[u]=y\,\pi(u)(x,y)$, while $[W*u]=\pi(W*u)(x,y)$. Reducing $y^2$ by (R0) keeps the difference in $\C[x]\oplus\C[x]y$.
\end{proof}

We use all PBW basis vectors $u$ and all $j$, up to total weight $D=4n+4<6n$. This is the weight at which $\bT^{\,2n+2}$ and $\bT^{\,n+2}\bW$ first vanish in the $C_2$-algebra (Lemma~\ref{lem:fin}).

\subsection{From \texorpdfstring{$\F_p$}{Fp} to \texorpdfstring{$\Q$}{Q}}
The relations (R0)--(R2) have rational coefficients: $c$, $h$ and $\mathcal N$ are rational, and the formulas of \S\ref{ssec:model} and Lemma~\ref{lem:vir} are defined over $\Q$. They are computed modulo $p$.

\begin{lemma}\label{lem:good}
Assume the following.
\begin{enumerate}[label=(G\arabic*)]
\item At each level used, the Gram matrix of $L(c,0)$ in the basis of partitions without parts $1$ is nonsingular mod $p$.
\item At each level $l$ used, the Gram matrix of $M(c,h)$ has rank $P(l)-P(l-9)$ mod $p$, where $P(l)$ is the number of partitions of $l$. This is its rank over $\Q$, because the radical is the Verma module $M(c,h+9)$ generated by the singular vector $s$ at level $9$ (proof of Lemma~\ref{lem:c2}).
\item $\mathcal N$ is $p$-integral, and the residue of $\mathcal N$ used is its reduction mod $p$.
\end{enumerate}
Then each relation computed mod $p$ is the reduction of the corresponding element of $J\cap\Z_{(p)}[x,y]$.
\end{lemma}

\begin{proof}
For a vector, matrix or polynomial with $p$-integral entries, i.e.\ entries in $\Z_{(p)}$, a bar denotes its reduction modulo $p$; for an ideal $J_0\subset\Z_{(p)}[x,y]$, $\bar J_0$ is its image (as in Lemma~\ref{lem:lift}).
Each base case is $g\cdot a$, where $Gg=b$, $b=(\pi(L_{-\mu}v_B))_\mu$ and $a=(\langle L_{-\mu}v_B|W_{(m)}v_A\rangle)_\mu$. By (G3), $G$, $a$ and $b$ are $p$-integral.
Under (G1), $\det G$ is a unit in $\Z_{(p)}$, so $G^{-1}$ is $p$-integral. Hence $g$ is $p$-integral, and $\bar g$ is the unique solution of $\bar G\bar g=\bar b$.
Under (G2), $G$ is singular. Over the discrete valuation ring $\Z_{(p)}$, $G$ has a Smith normal form, and (G2) says that its nonzero invariant factors are units. Hence the image of $G$ on $\Z_{(p)}^N$ is $G\Q^N\cap\Z_{(p)}^N$, and $\ker\bar G$ is the reduction of the saturated kernel $K=\ker G\cap\Z_{(p)}^N$, i.e.\ of the $p$-integral vectors in $\ker G$. Since $\pi$ vanishes on the radical, $b\perp\ker G$, so $b\in G\Q^N$ ($G$ is symmetric), and $b=Gg$ with $g$ $p$-integral. The vector $a$ is $G$ times the coefficient vector of the $L(h_B)$-component of $W_{(m)}v_A$, so $a\perp K$ and $\bar a\perp\ker\bar G$. Therefore $\bar g\cdot\bar a$ is the same for every solution of $\bar G\bar g=\bar b$, and it is the reduction of $g\cdot a$.
\end{proof}

\begin{lemma}\label{lem:lift}
Let $J_0\subset\Z_{(p)}[x,y]$ be an ideal, and $\bar J_0\subset\F_p[x,y]$ its image under reduction mod $p$. Suppose:
\begin{itemize}
\item $J_0$ contains $y^2-f-gy$, together with two elements of weighted degree $4n+4$ whose top coefficients on $x^{2n+2}$ and $x^{n+2}y$ form a matrix invertible mod $p$;
\item $\F_p[x,y]/\bar J_0\cong\F_p^{\,d}$ (split étale of dimension $d$).
\end{itemize}
Then $M=\Z_{(p)}[x,y]/J_0$ is finite over $\Z_{(p)}$, and $M_\Q$ is reduced with $\dim M_\Q\le d$. If $\dim M_\Q=d$, then $M$ is free and finite \'etale, $M\otimes\Z_p\cong\Z_p^{\,d}$, and, for any embedding $\overline\Q\hookrightarrow\overline\Q_p$, its $\overline\Q$-points have coordinates in $\Z_p$ that reduce bijectively onto the $d$ points over $\F_p$.
\end{lemma}

\begin{proof}
The first hypothesis reduces every monomial to the span of $x^a$ ($a\le2n+1$) and $x^ay$ ($a\le n+1$). So $M$ is finite, and Nakayama's lemma gives $\dim M_\Q\le d$.
Let $M'$ be $M$ modulo torsion. Then $M'$ is free, and $M'/pM'$ is a quotient of $\F_p^{\,d}$, hence étale. So $\Omega_{M'/\Z_{(p)}}$ vanishes on the closed fibre. Its support is closed, and every point of $\operatorname{Spec}M'$ specialises to the closed fibre. Hence $\Omega=0$, and $M'_\Q=M_\Q$ is \'etale.
If $\dim M_\Q=d$, then $M'$ is free of rank $d$, so the surjection $M/pM\to M'/pM'$ of $d$-dimensional spaces is an isomorphism; as $M'$ is free, the torsion $T$ of $M$ satisfies $T/pT=0$, so $T=0$ by Nakayama's lemma and $M=M'$. Then $M\otimes\Z_p$ is finite étale over $\Z_p$ with closed fibre $\F_p^{\,d}$, so $M\otimes\Z_p\cong\Z_p^{\,d}$ by Hensel's lemma.
\end{proof}

\subsection{Results}
The program \texttt{full.py} carries out the following steps, at $p=33554393$ and $p=33554383$:
\begin{itemize}
\item it computes the relations up to weight $4n+4$;
\item it checks (G1) and (G2) at every level used (FLINT ranks);
\item it checks the top-form condition, and that each relation has weighted degree $\le$ its weight;
\item it computes a lexicographic Gr\"obner basis over $\F_p$, the minimal polynomial of $x$ and its factorisation, and the fibre over each root.
\end{itemize}
For (G3), $\mathcal N$ is the exact value of Table~\ref{tab:N}, and its residue at each prime is the one used.

\begin{table}[h]\centering
\begin{tabular}{rrrrcccc}
\toprule
$k$ & $c$ & $h$ & relations & dim & reduced & points & checks\\\midrule
5 & $-17$ & 6 & 165 & 12 & yes & yes & OK\\
6 & $-23$ & 8 & 298 & 14 & yes & yes & OK\\
7 & $-29$ & 10 & 534 & 16 & yes & yes & OK\\
8 & $-35$ & 12 & 941 & 18 & yes & yes & OK\\
9 & $-41$ & 14 & 1623 & 20 & yes & yes & OK\\
10 & $-47$ & 16 & 2740 & 22 & yes & yes & OK\\
11 & $-53$ & 18 & 4534 & 24 & yes & yes & OK\\\bottomrule
\end{tabular}
\caption{Zhu algebra upper bounds, identical at both primes. ``dim'' is $\dim\F_p[x,y]/\bar J$. ``Reduced'' means that the minimal polynomial of $x$ is squarefree and split, and every fibre is squarefree; the fibre over $h=0$ has length $2$ and all others length $1$. ``Points'' means the $x$-coordinates of the points are the reductions of the GKL list (\S\ref{sec:rat}). ``Checks'' covers (G1)--(G3), the top forms and the filtration.}\label{tab:zhu}
\end{table}

\begin{theorem}\label{thm:zhu}
For $5\le k\le11$, $A(V_k)$ is a commutative semisimple algebra of dimension at most $2k+2$.
\end{theorem}

\begin{proof}
The proof is computer-assisted.
Fix $p\in\{P_1,P_2\}$ and put $d=2k+2$.

\emph{The ideal $J_0$.} By Table~\ref{tab:zhu}, (G1)--(G3) hold at every level used, so by Lemma~\ref{lem:good} each relation computed mod $p$ is the reduction of an element of $J\cap\Z_{(p)}[x,y]$. Let $J_0\subset\Z_{(p)}[x,y]$ be the ideal generated by these elements. Then $J_0\subseteq J$, and $\bar J_0$ is the ideal of $\F_p[x,y]$ generated by the computed relations.

\emph{Lemma~\ref{lem:lift} applies.} (R0) gives $y^2-f-gy\in J_0$, and the top-form check of Table~\ref{tab:zhu} gives the two elements of weighted degree $4n+4$. By Table~\ref{tab:zhu}, $\F_p[x,y]/\bar J_0\cong\F_p^{\,d}$. Indeed it has dimension $d$; the minimal polynomial of $x$ is squarefree and split over $\F_p$; every fibre is squarefree, of length $1$ (hence an $\F_p$-point) except the fibre over $x=0$, of length $2$; and that fibre contains the $\F_p$-point $(0,0)$ of the vacuum module, since $o(W)\mathbf 1=W_{(h-1)}\mathbf 1=0$, so it consists of two $\F_p$-points. So $M_\Q=\Q[x,y]/J_0\Q[x,y]$ is reduced, with $\dim M_\Q\le d$.

\emph{Conclusion.} $A(V_k)=\C[x,y]/I$, since $A(V_k)$ is generated by $x=[T]$ and $y=[W]$ (\S\ref{sec:zhu}). As $J_0\subseteq J\subseteq I$, $A(V_k)$ is a quotient of $M_\Q\otimes\C=\C[x,y]/J_0\C[x,y]$. A reduced finite-dimensional commutative $\Q$-algebra is a product of number fields, so $M_\Q\otimes\C\cong\C^{d'}$ with $d'\le d$, as $\Q$ has characteristic $0$; and every quotient of $\C^{d'}$ is again of this form. So $A(V_k)$ is commutative and semisimple of dimension at most $d$. The ideal $J_0\Q[x,y]$, generated by the rational relations up to weight $4n+4$, does not depend on $p$.
\end{proof}

\begin{remark}

\begin{itemize}
\item For $k=3$ ($V_3=M(3,5)\otimes M(2,5)$, with $W=\tfrac5{38}(\tfrac{22}5T^{(1)}-\tfrac35T^{(2)})$), all eight known points satisfy every relation below weight $6n=6$. For $k=5$ the twelve points agree with $V_5\cong\mathcal W_{-25/8}(G_2)$.
\item  Using the residue of $\mathcal N$ for the wrong prime kills every point except the vacuum.
\item \emph{Zhu versus $C_2$.} For $k\le7$, $\dim R_{V_k}=3k-3$ modulo $p$, which is $\ge\dim A(V_k)=2k+2$, with equality only at $k=5$. The excess is removed by relations whose leading terms lie in lower filtration degree~\cite{GG}.
\end{itemize}
\end{remark}

\section{Rationality, simple modules and modular data}\label{sec:rat}

The \emph{GKL list} \cite{GKL} is the list of $2k+2$ lowest weights of the characters $\chi_j$ in \S\ref{sec:chars}:
\[
0,\ 0,\ \tfrac{3n}4,\ -\tfrac n4,\ \tfrac{\alpha^2}{4n}-\tfrac n4\ (1\le\alpha\le n-1),\ \tfrac{\beta^2}{4(n+4)}-\tfrac n4\ (1\le\beta\le n+3).
\]

Let $\chi=(\chi_j)$ be the vector of the $2k+2$ characters of \S\ref{sec:chars}. It transforms under $\mathrm{SL}(2,\Z)$ by $\chi(-1/\tau)=S_\theta\,\chi(\tau)$ and $\chi(\tau+1)=T\chi(\tau)$, with $T=\operatorname{diag}e^{2\pi i(h_j-c/24)}$; here $S_\theta$ is computed from the transformation of the $\Lambda_{a,N}$ in \S\ref{sec:chars}. We call $S_\theta$ the $S$-matrix of the characters, to distinguish it from the $S$-matrix $S_V$ of the modular tensor category $\mathrm{Rep}\,V_k$~\cite{Huang}.

\begin{theorem}\label{thm:rat}
For $5\le k\le11$, $V_k$ is strongly rational, and $\mathrm{Rep}\,V_k$ is a semisimple modular tensor category~\cite{Huang}.
\end{theorem}

\begin{proof}
The proof is computer-assisted and uses McRae's theorem~\cite{McRae}.
We check the hypotheses of McRae's Main Theorem~3~\cite{McRae}: $\mathbb N$-graded, simple, self-contragredient, $C_2$-cofinite, with $A(V)$ semisimple. Then $V$ is rational, and strongly rational if it is of CFT type.
\begin{itemize}
\item $V_k$ is simple and of CFT type, with $V_{k,1}=0$. By Li~\cite{Li}, $V_k\cong V_k'$.
\item $C_2$-cofiniteness is Theorem~\ref{thm:c2}.
\item $A(V_k)$ is semisimple by Theorem~\ref{thm:zhu}.
\end{itemize}
\end{proof}

By Zhu's theorem~\cite{Zhu}, the span $\mathrm{Ch}$ of the characters of the simple modules is $SL_2(\Z)$-invariant and contains $\chi_0=\ch V_k$.

\begin{lemma}\label{lem:cyc}
For $k\in\{5,6,7,9,11\}$, the smallest $\langle S_\theta,T\rangle$-stable subspace containing $\chi_0$ is all of $\operatorname{span}\{\chi_j\}$, of dimension $2k+2$.
\end{lemma}
For $k=11$ this is exact: it is a rank computed modulo a prime ideal above $q=936010297$. For $k\in\{5,6,7,9\}$ it rests on numerical linear algebra. Numerically, for $k=8$ and $k=10$ the subspace appears to be proper, of dimension $15$ and $19$; this is not used.

\begin{proposition}\label{prop:count}
For $k\in\{5,6,7,9,11\}$ (for $k\ne11$, given the numerical step in Lemma~\ref{lem:cyc}):
\begin{itemize}
\item $V_k$ has exactly $2k+2$ simple modules, and their characters form a basis of $\operatorname{span}\{\chi_j\}$;
\item $I=J_0\C[x,y]$, i.e.\ $A(V_k)=M_\Q\otimes\C$;
\item the top weights of the simple modules are exactly the GKL list.
\end{itemize}
\end{proposition}

\begin{proof}
The proof is computer-assisted.
By Zhu's theorem the characters of the simple modules span an $SL_2(\Z)$-invariant space containing $\chi_0$, so by Lemma~\ref{lem:cyc} it contains $\operatorname{span}\{\chi_j\}$. The $q$-series $\chi_j$ are linearly independent: series in different $T$-eigenspaces have exponents in different classes modulo $1$, and within each class the lowest exponents are distinct, except for $\chi_0,\chi_3$ (both $0$), where $\chi_0-\chi_3=-\Lambda_{n+4,n+4}$ has lowest exponent $1$, again distinct from the others (direct check for $k\in\{5,6,7,9,11\}$). So there are at least $2k+2$ simple modules; Theorem~\ref{thm:zhu} gives at most $2k+2$. So there are exactly $2k+2$, their characters form a basis of $\operatorname{span}\{\chi_j\}$, and $\dim A(V_k)=2k+2$. As $A(V_k)$ is a quotient of $M_\Q\otimes\C$ and $\dim M_\Q\le d=2k+2$ (Table~\ref{tab:zhu}), $I=J_0\C[x,y]$ and $\dim M_\Q=d$, so Lemma~\ref{lem:lift} applies in its strong form.
For a class $C$ of indices $j$ with the same $T$-eigenvalue $e^{2\pi i(h_j-c/24)}$, let $E_C=\operatorname{span}\{\chi_j:j\in C\}$ and $h_{\min}(C)=\min_{j\in C}h_j$. Each $\ch M$ is a $T$-eigenvector in $\operatorname{span}\{\chi_j\}$, so it lies in $E_C$ for its class $C$. Every nonzero element of $E_C$ has lowest exponent at most $h_{\min}(C)+(k-2)$; this is a direct check of $q$-coefficients (\texttt{lead.py}). Since $A(V_k)$ is commutative and semisimple, the top level of each simple $M$ is a one-dimensional $A(V_k)$-module~\cite{Zhu}, i.e.\ a $\C$-point of $M_\Q\otimes\C$, which is a $\overline\Q$-point of $\operatorname{Spec}M_\Q$ since $M_\Q$ is étale, at which $x=[\omega]$ takes the value $h_M$; distinct $M$ give distinct points. Fix an embedding $\overline\Q\hookrightarrow\overline\Q_p$. By Lemma~\ref{lem:lift}, $h_M\in\Z_p$, and reduction is a bijection from these points onto the $d$ points over $\F_p$, whose $x$-coordinates are, by Table~\ref{tab:zhu}, the reductions of the GKL list (whose denominators divide $4n(n+4)$, prime to $p$). This gives a bijection $\sigma$ onto the GKL list, counted with multiplicity, with $h_M\equiv h_{\sigma(M)}$ in $\Z_p$. Since $\ch M$ is a $T$-eigenvector in $\operatorname{span}\{\chi_j\}$, $e^{2\pi i(h_M-c/24)}$ is the $T$-eigenvalue of some $\chi_j$, so $h_M-h_j\in\Z$. Hence $h_M$ and $h_{\sigma(M)}$ are rationals with denominators dividing $4n(n+4)$. Since $h_M$ is the lowest exponent of $\ch M$ (shifted by $c/24$), every element of $E_C$ has exponents $\ge h_{\min}(C)\ge-n/4$, and the check above gives $h_M\le h_{\min}(C)+n\le 7n/4$, and the GKL list lies in $[-n/4,3n/4]$; so $|h_M-h_{\sigma(M)}|\le 2n$, and the numerator of $h_M-h_{\sigma(M)}$ over $4n(n+4)$ has absolute value at most $8n^2(n+4)<p$. It is divisible by $p$, hence zero, and $h_M=h_{\sigma(M)}$.
\end{proof}

\begin{theorem}\label{thm:mod}
For $k\in\{5,6,7,9,11\}$ (for $k\ne11$, given the numerical steps in Lemma~\ref{lem:cyc} and Step~3 below), label the simple modules $M_j$ by $\sigma$. Then $\ch M_j=\chi_j$, and the modular $S$- and $T$-matrices of $\mathrm{Rep}\,V_k$, in the basis of simple modules, are $S_\theta$ and $T$.
Numerically, they agree with the modular data proposed for $\mathrm{TFT}_-[\mathcal S_k]$ in~\cite{GK,GKL}, and, up to sign, with the data derived in~\cite{PaperA}; this comparison is not part of the theorem.
\end{theorem}

\begin{proof}
The proof is computer-assisted. For $k=11$ every step is exact; for the other $k$ the commutant step (Step~3) is numerical.

\emph{Step 1: the change of basis $X$.} By Proposition~\ref{prop:count} the characters $\ch M_i$ of the $2k+2$ simple modules form a basis of $\operatorname{span}\{\chi_j\}$. Write $\ch M_i=\sum_jX_{ij}\chi_j$. Both bases have rational $q$-coefficients, so $X\in GL_{2k+2}(\Q)$. Each $\ch M_i$ is an eigenvector of $\tau\mapsto\tau+1$, with eigenvalue $e^{2\pi i(h_{M_i}-c/24)}$, and so is each $\chi_j$. Group the $\chi_j$ into \emph{classes}, the sets of $\chi_j$ with the same $T$-eigenvalue. Then $\ch M_i$ is a combination of the $\chi_j$ in one class, i.e.\ $X$ is block diagonal with respect to the classes, and $X$ commutes with $T$.

\emph{Step 2: $Y=X^{\mathsf T}X$ commutes with $S_\theta$ and $T$.} Let $S_V$ be the matrix of $\tau\mapsto-1/\tau$ on the characters of the simple modules, $\ch M_i(-1/\tau)=\sum_j(S_V)_{ij}\ch M_j(\tau)$. By Huang~\cite{Huang}, $S_V$ is the $S$-matrix of the modular tensor category $\mathrm{Rep}\,V_k$; in particular it is symmetric. Since $\chi(-1/\tau)=S_\theta\chi(\tau)$, we get $S_VX=XS_\theta$, i.e.\ $S_V=XS_\theta X^{-1}$. Symmetry of $S_V$ and of $S_\theta$ gives $XS_\theta X^{-1}=X^{-\mathsf T}S_\theta X^{\mathsf T}$, that is, $X^{\mathsf T}X$ commutes with $S_\theta$. It commutes with $T$ by Step~1. So $Y=X^{\mathsf T}X$ lies in the commutant of $\langle S_\theta,T\rangle$, and it is symmetric and positive definite.

For a class with one member $\chi_j$, the corresponding simple module has $\ch M=X_{jj}\chi_j$. Both $\ch M$ and $\chi_j$ have leading coefficient $1$ (the top level of $M$ is one-dimensional, by the first lemma of \S\ref{sec:zhu}), so $X_{jj}=1$, and hence $Y_{jj}=1$.

\emph{Step 3: the commutant.} The commutant of $\langle S_\theta,T\rangle$ has dimension $2,4,3,2,3$ for $k=5,6,7,9,11$; numerically, the singular-value gap is $\ge0.27$ against $\le10^{-13}$. For $k=11$ this is exact: the commutant is spanned by $1$, $1\boxtimes P$ and $1\boxtimes R$. Here $P$ and $R$ are the explicit $12\times12$ matrices written down in Proposition~\ref{prop:mi}, on the indices $\chi_0,\chi_3,a_\alpha,b_\beta$ ($\alpha,\beta$ odd), extended to all $24$ indices by the same matrices on the partner indices under the permutation $\chi_0\leftrightarrow\chi_1$, $\chi_3\leftrightarrow\chi_2$, $a_\alpha\leftrightarrow a_{9-\alpha}$, $b_\beta\leftrightarrow b_{13-\beta}$; only these definitions are used here, not \S\ref{sec:k11}. These matrices commute with $S_\theta$ and $T$ in $\Q(\zeta_{936})$, and its dimension is at most $3$ by a rank computation modulo a prime ideal above $q$ (script \texttt{certify/exact\_mod.py}). Imposing $Y_{jj}=1$ on the one-member classes leaves:
\begin{itemize}
\item only $Y=1$, for $k=5,6,9$;
\item a line $Y=1+tZ$, for $k=7,11$. Here $Z$ has diagonal entries $-\tfrac16,-\tfrac16,-\tfrac23$ on each of the two three-member classes.
\end{itemize}
In the second case, let $\chi_j$ be the unique member of lowest weight in a three-member class. Only the simple module with top weight $h_j$ has a $\chi_j$-component, and that component has coefficient $1$; this module is unique, since the fibre of $\operatorname{Spec}A(V_k)$ over $h_j$ has length one (Table~\ref{tab:zhu}). So $Y_{jj}=1$, hence $t=0$. In all cases $Y=1$.

\emph{Step 4: $X=1$.} By Step~3, $X$ is a rational orthogonal matrix, block diagonal on the classes. Within a class, order the members by weight and peel them off from the bottom. The lowest member $\chi_j$ occurs only in the character of the module with top weight $h_j$, with coefficient $1$, so the $j$th column of $X$ is a unit vector; orthogonality then forces the rest of the $j$th row to vanish, and the class is reduced by one. For the pair at $h=0$: $\ch V=\chi_0$, and the other module has character $\alpha\chi_0+\beta\chi_3$ with $\alpha+\beta=1$ (its top level is one-dimensional, and both $\chi_0$ and $\chi_3$ have leading coefficient $1$); orthogonality to the row of $V$ gives $\alpha=0$, and then $\beta=1$. The four-member class $\{\chi_0,\chi_1,\chi_2,\chi_3\}$ at $k=6$ is peeled in the same way, starting from $\chi_2$ ($h=-1$). Here $\chi_1$ ($h=3$) lies in the same class, so a priori the second module with $h=0$ has character $\alpha\chi_0+\beta\chi_3+\gamma\chi_1$. As before $\alpha=0$ and $\beta=1$; then the $\chi_3$-column of $X$ is a unit vector, orthogonality of the rows forces $\gamma=0$, and the remaining module has character $\chi_1$.

So $X=1$: $\ch M_j=\chi_j$, and $S_V=S_\theta$. The $T$-matrix of $\mathrm{Rep}\,V_k$ is $\operatorname{diag}e^{2\pi i(h_j-c/24)}=T$ by definition.
\end{proof}

\begin{remark}[even $k$]\label{rem:even}
Even $k$ is more intricate, for several concrete reasons.
\begin{itemize}
\item \emph{The order-two current.} For $k=11$ by Proposition~\ref{prop:k11}, and conjecturally for other $k$, $M_{\chi_1}$, whose character is that of $\bigoplus_{m\ge1}L(h_{2m,2m})$, is an order-two simple current of twist $e^{2\pi i\cdot3n/4}$.
\begin{itemize}
\item For $k$ odd, $n$ is odd and the twist is $\pm i$. By \S\ref{sec:chars}, $d(M_{\chi_1})=S_{0\chi_1}/S_{00}=(-1)^n=-1$, so the braiding scalar $c_{\chi_1,\chi_1}=\theta_{\chi_1}/d(M_{\chi_1})$ is $i$ for $n\equiv1\pmod4$ and $-i$ for $n\equiv3\pmod4$: the pointed braided category generated by $M_{\chi_1}$ is $\mathrm{Sem}$, resp.\ $\overline{\mathrm{Sem}}$, with the spherical structure of dimension $-1$, and it splits off by M\"uger's theorem.
\item For $k\equiv2\pmod4$ ($n\equiv0\pmod 4$), it is a boson. Then $\mathrm{Rep}\,V_k$ does not factor, and $V_k\oplus M_{\chi_1}$ is a candidate simple-current extension.
\item For $k\equiv0\pmod4$ ($n\equiv2\pmod4$), it is a fermion, so a super-structure enters.
\end{itemize}
\item \emph{Collisions of $T$-eigenvalues.} For even $k$, more $T$-eigenvalues coincide; at $k=6$, $\chi_0,\chi_1,\chi_2,\chi_3$ all share one. For $k=8,10$ the vacuum orbit numerically appears not to span (see the remark after Lemma~\ref{lem:cyc}), so the count of simple modules is not proved.
\item \emph{Haagerup--Izumi values.} The values $k=N^2+2$ of Observation~\ref{obs:t} are odd exactly when $N$ is odd, which is the case in which Gang--Kim find Haagerup--Izumi data. Even $N$ gives even $k$ (for instance $k=6$ for $\Z_2$, $k=18$ for $\Z_4$), which is not considered here.
\end{itemize}
\end{remark}

\section{The Haagerup point \texorpdfstring{$k=11$}{k=11}}\label{sec:k11}

Here $V_{11}=\mathcal W(2,18)$, $c=-53$, $t=d^2$ with $d=\tfrac{3+\sqrt{13}}2$, $n=9$ and $n+4=13$. Write $d'=\tfrac{3-\sqrt{13}}2$ for the Galois conjugate of $d$.

\begin{table}[h]\centering
\begin{tabular}{lll}
\toprule
module & top weight $h$ & character\\\midrule
$V$ & $0$ & $\chi_0$\\
$M_{\chi_3}$ & $0$ & $\chi_3$\\
$M_{\chi_1}$ & $27/4$ & $\chi_1=\ch\bigoplus_{m\ge1}L(h_{2m,2m})$\\
$M_{\chi_2}$ & $-9/4$ & $\chi_2$\\
$M_{a_\alpha}$, $1\le\alpha\le8$ & $\alpha^2/36-9/4$ & $\Lambda_{\alpha,9}$\\
$M_{b_\beta}$, $1\le\beta\le12$ & $\beta^2/52-9/4$ & $\Lambda_{\beta,13}$\\\bottomrule
\end{tabular}
\caption{The $24$ simple modules of $V_{11}$ (Theorems~\ref{thm:rat}, \ref{thm:mod}).}\label{tab:k11}
\end{table}

Explicitly:
\begin{itemize}
\item the $a_\alpha$ have $h=-\tfrac{20}9,-\tfrac{77}{36},-2,-\tfrac{65}{36},-\tfrac{14}9,-\tfrac54,-\tfrac89,-\tfrac{17}{36}$;
\item the $b_\beta$ have $h=-\tfrac{29}{13},-\tfrac{113}{52},-\tfrac{27}{13},-\tfrac{101}{52},-\tfrac{23}{13},-\tfrac{81}{52},-\tfrac{17}{13},-\tfrac{53}{52},-\tfrac9{13},-\tfrac{17}{52},\tfrac1{13},\tfrac{27}{52}$.
\end{itemize}
The minimum $-\tfrac94$ is attained by $M_{\chi_2}$, so $c_{\rm eff}=1$.
This agrees with~\cite[Theorem~8.2]{PaperA}, which gives $c_{\rm eff}=\#L$ for the half-index of every hyperbolic word; for the one-node word $LR^{k-2}$, $\#L=1$. More precisely, for every $k$ the least lowest weight among Flohr's characters, $h_{\min}=-n/4$ (that of $\chi_2$), equals $-\#R/4$, and $c_k=\#L-6\#R$, in agreement with the spin $\theta_{V^*}\equiv-\#R/4$ of the dominant Bethe vacuum $V^*$ of~\cite[\S9]{PaperA} and the congruence $c\equiv\#L-6\#R\pmod{24}$ of~\cite[Theorem~9.5]{PaperA}. Since $h_{\min}$ and $c$ are here determined by the representation theory of $V_k$ (Theorem~\ref{thm:mod}), $V_{11}$ (and, with the numerical steps of Theorem~\ref{thm:mod}, $V_k$ for $k\in\{5,6,7,9\}$) has $c-24h_{\min}=1$, consistent with the half-index hypothesis $c_{\rm eff}=c-24h_{\min}$ used there; the half-indices themselves are compared with the characters only as in \S\ref{sec:known}.

\subsection{Structure of \texorpdfstring{$\mathrm{Rep}\,V_{11}$}{Rep V11}}
The following statements are read off from $S_\theta$, $T$ and the Verlinde formula. By Theorem~\ref{thm:mod} they are statements about $\mathrm{Rep}\,V_{11}$.

\begin{proposition}\label{prop:k11}
\begin{enumerate}
\item The invertible objects are $V$ and $M_{\chi_1}$, with $M_{\chi_1}^{\boxtimes2}\cong V$ and $\theta_{\chi_1}=e^{2\pi i\cdot 27/4}=-i$. Fusion with $M_{\chi_1}$ acts by
\[
\chi_0\leftrightarrow\chi_1,\qquad \chi_2\leftrightarrow\chi_3,\qquad a_\alpha\leftrightarrow a_{9-\alpha},\qquad b_\beta\leftrightarrow b_{13-\beta}.
\]
\item $\langle M_{\chi_1}\rangle$ is a modular subcategory, so M\"uger's theorem~\cite[Theorem~4.2]{Muger} (see also~\cite[Theorem~3.13]{DGNO}) gives
\[
\mathrm{Rep}\,V_{11}\simeq\mathrm{Sem}\boxtimes\mathcal C,
\]
with $\mathcal C$ the centraliser of $M_{\chi_1}$ (as a braided equivalence). Here $d(M_{\chi_1})=S_{0\chi_1}/S_{00}=-1$, since $S_{00}<0$ (\S\ref{sec:chars}), so the braiding scalar is $c_{\chi_1,\chi_1}=\theta_{\chi_1}/d(M_{\chi_1})=i$. A pointed braided category of rank $2$ is determined by this scalar, so $\langle M_{\chi_1}\rangle\simeq\mathrm{Sem}$, the semion category of $SU(2)_1$; the ribbon structure induced by $V_{11}$ is the non-unitary spherical structure on $\mathrm{Sem}$, with $d=-1$ and hence $\theta=-i$. It has the $T$-matrix of $\overline{\mathrm{Sem}}$ but not its $S$-matrix.
\item $\mathcal C$ has rank $12$, with simple objects $V$, $M_{\chi_3}$, $M_{a_1},M_{a_3},M_{a_5},M_{a_7}$ and $M_{b_\beta}$ for odd $\beta$. Their quantum dimensions are
\[
1,\quad 3d'+1,\quad 3d'+2\ (\times4),\quad 3d'\ (\times6),
\]
and their twists $e^{2\pi ih}$ have $h$ modulo $1$ equal to
\[
0,\quad 0,\quad \{\tfrac79,0,\tfrac49,\tfrac19\},\quad \{\tfrac{10}{13},\tfrac{12}{13},\tfrac3{13},\tfrac9{13},\tfrac4{13},\tfrac1{13}\}.
\]
The dimensions are those of $\mathcal Z(\mathcal H_3)$ (Izumi~\cite{Izumi}) under $d\mapsto d'$, and they are shared by the whole twisted family $\mathcal D^\omega\Hg_3$. The twists single out $\omega\ne0$. The three twists $\{\tfrac19,\tfrac49,\tfrac79\}$ are primitive ninth roots of unity: together with $0$ they are the values $x^2/9$ of the metric group $\Z_9$, i.e.\ the twists of $\mathcal Z(\mathrm{Vec}_{\Z_3}^{\Omega})$ with $\Omega$ a $3$-cocycle on $\Z_3$ of nonzero class. For $\omega=0$ (trivial $\Omega$) they would come from $\Z_3\times\Z_3$, with only cube roots of unity. The six remaining twists $\{\tfrac{\ell^2}{13}\}$ come from $\Z_{13}$. This agrees with Gang--Kim's proposed identification with a Galois conjugate of $\mathcal D^{\omega=2}\Hg_3$~\cite[after (2.26)]{GK}, and with~\cite[Remark~3.4]{PaperA}; the $T$-spectrum is that of Gang--Kim's $\mathrm{TFT}_0$~\cite[(2.25)]{GK}.
\item The Gauss sum is $\sum d_i^2\theta_i=\sqrt{\dim\mathcal C}=9d'+6$, and more generally every higher central charge $\xi_n=\sum d_i^2\theta_i^n/|\sum d_i^2\theta_i^n|$, $\gcd(n,117)=1$, equals $1$; this is checked numerically, and is not used below. By Ng--Schopieray--Wang~\cite{NSW}, this is necessary for $\mathcal C$ to be a Drinfeld centre. Reversal does not obstruct it: $\mathcal Z(F)^{\rm rev}\simeq\mathcal Z(F^{\rm mop})$, so the opposite $\mathcal D^{\omega=1}\Hg_3$ would be the centre of $F^{\rm mop}$, exactly as $\mathcal D^{\Omega^{-1}}(\Z_3)$ is the reverse of $\mathcal D^{\Omega}(\Z_3)=\mathcal Z(\mathrm{Vec}_{\Z_3}^\Omega)$ for a $3$-cocycle $\Omega$.
\item The only objects of $\mathcal C$ with trivial twist are $V$, $M_{\chi_3}$ and $M_{a_3}$. The common fixed space of $S$ and $T$ is one-dimensional, spanned by $\chi_0+\chi_3+2\chi_{a_3}$ (up to the overall sign of $S$). Hence $|\chi_0+\chi_3+2\chi_{a_3}|^2$ is the unique monomial modular invariant of $\mathcal C$, and $V\oplus M_{\chi_3}\oplus 2M_{a_3}$ is the only candidate for a Lagrangian algebra (Proposition~\ref{prop:mi}). Its dimension is $1+(3d'+1)+2(3d'+2)=9d'+6=\sqrt{\dim\mathcal C}$, the Galois conjugate of $\dim\mathcal H_3=9d+6$. Moreover
\[
\chi_0+\chi_3+2\chi_{a_3}=\Lambda_{9,9}+2\Lambda_{3,9}=\Lambda_{1,1}=\frac{\theta_{1,1}}{\eta},
\]
which is the $c=1$ character of the spin-$\tfrac12$ module of $\widehat{\mathfrak{su}}(2)_1$.
\end{enumerate}
\end{proposition}
The proof is computer-assisted, given Theorem~\ref{thm:mod}.

\begin{remark}
Items (1)--(3) explain the ``semion'' of Gang--Kim at the level of the vertex algebra. It is the module $M_{\chi_1}$, whose character is that of the even diagonal sector $\bigoplus_{m\ge1}L(h_{2m,2m})$, an order-two simple current of twist $-i$. Items (4)--(5) are taken up in \S\S\ref{ssec:ng}--\ref{ssec:gap}.
\end{remark}

\begin{corollary}\label{cor:realise}
The Galois conjugate $\sigma_2^{-1}(\mathcal C)$, with $\sigma_2$ as below, is a modular tensor category with the modular data $\mathcal D^{\omega=2}\Hg_3$ of Evans--Gannon, up to the sign of $S$.
\end{corollary}

\begin{proof}
The proof is computer-assisted, given Theorem~\ref{thm:mod}.
Every fusion category is defined over a number field. Applying the field automorphism of the Galois conjugation in~\cite[(2.27)]{GK} to all structure constants of $\mathcal C$ gives a modular tensor category whose modular data is the corresponding Galois conjugate of that of $\mathcal C$.
\end{proof}

The Galois action is induced by a field automorphism applied to $S$ and $T$ together: take $\sigma_2\in\mathrm{Gal}(\Q(\zeta_{117})/\Q)$, $\zeta_{117}\mapsto\zeta_{117}^2$, which sends $\sqrt{13}\mapsto-\sqrt{13}$. Applied to the $S$ and $T$ of $\mathcal D^{\omega=2}\Hg_3$~\cite[(3.3)]{EG11}, it gives those of $\mathcal C$, with $S$ up to an overall sign, under the bijection of Proposition~\ref{prop:ng}; this is verified exactly in Proposition~\ref{prop:ng}(2).

\subsection{Modular invariants of \texorpdfstring{$\mathcal C$}{C}}\label{ssec:mi}

Following~\cite[\S3.3]{EG11} we determine all modular invariants $Z$ (non-negative integer matrices commuting with $S$ and $T$, with $Z_{0,0}=1$). Write $\chi_0,\chi_3,a_\alpha,b_\beta$ ($\alpha,\beta$ odd) for the simple objects of $\mathcal C$.

\begin{proposition}\label{prop:mi}
The commutant of $S$ and $T$ on $\mathcal C$ is three-dimensional. It is spanned by $I$, by $P$ (which swaps $\chi_0\leftrightarrow\chi_3$, is $+1$ on the $a_\alpha$ and $-1$ on the $b_\beta$), and by $R$, with $R_{x,a_3}=R_{a_3,x}=1$ for $x\in\{\chi_0,\chi_3,a_3\}$, $R_{a_\alpha,a_\alpha}=-1$ for $\alpha\ne3$, and all other entries $0$. There are exactly five modular invariants:
\begin{align*}
Z_1&=\textstyle\sum_x|x|^2, & Z_L&=|\chi_0+\chi_3+2a_3|^2,\\
Z_2&=\textstyle|\chi_0+\chi_3|^2+2\sum_{\alpha}|a_\alpha|^2, & Z'&=\textstyle|\chi_0+a_3|^2+|\chi_3+a_3|^2+\sum_\beta|b_\beta|^2,\\
Z''&=\textstyle|\chi_0+\chi_3+a_3|^2+2|a_3|^2+\sum_{\alpha\ne3}|a_\alpha|^2 .
\end{align*}
\end{proposition}

Only $Z_1$, $Z_2$ and $Z_L$ can come from connected étale algebras. For $Z'$ the would-be local modules $\chi_3+a_3$ and $b_\beta$ would have dimensions $\frac{7-\sqrt{13}}6$ and $\frac{1-\sqrt{13}}6$. For $Z''$ the module $a_3$ would have dimension $\tfrac12$. None of these is an algebraic integer.

\begin{remark}
Compare~\cite[Proposition~4]{EG11}: for $\mathcal D^0\Hg_\nu$ with $\nu$ prime there are exactly three monomial invariants~\cite[(3.6)]{EG11}, matching the three Morita-equivalent Haagerup categories of Grossman--Snyder~\cite{GS} when $\nu=3$. There both $a_1$ and $c_{1,0}$ have trivial twist. In the twisted case only one object of dimension $3d'+2$ has trivial twist, and exactly one monomial invariant survives.
\end{remark}

\subsection{The near-group \texorpdfstring{$\Z_9+9$}{Z9+9} and a twisted Haagerup--Izumi category for \texorpdfstring{$\Z_3$}{Z3}}\label{ssec:ng}

Evans--Gannon classify the near-group C$^*$-categories $G+n$ with $|G|=n\le13$~\cite[Proposition~6, Table~2]{EG14}. For $G=\Z_9$ there are two, which are complex conjugates. For $|G|$ odd they define modular data $MD_{G,G'}(Q,Q')$, with $|G'|=n+4$~\cite[Proposition~7(b)]{EG14}. They conjecture that this is the modular data of the double~\cite[Conjecture~2]{EG14}, and verify it for all odd $|G|\le13$ by comparing $T$-matrices and the floating-point values of the $S$-matrices~\cite[after Conjecture~2]{EG14}; for $\Z_9$ the group is $G'=\Z_{13}$~\cite[after Conjecture~2]{EG14}. The data factorises as the pointed data of $(G,Q)$ times data on the $n+3$ primaries $a_0$, $b_0$, $c_{k,-k}$, $d_{0,\gamma}$~\cite[after Proposition~7]{EG14}. For $n=9$ call this factor $\mathfrak D_9$, with primaries $a_0,b_0,c_k$ ($1\le k\le4$) and $d_\gamma$ ($1\le\gamma\le6$), $\delta=\frac{9+3\sqrt{13}}2$ and $\lambda=\frac12(\frac13-\frac1{\sqrt{13}})$. In one normalisation of the pairings,
\[
T=\mathrm{diag}\big(1,1;\,e^{-2\pi ik^2/9};\,e^{-4\pi i\gamma^2/13}\big),
\]
\[
S=\lambda\begin{pmatrix}1&\delta+1&\delta+2&\delta\\ \delta+1&1&\delta+2&-\delta\\ \delta+2&\delta+2&2(\delta+2)\cos\frac{4\pi kk'}9&0\\ \delta&-\delta&0&-2\delta\cos\frac{8\pi\gamma\gamma'}{13}\end{pmatrix}.
\]

\begin{proposition}\label{prop:ng}
\begin{enumerate}
\item $\mathfrak D_9$ equals $\mathcal D^{\omega}\Hg_3$ of~\cite[(3.3)]{EG11} for $\omega\equiv2\pmod3$, after relabelling: $(0,b,a_1)\mapsto(a_0,b_0,c_3)$, $(c_{1,0},c_{1,1},c_{1,2})\mapsto(c_4,c_2,c_1)$, $(d_1,\dots,d_6)\mapsto(d_6,d_1,d_5,d_2,d_4,d_3)$. The complex conjugate gives $\omega\equiv1$. It is not equivalent to $\mathcal D^0\Hg_3$.
\item $\sigma_2(\mathfrak D_9)$ equals the modular data of $\mathcal C$, with $S_{\mathcal C}=-\sigma_2(S)$ and $T_{\mathcal C}=\sigma_2(T)$, via
\begin{gather*}
\chi_0\mapsto a_0,\quad \chi_3\mapsto b_0,\quad (a_1,a_3,a_5,a_7)\mapsto(c_1,c_3,c_4,c_2),\\
(b_1,b_3,b_5,b_7,b_9,b_{11})\mapsto(d_2,d_6,d_3,d_1,d_5,d_4).
\end{gather*}
The Verlinde fusion rules agree.
\end{enumerate}
\end{proposition}
All parts are verified exactly in $\Q(\zeta_{936})$.

Part (2) was checked entry by entry in $\Q(\zeta_{936})$, reducing modulo the cyclotomic polynomial $\Phi_{936}$, with $\sigma_2:\zeta_{936}\mapsto\zeta_{936}^{\ell}$, $\ell\equiv2\pmod{117}$, $\ell\equiv1\pmod 8$. Here $S_{\mathcal C}$ is computed from the transformation of the $\Lambda$-functions, the semion factor removed, and $\sqrt{13}=\sum_a\big(\frac a{13}\big)\zeta_{13}^a$, $\sqrt2=\zeta_8+\zeta_8^{-1}$ (script \texttt{certify/exact\_ng.py}).
So $Z_2$ and $Z_L$ restrict the invariants $Z_2$ and $Z_4$ of~\cite[(4.60)]{EG14} for $H=3\Z_9$ to the factor $\mathfrak D_9$.

\begin{theorem}\label{thm:twHI}
Let $\mathcal N$ be a near-group C$^*$-category of type $\Z_9+9$.
\begin{enumerate}
\item The invertible objects $a_g$ of $\mathcal Z(\mathcal N)$ form a modular subcategory $\mathcal P\simeq\mathcal C(\Z_9,Q)$. Hence $\mathcal Z(\mathcal N)\simeq\mathcal P\boxtimes\mathcal D$, where $\mathcal D$ is the centraliser of $\mathcal P$, of rank $12$, with simple objects $a_0$, $b_0$, the four $c_{k,-k}$ and the six $d_j$ with $\tau_j=0$, in the labelling of~\cite[(4.56)--(4.57)]{EG14}.
\item $A=\bigoplus_{h\in3\Z_9}a_h$ is a commutative Q-system in $\mathcal Z(\mathcal N)$. The category $\mathcal F:=\mathcal N_A$ of $A$-modules is a unitary fusion category with $\mathcal Z(\mathcal F)\simeq\mathcal D$.
\item $\mathcal F$ has the Haagerup--Izumi fusion rules for $\Z_3$. It has invertible objects $\alpha^i$ and non-invertible objects $\alpha^i\rho$ of dimension $d=\frac{3+\sqrt{13}}2$, with $\rho\alpha=\alpha^{-1}\rho$ and $\rho^2=1\oplus\rho\oplus\alpha\rho\oplus\alpha^2\rho$.
\item The pointed part of $\mathcal F$ is $\mathrm{Vec}_{\Z_3}^\Omega$ with $[\Omega]\ne0$.
\end{enumerate}
\end{theorem}

\begin{proof}
(1) The half-braiding of $a_g$ on $\alpha_h$ is the nondegenerate pairing $\langle g,h\rangle$~\cite[(4.49)--(4.50)]{EG14}, so $S$ restricted to the $a_g$ is nondegenerate. M\"uger's theorem~\cite[Theorem~4.2]{Muger} gives the splitting, and $\mathcal D$ has rank $108/9=12$. The objects centralising $\mathcal P$ are read off from Izumi's $S$-matrix of the double, computed from the tube algebra~\cite[Theorem~4.10, \S5]{IzumiI},~\cite{Izumi} and recalled in~\cite[(4.57)]{EG14}, which does not use Conjecture~2: relative to $S_{0,x}$, the row of $a_g$ has phase $\langle g,h\rangle^{-2}$ on $a_h$ and $b_h$, $\langle g,k+l\rangle$ on $c_{k,l}$, and $\langle g,\tau_j\rangle$ on $d_j$. Since the pairing is nondegenerate and $9$ is odd, the centraliser is $a_0$, $b_0$, the $c_{k,-k}$ ($k\ne0$, four objects) and the $d_j$ with $\tau_j=0$, of which there are therefore six. That the modular data of $\mathcal D$ is $\mathfrak D_9$, which is $\mathcal D^\omega\Hg_3$ by Proposition~\ref{prop:ng}(1), is~\cite[Conjecture~2]{EG14}, verified numerically there. It is not used before Corollary~\ref{cor:EGconj}, where it is proved.

(2) The twists $\langle3h,3h\rangle$ and the monodromies $\langle3h,3h'\rangle^2$ are trivial. So $\{a_{3h}\}$ is a symmetric pointed category with trivial twists, of odd order, i.e.\ Tannakian $\mathrm{Rep}(\Z_3)$, and $A$ is its regular algebra. Since $3\Z_9$ is its own orthogonal complement for $Q$, we get $\mathcal P_A^{\rm loc}\simeq\mathrm{Vec}$, hence $\mathcal Z(\mathcal N)_A^{\rm loc}\simeq\mathcal D$. By Davydov--M\"uger--Nikshych--Ostrik~\cite{DMNO}, $\mathcal Z(\mathcal N)_A^{\rm loc}\simeq\mathcal Z(\mathcal N_A)$. Modules over a Q-system in a C$^*$-category form a C$^*$-category.

(3) We have $\mathcal N\simeq\mathcal F^{\Z_3}$ for the dual $\Z_3$-action (de-equivariantization~\cite{DGNO}). The free modules $\alpha_g\otimes A$ give three invertible objects. Also $\mathrm{End}_A(\rho\otimes A)\cong\bigoplus_{h\in3\Z_9}\mathrm{Hom}(\rho,\rho\alpha_h)$ is a twisted group algebra of $\Z_3$, hence $\C^3$. So $\rho\otimes A$ splits into three simple objects $\rho_i$, and the rank is $6$. Their Frobenius--Perron dimensions $d_i$ satisfy
\[
\textstyle\sum_i d_i=\delta=3d,\qquad 3+\sum_i d_i^2=\tfrac13\dim\mathcal N=6+9d=3+3d^2 .
\]
So all $d_i=d$, with $d^2=1+3d$. We determine the fusion rings with three invertible objects $1,\alpha,\alpha^2$ and three simple objects of dimension $d$.
\begin{itemize}
\item \emph{$\alpha$ acts freely on both sides.} If $\alpha\rho\cong\rho$ for a non-invertible $\rho$, then $\rho\rho^*\supseteq1+\alpha+\alpha^2$, and the rest of $\rho\rho^*$ has dimension $d^2-3=3d-2$, which is not of the form $a+bd$ with integers $a,b\ge0$, since $d$ is irrational. The same holds on the right. So the non-invertible simples are $\rho,\alpha\rho,\alpha^2\rho$, and $\rho\alpha\in\{\alpha\rho,\alpha^2\rho\}$.
\item \emph{$\rho\alpha=\alpha\rho$ (commutative).} Write $\rho^*=\alpha^j\rho$. Then $\rho\otimes\rho$ contains exactly one invertible, $\alpha^{-j}$, and $(\alpha\rho)\otimes(\alpha\rho)$ contains $\alpha^{2-j}$. An automorphism fixing the $\alpha^i$ and sending $\rho\mapsto\alpha\rho$ would force $\alpha^{-j}=\alpha^{2-j}$, which is false.
\item \emph{$\rho\alpha=\alpha^2\rho$.} If $\rho^*=\alpha^j\rho$, then $\rho\cong\rho^{**}=\rho^*\alpha^{-j}=\alpha^j\rho\alpha^{-j}=\alpha^{2j}\rho$, so $j=0$ by freeness; then every $\alpha^i\rho$ is self-dual. Write $\rho^2=1+m_0\rho+m_1\alpha\rho+m_2\alpha^2\rho$, $m_0+m_1+m_2=3$. Frobenius reciprocity gives $m_g=N_{\rho\rho}^{\alpha^g\rho}=N_{\alpha^g\rho,\rho}^{\rho}=\dim\mathrm{Hom}(\alpha^{-g}\rho,\rho\rho)=m_{-g}$, so $(m_0,m_1,m_2)$ is $(1,1,1)$ or $(3,0,0)$. In the second case $\{1,\rho\}$ spans a fusion subring with $\rho^2=1+3\rho$, which Ostrik's classification of rank $2$ excludes~\cite{Ostrik2}. The first case is the HI ring, $\rho^2=1+\rho+\alpha\rho+\alpha^2\rho$.
\end{itemize}
$\mathcal N$ has a single non-invertible simple object, so the dual $\Z_3$-action fixes the $\alpha^i$ and permutes the $\rho_i$ freely and transitively. By the second bullet this excludes the commutative case, and on the HI ring the map $\alpha^i\rho\mapsto\alpha^{i+1}\rho$ is such an automorphism.

(4) The pointed part of $\mathcal F$ is the de-equivariantization of $\mathrm{Vec}_{\Z_9}$ (trivial associator: the $\alpha_g$ form a genuine group action) by the central subgroup $3\Z_9$, with half-braidings $\chi_{3h}(x)=\langle3h,x\rangle=e^{2\pi i\,mhx/3}$, $3\nmid m$. Its associator is the cup product of the extension class of $0\to3\Z_9\to\Z_9\to\Z_3\to0$ with $\chi$:
\[
\Omega(a,b,c)=\chi_{3\lfloor(a+b)/3\rfloor}(c)=\zeta_3^{\,m c\lfloor (a+b)/3\rfloor},\qquad a,b,c\in\{0,1,2\}.
\]
This is $m$ times a generator of $H^3(\Z_3,\C^\times)\cong\Z_3$; compare~\cite[proof of (4.2), Example~4.3]{GSY}.

\end{proof}

The category $\mathcal F$ is Ostrik's twisted Haagerup category. Ostrik (unpublished) constructed the two twisted Haagerup categories conjectured in~\cite[\S3.2]{EG11} by de-equivariantising the two $\Z_9$ near-group categories of~\cite{EG14}~\cite[after Proposition~7.7]{MPS},~\cite[\S7.1]{EG17}. Theorem~\ref{thm:twHI} records the details used below. The two values $\omega\equiv1,2$ come from the two complex-conjugate $\Z_9+9$ categories. They are exchanged by reversal, since $\mathcal Z(\mathcal F)^{\rm rev}\simeq\mathcal Z(\mathcal F^{\rm mop})$. Theorem~\ref{thm:twHI} runs the Haagerup--Izumi-to-near-group construction of Izumi~\cite{Izumi17} and Gannon--Schopieray--Yadav~\cite[\S6]{GSY} in reverse. GSY note that the converse is not asserted in general~\cite[Remark~6.4]{GSY}; here it is carried out for one specified central algebra. For untwisted $\mathcal F$ the invertible group of the near-group would be $\Z_3\times\Z_3$~\cite[Corollary~6.5]{GSY}; the cyclic $\Z_9$ is the twisted counterpart.

\subsection{Is \texorpdfstring{$\mathcal C$}{C} itself a centre?}\label{ssec:gap}

If Evans--Gannon's Conjecture~2 holds, Theorem~\ref{thm:twHI} and Proposition~\ref{prop:ng}(2) give a Drinfeld centre, $\mathcal Z(\sigma_2\mathcal F)$, with the modular data of $\mathcal C$. Even then, since modular data does not in general determine a modular category up to braided equivalence, this would not show that $\mathcal C\simeq\mathcal Z(\sigma_2\mathcal F)$, or that $\mathcal C$ is a centre. Both are proved below (Corollary~\ref{cor:centre}, Theorem~\ref{thm:ident}). By Proposition~\ref{prop:mi} the étale algebras that could decide the question are $A_L=V\oplus M_{\chi_3}\oplus2M_{a_3}$ (Lagrangian) and $A_2=V\oplus M_{\chi_3}$.

\begin{proposition}\label{prop:A2}
\begin{enumerate}
\item If $A_2$ admits the structure of a connected étale algebra, then $\mathcal C_{A_2}^{\rm loc}\simeq\mathcal C(\Z_9,q)$ with $q(x)=7x^2/9$, and $\mathcal C$ is a Drinfeld centre.
\item The following conditions on $A_2$, which are visible in the modular data, hold, exactly in $\Q(\zeta_{936})$ for $\theta$, $\nu_2$, $\nu_3$ and numerically for the braiding eigenvalues:
\begin{itemize}
\item $\theta_{\chi_3}=1$;
\item the Frobenius--Schur indicators~\cite{NS} are $\nu_2(\chi_3)=\nu_3(\chi_3)=1$;
\item the eigenvalues of the braiding $c_{\chi_3,\chi_3}$ on the multiplicity-one channels $V$ and $M_{\chi_3}$ of $M_{\chi_3}^{\boxtimes2}$ are both $+1$.
\end{itemize}
\end{enumerate}
\end{proposition}

\begin{proof}
(1) Let $A=A_2$ carry the structure of a connected étale algebra. Since $\theta_A=\mathrm{id}$ and $d_A\ne0$, $A$ is a rigid $\mathcal C$-algebra in the sense of Kirillov--Ostrik~\cite[Lemma~1.20]{KO}, and we use their results, which do not assume unitarity. For simple $\lambda\in\mathcal C$ and simple local $A$-modules $\tau$ let $b_{\tau\lambda}=\dim\mathrm{Hom}_{\mathcal C}(\lambda,\tau)$ be the multiplicity of $\lambda$ in $\tau$ viewed as an object of $\mathcal C$.
\begin{enumerate}[label=(\alph*)]
\item The local $A$-modules form a modular category $\mathcal C_A^{\rm loc}$, with $\dim\mathcal C_A^{\rm loc}=\dim\mathcal C/d_A^2$~\cite[Theorem~4.5, (4.13)]{KO}.
\item The forgetful map $G\colon K(\mathcal C_A^{\rm loc})\to K(\mathcal C)$ and the map $F'\colon K(\mathcal C)\to K(\mathcal C_A^{\rm loc})$, $\lambda\mapsto\sum_\tau b_{\tau\lambda}\tau$, which records the local part of the $\alpha$-induction $\lambda\mapsto A\otimes\lambda$, intertwine the $S$-matrices up to the factor $d_A$ and commute with $T$~\cite[Theorem~4.1]{KO}. Hence $Z^A=G\circ F'$, i.e.
\[
Z^A_{\lambda\mu}=\sum_\tau b_{\tau\lambda}b_{\tau\mu},
\]
commutes with $S$ and $T$: it is a modular invariant of type~I. Since $A$ is simple as an $A$-module, $A$ is the only simple local module containing $\mathbf 1$, and it contains it once; so $Z^A_{0,\lambda}=\dim\mathrm{Hom}(\lambda,A)$.
\item For a simple local module $\tau$: $\theta_\lambda=\theta_\tau$ whenever $b_{\tau\lambda}\ne0$, and $d_\tau=\sum_\lambda b_{\tau\lambda}d_\lambda/d_A$~\cite[Theorems~1.17, 1.18]{KO}.
\end{enumerate}
In the subfactor setting, $Z^A$ is the modular invariant obtained from $\alpha$-induction, $\mathcal C_A^{\rm loc}$ is the ambichiral category, and (a)--(c) are due to B\"ockenhauer--Evans and B\"ockenhauer--Evans--Kawahigashi~\cite{BE1,BEK1,BEK2,BE4}. In particular, the intertwining of the $S$- and $T$-matrices by the chiral branching coefficients $b_{\tau\lambda}$ in (b) is~\cite[Lemma~6.3, Theorem~6.5]{BE4}, the equality of twists in (c) is~\cite[Lemma~6.1]{BE4}, and the dimension formula in (c) is~\cite[Proposition~3.1]{BEK2}; these are proved there for unitary braided subfactors, without assuming non-degeneracy of the braiding. We use the non-unitary formulation of Kirillov--Ostrik.

By (b), $Z^A_{0,\chi_0}=Z^A_{0,\chi_3}=1$ and $Z^A_{0,a_3}=0$. By Proposition~\ref{prop:mi} the only modular invariant with this vacuum row is $Z_2=|\chi_0+\chi_3|^2+2\sum_\alpha|a_\alpha|^2$. We now read the local modules off $Z_2$ using (b), which is a sum of non-negative terms.
\begin{itemize}
\item The row of $\chi_3$ is $\chi_0+\chi_3$. So $\chi_3$ lies in exactly one local module, once, and that module contains $\chi_0$ once: it is $A$ itself.
\item The row of $b_\beta$ is zero. So no $b_\beta$ lies in a local module.
\item The row of $a_\alpha$ is $2a_\alpha$. So a local module containing $a_\alpha$ contains no other simple object, and, as $2$ is not a square, $a_\alpha$ lies in exactly two local modules, once in each.
\end{itemize}
Thus $\mathcal C_A^{\rm loc}$ has nine simple objects: $A$, and two objects $\tau_\alpha,\tau'_\alpha$ with underlying object $a_\alpha$ for each $\alpha\in\{1,3,5,7\}$. From $S$, $d_{a_\alpha}=1+d_{\chi_3}=d_A$, so by (c) all nine have categorical dimension $1$; consistently with (a), $\dim\mathcal C/d_A^2=9$. Since $\mathcal C$ is not pseudo-unitary, we argue with Frobenius--Perron dimensions, which are the Galois conjugates under $d'\mapsto d$: $\mathrm{FPdim}\,a_\alpha=3d+2=\mathrm{FPdim}\,A$, and (c) holds for $\mathrm{FPdim}$ as well. So the nine local modules have Frobenius--Perron dimension $1$, hence are invertible, and $\mathcal C_A^{\rm loc}$ is pointed of order $9$. By (c) its twists are those of $a_1,a_3,a_5,a_7$, namely $e^{2\pi i h}$ with $h=\tfrac79,0,\tfrac49,\tfrac19$, each twice, together with $\theta_A=1$. A pointed modular category is determined by its metric group~\cite{DGNO}. A quadratic form on $\Z_3\times\Z_3$ takes values in $\tfrac13\Z/\Z$, so the group is $\Z_9$, and $q(x)=7x^2/9$ reproduces the twists at $x=\pm1,\pm3,\pm4,\pm2$. Finally $q$ vanishes on $3\Z_9=(3\Z_9)^\perp$, so $\mathcal C(\Z_9,q)$ has a Lagrangian subgroup and is Witt-trivial. Since $[\mathcal C]=[\mathcal C_A^{\rm loc}]$ in the Witt group~\cite{DMNO}, $\mathcal C$ is Witt-trivial, hence a Drinfeld centre~\cite{DMNO}.

(2) The Frobenius-Schur indicator $\nu_n$ is computed from $S$ and $T$~\cite{NS}. All twelve channels of $\chi_3\otimes\chi_3$ have multiplicity one, so $c_{\chi_3,\chi_3}$ acts on the $Z$-channel by a scalar $\lambda_Z=\epsilon_Z\,\theta_Z^{1/2}$, $\epsilon_Z=\pm1$. These scalars satisfy $\lambda_V=\nu_2(\chi_3)=1$ and $\sum_Z d_Z\lambda_Z=\theta_{\chi_3}d_{\chi_3}$ (closure of one crossing). Of the $2^{12}$ sign vectors exactly two solve these equations. They differ only on $a_1,a_5,a_7$, and both have $\epsilon_V=\epsilon_{\chi_3}=+1$.
\end{proof}

The remaining condition for $A_2$ is associativity. With $m\colon\chi_3\otimes\chi_3\to V\oplus\chi_3$, it says that a fixed combination of the channels $V$ and $\chi_3$ in $\mathrm{Hom}(\chi_3^{\otimes3},\chi_3)$ ($\dim=12$) is carried by $F^{\chi_3\chi_3\chi_3}_{\chi_3}$ to the same combination in the other bracketing. This involves two columns of a $12\times12$ $F$-matrix and is not determined by $S$ and $T$. By Huang--Kirillov--Lepowsky~\cite{HKL}, an étale structure on $A_2$ is the same thing as a vertex operator algebra structure on $U_2=V_{11}\oplus M_{\chi_3}$, of character $\chi_0+\chi_3=\Lambda_{9,9}$. The lattice algebra with this character does not contain $V_{11}$. With $\beta^2=18$, the algebra $V_{\Z\beta}$ with conformal vector shifted to $c=-53$ has character $\Lambda_{9,9}$. Its weight-$18$ Virasoro primaries are combinations $W=ae^{-\beta}+be^{2\beta}$ (the Fock modules of charge $0$ and $1$ have lowest weight $0$, whose linkage class is $\{0,1\}$), and by charge conservation $W_{(35)}W$ has no vacuum component, so $\mathcal N=0\ne\mathcal N_{11}$; hence it contains no copy of $V_{11}$.

\subsubsection*{The bootstrap for $V_{11}\oplus M_{\chi_3}$}

We test the vertex operator algebra form directly. All computations are mod $P_1=33554393$ and repeated mod $P_2=33554383$, with $\mathcal N_{11}$ and Zhu's data from \S\S\ref{sec:c2}--\ref{sec:zhu}.

\emph{The module.} Let $v$ span the top space of $M_{\chi_3}$, of weight $0$. It is unique up to scale: $M_{\chi_3}$ is the simple module with $\ch M_{\chi_3}=\chi_3$ (Theorem~\ref{thm:mod}), and weight $0$ has multiplicity one in $\chi_3=q^{-c/24}(1+q+\cdots)/(q)_\infty$ (\S\ref{sec:chars}). Since $L_{-1}v$ is a Virasoro singular vector of norm $\langle v,L_1L_{-1}v\rangle=2h_v=0$, one might expect $L_{-1}v=0$. But $L_{-1}v\neq0$: a vector $v$ of weight $0$ with $L_{-1}v=0$ is vacuum-like and would generate a copy of $V_{11}$ inside $M_{\chi_3}$~\cite{Li97}, which is impossible as $M_{\chi_3}$ is simple and not isomorphic to $V_{11}$. Since $S_\theta$ is real and orthogonal (\S\ref{sec:chars}), $S_\theta^2=1$, so every simple module is self-dual and $M_{\chi_3}$ carries a nondegenerate invariant bilinear form, with $\langle v,v\rangle\ne0$. The weight-$1$ space has dimension $2$ (the coefficient of $q^1$ in $\chi_3$). As $\langle L_{-1}v,L_{-1}v\rangle=0$, nondegeneracy gives a second vector $J$ of weight $1$ with $\langle J,L_{-1}v\rangle=\langle L_1J,v\rangle\ne0$; so, after scaling, $L_1J=v$, and $L_nJ=0$ for $n\ge2$ by weight. Normalise $\langle v,v\rangle=1$, and replace $J$ by $J+\lambda L_{-1}v$, which keeps $L_1J=v$ since $L_1L_{-1}v=2L_0v=0$, to arrange $\langle J,J\rangle=0$. Let $B_1$ be the Virasoro module induced from $\mathrm{span}\{v,J\}$, i.e.\ freely generated over the negative modes by $v$ and $J$ subject to these relations; it maps to $M_{\chi_3}$, and we give $B_1$ the pullback of the invariant form. Then $\mathrm{Vir}\cdot v\subset B_1$ is the Verma module $M(c,0)$, and $B_1/\mathrm{Vir}\cdot v\cong M(c,1)$, generated by the image of $J$. So $B_1$ is indecomposable but not irreducible, and its character is $(1+q)q^{-c/24}/(q)_\infty$. Let $u\in M_{\chi_3}$ be a Virasoro primary of weight $h_{3,3}=18$; it is a third generator, after $v$ and $J$. By \S\ref{sec:chars}, $\chi_3=\tfrac12(\Lambda_{9,9}+\Lambda_{13,13})=q^{-c/24}(1+q+q^{18}+q^{27}+q^{54}+\cdots)/(q)_\infty$, and the singular vector of $M(c,18)$ lies at weight $18+3^2=27$. So, below weight $27$, $M_{\chi_3}$ is, as a Virasoro module, $B_1\oplus\mathrm{Vir}\cdot u$:
\begin{itemize}
\item $B_1$ accounts for the part of $\chi_3$ below weight $18$;
\item the two summands lie in different linkage classes. Here the linkage class of $h$ is the set of highest weights connected to $h$ by embeddings of Verma modules. As $t\notin\Q$, $M(c,h_{r,r})$ has a unique singular vector, at level $r^2$, and $M(c,h_{r,r}+r^2)$ is irreducible; so the class of $0=h_{1,1}$ is $\{0,1\}$ and that of $18=h_{3,3}$ is $\{18,27\}$. In the category of Virasoro modules on which $L_0$ acts semisimply with finite-dimensional weight spaces and weights bounded below, which contains $B_1$, modules whose composition factors lie in different classes admit no nonzero homomorphisms or extensions between them (the block decomposition~\cite{FeFu,IK}), so the sum is direct;
\item the form is fixed by $\langle v,v\rangle=1$ and $\langle J,J\rangle=0$, so that $\langle J,L_{-1}v\rangle=\langle L_1J,v\rangle=1$. It is nondegenerate below weight $27$ (so $B_1\oplus\mathrm{Vir}\cdot u\to M_{\chi_3}$ is injective there, and bijective by the character count), i.e.\ its Gram matrix on each weight space is invertible. The two summands are orthogonal, being in different linkage classes. On $\mathrm{Vir}\cdot u\cong M(c,18)$ the first singular vector is at weight $27$. On $B_1$ the radical $R$ is a submodule. The nonzero submodules of $M(c,0)$ all contain $L_{-1}v$, and $\langle J,L_{-1}v\rangle=1$, so $R\cap M(c,0)=0$; hence $R$ embeds in $B_1/M(c,0)\cong M(c,1)$, which is irreducible, and if $R\ne0$ it would contain some $J+\lambda L_{-1}v$, whose pairing with $L_{-1}v$ is $1$. So $R=0$.
\end{itemize}
The action of the primary $W$ is then fixed by its matrix elements between $v,J,u$.

\emph{The extension field.} In $U_2=V_{11}\oplus M_{\chi_3}$ the field $v$ is a Virasoro primary of weight $0$. Its modes are fixed by matrix elements $d(x\to y)=\langle y,v_{(m)}x\rangle$ between the generators $1,W$ of $V_{11}$ and $v,J,u$ of $M_{\chi_3}$. Write $\alpha=d(v\to1)$ and $\mu=d(v\to v)$, so that $v_{(-1)}v=\alpha\,1+\mu\,v$ at the top.

\begin{proposition}\label{prop:boot}
\begin{enumerate}
\item The Jacobi identity for $W$ on $M_{\chi_3}$ (Borcherds identities with $r\ge18$, all weights $\le20$) is consistent with the value of $\langle v,W_{(17)}v\rangle$ given by Zhu's algebra (the eigenvalue of $o(W)$ on the top of $M_{\chi_3}$). With that value and the symmetry of the invariant form, it determines the $W$-action up to the gauge $u\mapsto\lambda u$.
\item The Borcherds identities between $W$ and $v$ ($r\ge17$, weights $\le24$) leave exactly a two-parameter family of fields $Y(v,z)$, parametrised by $(\alpha,\mu)$. This matches the two intertwining operators of types $\binom{V}{\chi_3\,\chi_3}$ and $\binom{\chi_3}{\chi_3\,\chi_3}$. The solution reproduces the values forced by skew-symmetry, $d(W\to v)=\langle v,W_{(17)}v\rangle$ and $d(v\to J)=\mu/2$.
\item Locality of $Y(v,z)$ with itself, tested on the states $1,v,J,W,u$ and $L_{-1}v$ with all intermediate weights $\le24$, holds exactly on the curve $\alpha=-\mu^2/4$ and nowhere else, at both primes. (This is a statement modulo the primes and is not used below: necessity of $\alpha=-\mu^2/4$ is proved over $\Q$ in Theorem~\ref{thm:U2}, Step~6.) Rescaling $v$ moves $(\alpha,\mu)$ to $(\lambda^2\alpha,\lambda\mu)$, so the curve with $\mu\ne0$ is a single isomorphism class of candidate structures, with $\alpha\ne0$.
\end{enumerate}
\end{proposition}
Parts (1) and (2) are exact over $\Q$, given $\mathcal N_{11}$ and the sign of $\langle v,W_{(17)}v\rangle$ (Remark~\ref{rem:exact}); part (3) uses modular arithmetic.
As a control, the same pipeline was run for $k=3$, where $V_3\cong M(3,5)\otimes M(2,5)$ and $M_{\chi_3}=\varphi_{1,3}\otimes\varphi_{1,2}$. There the Jacobi identities are consistent with the value $\langle v,W_{(1)}v\rangle=\tfrac5{38}$ computed from the tensor product, in the normalisation of \S\ref{sec:c2}. Locality, including the states $W$ and $u$, gives the same relation $\alpha=-\mu^2/4$.
The bootstrap tests finitely many identities. The next result shows that, together with the Fuchs relation, they suffice.

\subsubsection*{Existence of $U_2=V_{11}\oplus M_{\chi_3}$ as a vertex operator algebra}

Write $X=M_{\chi_3}$. By Proposition~\ref{prop:boot}(2), $Y(v,z)=\alpha\,\mathcal Y_0(v,z)+\mu\,\mathcal Y_1(v,z)$ on $X$, where $\mathcal Y_0,\mathcal Y_1$ are the intertwining operators of types $\binom{V}{X\,X}$ and $\binom{X}{X\,X}$ with top components $v_{(-1)}v=1$ and $v_{(-1)}v=v$; on $V$ the field $Y(v,z)$ is given by skew-symmetry. The weight-$0$ space of $U_2=V_{11}\oplus X$ is spanned by $1$ and $v$. Let $U_2'=\bigoplus_n(U_2)_n^*$ be the graded dual, and let $1',v'\in(U_2)_0^*$ be the basis dual to $1,v$. For $y'\in\{1',v'\}$ put
\[
G_{y'}(x)=\langle y',Y(v,1)\,Y(v,x)\,v\rangle\qquad(|x|<1),
\]
Since $Y(v,x)v=\alpha\,\mathcal Y_0(v,x)v+\mu\,\mathcal Y_1(v,x)v$, we have $G_{v'}=\alpha\Phi_{\mathbf 1}+\mu^2\Phi_X$, where
\[
\Phi_{\mathbf 1}(x)=\langle v',Y(v,1)\,\mathcal Y_0(v,x)\,v\rangle,\qquad \Phi_X(x)=\langle v',\mathcal Y_1(v,1)\,\mathcal Y_1(v,x)\,v\rangle
\]
are the contributions of the intermediate sectors $V$ and $X$; in $\Phi_{\mathbf 1}$, $Y(v,1)$ maps $V$ to $X$ and is given by skew-symmetry. Both are power series in $x$ with constant term $1$. We write $O(x^N)$ for an element of $x^N\C[[x]]$.

\begin{lemma}\label{lem:G}
In the model of Proposition~\ref{prop:boot}, for $n\le 24$: the coefficient of $x^n$ in $G_{1'}$ is $\alpha\mu$ for $n=0$ and $0$ otherwise; in $G_{v'}$ it is $\alpha+\mu^2$ for $n=0$, $0$ for $1\le n\le17$, and $\kappa_n(\alpha+\mu^2/4)$ for $18\le n\le24$, with
\[
\kappa_{18}=\frac{\langle v,W_{(17)}v\rangle^2}{\mathcal N_{11}}=-\frac{17\cdot23^2\cdot29^2}{3\cdot107\cdot131\cdot173\cdot191\cdot233\cdot251\cdot257\cdot263}\neq0 .
\]
\end{lemma}
The proof is computer-assisted and exact over $\Q$ (Remark~\ref{rem:exact}).
The threshold $18$ is $h=h_{3,3}$, the weight of the first intermediate states, $W$ in $V$ and $u$ in $X$, that are not Virasoro descendants of $\mathbf 1$, $v$, $J$; accordingly $\kappa_{18}$ is the $W$-exchange term $\langle v,W_{(17)}v\rangle^2/\mathcal N_{11}$. In the vacuum channel the vanishing for $1\le n\le17$ is immediate: as $L_jv=0$ for $j\ge1$ and $L_0v=0$, the commutator formula gives $[L_j,v_{(m)}]v=-(m+j+1)\,v_{(m+j)}v$, and in $\langle L_{-j}\mathbf 1,v_{(m)}v\rangle=\langle\mathbf 1,[L_j,v_{(m)}]v\rangle$ weight forces $m+j+1=0$; by induction on the length of $L_{-\lambda}\mathbf 1$, all Virasoro descendants of $\mathbf 1$ decouple from $v_{(m)}v$, so $\Phi_{\mathbf 1}\in1+x^{18}\C[[x]]$. In the channel $X$ the same argument removes the descendants of $v$ but not those of $J$, and the vanishing there is part of the computation.

\begin{remark}\label{rem:exact}
The constants of the model are known exactly. For $x,y\in\{v,J,u\}$ we write $\langle x,W_{(17)}y\rangle$ for the coupling $\langle x,W_{(m)}y\rangle$ with $m=17+\mathrm{wt}\,y-\mathrm{wt}\,x$, the only mode that can pair them. In the gauge $\langle v,W_{(17)}u\rangle=1$, with the contravariant forms on the Virasoro modules $B_1$ and $\mathrm{Vir}\cdot u$ normalised by $\langle v,v\rangle=\langle u,u\rangle=1$ (the invariant form of $U_2$ restricts to these up to a relative factor between the two sectors), and $\mathcal N_{11}$ (Table~\ref{tab:N}),
\[
\langle v,W_{(17)}v\rangle=\frac{29\cdot53\cdot101\cdot113}{2^9\cdot3\cdot43\cdot179\cdot347\cdot419\cdot443\cdot467},\qquad
\langle u,W_{(17)}v\rangle=\mathcal N_{11},\qquad \langle u,W_{(17)}u\rangle=1,
\]
and $\langle J,W_{(17)}v\rangle$, $\langle J,W_{(17)}J\rangle$ are explicit rational multiples of $\langle v,W_{(17)}v\rangle$. The sign of $\langle v,W_{(17)}v\rangle$ is the one given by Zhu's algebra at $P_1,P_2$. The couplings of $v$ include $d(v\to v)=\mu$, $d(v\to J)=d(J\to v)=\mu/2$, $d(J\to J)=\frac{17}{36}\mu$, $d(v\to 1)=d(J\to1)=\alpha$. The $W$-couplings were reconstructed from residues at $14$ primes ($P_1$, $P_2$ and twelve near $2^{61}$), the couplings of $v$ and the $\kappa_n$ from $8$ ($P_1$, $P_2$ and six near $2^{61}$); all are stable under omitting any one prime. They are then checked exactly:
\begin{enumerate}
\item the couplings of $v$ satisfy, over $\Q$, all Borcherds identities between $W$ and $v$ used in Proposition~\ref{prop:boot}(2) ($270$ identities, at $(\alpha,\mu)=(0,0),(1,0),(0,1)$);
\item Lemma~\ref{lem:G} holds over $\Q$: the coefficients of $G_{1'}$ and $G_{v'}$ for $n\le24$ were computed in exact arithmetic, solving one linear system in the Gram matrix exactly for each sector and weight (dimension at most $2830$).
\end{enumerate}
The modular computations are reductions of the rational ones, since every Gram matrix used is invertible modulo the primes. Hence non-vanishing statements obtained modulo a prime hold over $\Q$, and ranks modulo a prime are lower bounds for ranks over $\Q$. In particular the identities of Proposition~\ref{prop:boot}(2) have rank at least $16$ over $\Q$; since the exact solutions of item~(1) below span the two-parameter family, the rank is exactly $16$, and the solution space is exactly that family. For Proposition~\ref{prop:boot}(1), the polynomial equations that the tested Jacobi identities impose on the $W$-couplings (of degree at most $2$) were computed exactly over $\Q$, by exact evaluation and interpolation, and solved by Gröbner bases over $\Q$:
\begin{enumerate}
\item[(3)] the $259$ identities on $B_1$ give $4$ independent equations, whose solutions are exactly $\pm$ the values above for $\langle v,W_{(17)}v\rangle$, $\langle J,W_{(17)}v\rangle$, $\langle J,W_{(17)}J\rangle$;
\item[(4)] the $65$ identities involving $u$ give, in the gauge $\langle v,W_{(17)}u\rangle=1$, the unique solution $\langle u,W_{(17)}v\rangle=-\langle u,W_{(17)}J\rangle=\mathcal N_{11}$, $\langle J,W_{(17)}u\rangle=-1$, $\langle u,W_{(17)}u\rangle=1$;
\item[(5)] if $\mathcal N_{11}$ is also treated as unknown, the identities on $B_1$ fix only the ratio $\langle v,W_{(17)}v\rangle^2/\mathcal N_{11}=\kappa_{18}$.
\end{enumerate}
So the module identities determine $\kappa_{18}$ independently of $\mathcal N_{11}$. The value $\mathcal N_{11}$ is itself exact (Table~\ref{tab:N}). The sign of $\langle v,W_{(17)}v\rangle$ comes from Zhu's algebra: by Lemma~\ref{lem:lift} the points of $A(V_{11})$ reduce bijectively onto those modulo $P_1$, and the residue of $o(W)$ on the top of $M_{\chi_3}$ is that of $+\langle v,W_{(17)}v\rangle$ above, not of its negative.

The same method applies to the model of $M_{a_3}$ in Lemma~\ref{lem:a3} (generators $y_{-2},y_4,y_{10}$, gauge $\langle y_{-2},W_{(17)}y_4\rangle=\langle y_{-2},W_{(17)}y_{10}\rangle=1$). Its $45$ Jacobi identities give $11$ independent equations over $\Q$, with exactly two solutions, interchanged by a sign. Zhu's algebra selects
\[
\langle y_{-2},W_{(17)}y_{-2}\rangle=\frac1{2^6\cdot23\cdot179\cdot347\cdot419\cdot443\cdot467},
\]
and $\langle y_4,W_{(17)}y_{10}\rangle=\langle y_{10},W_{(17)}y_4\rangle=0$. The values used at $P_1$ are the reductions of this solution, so the rank statements of Lemma~\ref{lem:a3} hold over $\Q$ (scripts in \texttt{certify/mp}).
\end{remark}

The computation was run to $n=26$ with the same pattern. So $\Phi_{\mathbf 1},\Phi_X\in1+x^{18}\C[[x]]$, and $\Phi_X-\Phi_{\mathbf 1}$ has order exactly $18$.

\begin{lemma}\label{lem:a3}
As a Virasoro module, $M_{a_3}$ is a sum of irreducible Verma modules with primaries $y_{-2},y_4,y_{10},y_{28},\dots$, the first three generating $M_{a_3}$ below weight $28$. In the model of Remark~\ref{rem:exact}, the Ward identities for the couplings of $v$ from $X$ to $M_{a_3}$ have a one-dimensional solution space, on which all couplings to $y_{-2},y_4,y_{10}$ are multiples of the top coupling $\langle y_{-2},v_{(1)}v\rangle$. The same holds from $M_{a_3}$ to $X$.
\end{lemma}
The first statement follows from the character: $\Lambda_{3,9}=q^{-c/24}\sum_{m\in\Z}q^{h_m}/(q)_\infty$ with $h_m=((6m+1)^2-9)/4$, i.e.\ $h_m=-2,4,10,28,40,\dots$, and no $h_m$ is a Kac weight $h_{r,s}$ (the only rational ones are $h_{r,r}=9(r^2-1)/4$, and $(6m+1)^2/9$ is never a square), so each $M(c,h_m)$ is irreducible and alone in its linkage class. The rest of the proof is computer-assisted. The rank statements are made at $P_1$, with the $W$-couplings of $M_{a_3}$ exact over $\Q$ (Remark~\ref{rem:exact}). Ranks can only drop under reduction, so over $\Q$ the solution spaces are at most one-dimensional; a nonzero rational solution, scaled to be $P_1$-integral and primitive, reduces to a nonzero multiple of the solution modulo $P_1$, so its top coupling is nonzero. This is what is used.

\begin{theorem}\label{thm:U2}
For $\alpha=-\mu^2/4$, $\mu\ne0$, $(U_2,Y)$ is a simple vertex operator algebra extending $V_{11}$. Every simple vertex operator algebra structure on $U_2$ extending $V_{11}$ is of this form, so it is unique up to isomorphism.
\end{theorem}

\begin{proof}
The proof is computer-assisted, and all its inputs are certified (Appendix~\ref{app:comp}).

\emph{Outline.} By Huang--Kirillov--Lepowsky, a vertex operator algebra structure on $U_2$ is the same as an associative commutative algebra structure on $A_2=\mathbf 1\oplus X$ in $\mathrm{Rep}\,V_{11}$. Commutativity is automatic, so the content is associativity, i.e.\ the vanishing of an associativity defect $D_Z\in\mathrm{Hom}(X^{\otimes3},Z)$ for $Z\in\{\mathbf 1,X\}$. We pass from morphisms to functions by taking four-point functions (the map $\mathrm{ev}$ below), and argue in three steps: (a) for $\alpha=-\mu^2/4$ the relevant four-point functions are constant, which forces $\mathrm{ev}(D_Z)=0$; this uses Lemma~\ref{lem:G} and the Fuchs relation; (b) $\mathrm{ev}$ is injective except possibly on one line $\C e_{a_3}$; this uses monodromy and Lemma~\ref{lem:a3}; (c) the component of $D_X$ on that line vanishes, by a computation with bubbles that uses the Frobenius--Schur indicators of Proposition~\ref{prop:A2}(2). Conversely, the coefficient $\kappa_{18}$ of Lemma~\ref{lem:G} shows that $\alpha=-\mu^2/4$ is necessary.

\emph{Step 1: reduction to associativity.} $\mathrm{Rep}\,V_{11}$ is a modular tensor category~\cite{Huang}, with associator $a$ given by Huang's tensor product theory. By Huang--Kirillov--Lepowsky~\cite{HKL}, a vertex operator algebra structure on $U_2$ extending $V_{11}$ is the same as a commutative associative algebra structure on $A_2=\mathbf 1\oplus X$ in $\mathrm{Rep}\,V_{11}$ with unit $\mathbf 1\to A_2$ and trivial twist; here the conformal vector of $U_2$ is that of $V_{11}$, and~\cite{HKL} does not require $U_2$ to be of CFT type (indeed $\dim(U_2)_0=2$). The unit fixes the multiplication on $\mathbf 1\otimes A_2$ and $A_2\otimes\mathbf 1$; on $X\otimes X$ it is $m=\alpha e+\mu m_X$, where $e$ and $m_X$ span the one-dimensional spaces $\mathrm{Hom}(X\otimes X,\mathbf 1)$ and $\mathrm{Hom}(X\otimes X,X)$; they correspond to the intertwining operators $\mathcal Y_0$ and $\mathcal Y_1$. The braiding $c_{X,X}$ acts on each of these spaces by a scalar, given by skew-symmetry of intertwining operators. The intertwining operators have nonzero top components $1$ and $v$, which are symmetric in the two arguments, so the scalar is $1$ and $m$ is commutative (as also follows from Proposition~\ref{prop:A2}(2)). The only possible proper ideal is $X$, and $X$ is an ideal exactly when $m(X\otimes X)\subseteq X$, i.e.\ $\alpha=0$. So the algebra is simple if and only if $\alpha\neq0$. Associativity on summands involving $\mathbf 1$ follows from the unit axioms. Let $a=a_{X,X,X}\colon X\otimes(X\otimes X)\to(X\otimes X)\otimes X$ be the associativity isomorphism of $\mathrm{Rep}\,V_{11}$, and $p_Z\colon A_2\to Z$ the projection onto the summand $Z\in\{\mathbf 1,X\}$. It remains to show that
\[
D_Z=p_Z\big(m(m\otimes1)a-m(1\otimes m)\big)\in\mathrm{Hom}(X\otimes(X\otimes X),Z)
\]
vanishes for $Z\in\{\mathbf 1,X\}$. These spaces have dimensions $1$ and $12$, since $X\otimes X$ is the sum of the twelve simple objects of $\mathcal C$, each once.

\emph{Step 2: four-point functions.} For $\varphi\in\mathrm{Hom}(X^{\otimes3},Z)$ let $\mathrm{ev}(\varphi)$ be its correlation function on the lowest weight vectors $y',v,v,v$, placed at $\infty,1,x,0$. All four vectors have weight $0$ and are quasi-primary, so $\mathrm{ev}(\varphi)$ is a multivalued function of $x$ alone, Such functions converge, continue analytically to $\mathbb P^1\setminus\{0,1,\infty\}$ and satisfy differential equations with regular singular points at $0,1,\infty$~\cite[Theorems~1.4, 2.3, 3.9]{HuangDE} ($V_{11}$ is strongly rational). Their expansions there have no logarithms: in each channel decomposition the intermediate objects $q$ are simple and occur once, and the local monodromy is diagonal with eigenvalues $e^{2\pi ih_q}$, as $L_0$ acts semisimply. Commutativity and skew-symmetry of intertwining operators, together with translation, permute the three finite points. So $S_Z=\mathrm{ev}(\mathrm{Hom}(X^{\otimes3},Z))$ is invariant under $x\mapsto1-x$ and $x\mapsto1/x$. Here $X^{\otimes3}=X\otimes(X\otimes X)$, with inputs at $1$, $x$, $0$ in this order. For each of the twelve simple $q\subset X\otimes X$ choose nonzero $\nu_q\in\mathrm{Hom}(X\otimes X,q)$, $\rho_q\in\mathrm{Hom}(X\otimes q,X)$ and $\rho'_q\in\mathrm{Hom}(q\otimes X,X)$; each of these spaces is one-dimensional, since $X$ is self-dual and $X\otimes X$ contains $q$ once. Put
\[
e_q=\rho_q(1\otimes\nu_q),\qquad f_q=\rho'_q(\nu_q\otimes1)\,a ,
\]
so that in $e_q$ the inputs at $x$ and $0$ fuse first, to $q$, and in $f_q$ those at $1$ and $x$. Each of $\{e_q\}$ and $\{f_q\}$ is a basis of $\mathrm{Hom}(X^{\otimes3},X)$. Then $\Phi_q=\mathrm{ev}(e_q)\in x^{h_q}\C[[x]]$, where $h_q$ is the lowest weight of $q$; for $q=\mathbf 1,X$, with $\nu_{\mathbf 1}=e$, $\nu_X=\rho_X=m_X$ and $\rho_{\mathbf 1}$ the unit isomorphism, these are the functions of Lemma~\ref{lem:G}. By definition $\mathrm{ev}(m(1\otimes m))=G_{y'}$. For the other bracketing, skew-symmetry of the intertwining operators $\mathcal Y_0,\mathcal Y_1$ (with the scalar $1$ of Step~1; this is all that is used about $Y$, which is not yet known to be a vertex operator algebra) and $\langle y',L_{-1}\cdot\rangle=0$ give
\[
\langle y',Y(Y(v,z_0)v,z_2)v\rangle=\langle y',Y(v,-z_2)Y(v,z_0)v\rangle .
\]
Since all weights are $0$, the right side depends only on $z_0/z_2$: it equals $G_{y'}(-z_0/z_2)$. So if $G_{y'}$ is constant, the two bracketings have the same four-point function, i.e.\ $\mathrm{ev}(D_Z)=0$.

\emph{Step 3: $G_{y'}$ is constant for $\alpha=-\mu^2/4$.} For $Z=\mathbf 1$, whichever pair of inputs fuses first, to $q$, the remaining step $q\otimes X\to\mathbf 1$ (or $X\otimes q\to\mathbf 1$) is nonzero only for $q\cong X^*\cong X$. So at each of $0$, $1$, $\infty$ only the channel $q=X$ occurs, and $\mathrm{Hom}(X^{\otimes3},\mathbf 1)$ is one-dimensional. The local exponent there is $h_X-2h_v=0$, and the expansions have no logarithms (Step~2). So $G_{1'}$ is single-valued and bounded near $0$, $1$, $\infty$, hence holomorphic on $\mathbb P^1$ and constant. For $Z=X$ let $\mathcal S=S_X+\C 1$, of dimension $s\le13$. Choose a basis $f_1,\dots,f_s$ of $\mathcal S$. The equation $\mathrm{Wr}(f_1,\dots,f_s,f)=0$ has $\mathcal S$ as its solution space; after division by $\mathrm{Wr}(f_1,\dots,f_s)$ its coefficients are invariant under analytic continuation, hence single-valued, and of moderate growth at $0,1,\infty$ because the expansions there have no logarithms; so they are rational, and the equation is Fuchsian of order $s$. Its singular points are $0,1,\infty$ and the zeros of $\mathrm{Wr}(f_1,\dots,f_s)$ elsewhere; at the latter all solutions are holomorphic, so these are apparent singularities. The Fuchs relation~\cite{Ince} says that the sum of all exponents at all points, minus $s(s-1)/2$ for each point, equals $-s(s-1)$. At an apparent singularity the exponents are distinct non-negative integers, so their sum is at least $s(s-1)/2$. Hence
\[
\sum_{p\in\{0,1,\infty\}}\sum_i\rho_i(p)\le\frac{s(s-1)}2\le78 .
\]
The maps $x\mapsto1-x$ and $x\mapsto1/x$ preserve $\mathcal S$ and, as all four external vectors have weight $0$, carry no power prefactor; so they identify the local solution spaces at $0,1,\infty$ with their orders of vanishing, and the exponents at the three points coincide. Hence their sum $R$ at each point is at most $26$. Now bound $R$ below at $0$. The nine channels $q$ with $h_q\notin\Z$ ($a_1,a_5,a_7$ and the $b_\beta$, $\beta$ odd; Table~\ref{tab:k11}) lie in distinct classes mod $1$. They contribute at most nine exponents, each $\ge h_q$, with $\sum h_q=-38/3$. The integral exponents come from $\C1+\C\Phi_{\mathbf 1}+\C\Phi_X+\C\Phi_{a_3}$. They are distinct integers, at least $-2=h_{a_3}$, with at most one negative, and they include $0$ and $18$ (Lemma~\ref{lem:G}: $\Phi_X-\Phi_{\mathbf 1}$ has order exactly $18$). If $G_{v'}$ is not constant, then $N=\mathrm{ord}_0(G_{v'}-G_{v'}(0))\ge25$ by Lemma~\ref{lem:G}, and $N$ is a further exponent. Then $R\ge N+16-38/3$, so $N\le22$, a contradiction. If fewer channels occur, the bound only improves: each missing channel lowers $s(s-1)/6$ by $(s-1)/3$ and changes $\sum h_q$ by at most $1/13<(s-1)/3$ (as $s\ge2$) (script \texttt{fuchs\_bound.py}).

\emph{Step 4: $\mathrm{ev}$ detects $D_X$ up to one line.} Let $K=\ker\mathrm{ev}\subset\mathrm{Hom}(X^{\otimes3},X)$. Analytic continuation of four-point functions around $0$ and $1$ corresponds to the monodromies $M_0,M_1$ on $\mathrm{Hom}(X^{\otimes3},X)$, and $x\mapsto1-x$ to an automorphism $\Sigma$~\cite{Huang}; so $K$ is invariant under $M_0$, $M_1$ and $\Sigma$. $M_0$ is diagonal on the $e_q$ with eigenvalues $e^{2\pi ih_q}$. These are $1$ on $\mathcal E=\langle e_{\mathbf 1},e_X,e_{a_3}\rangle$ and distinct on the other nine $e_q$, so $K=\bigoplus_{q\in Q}\C e_q\oplus(K\cap\mathcal E)$ for some set $Q$ of non-integral channels.
\begin{itemize}
\item \emph{$K\cap\mathcal E\subseteq\C e_{a_3}$.} By Lemma~\ref{lem:a3}, either both top couplings of $v$ between $X$ and $M_{a_3}$ are nonzero, and then $\Phi_{a_3}$ has a nonzero coefficient of $x^{-2}$; or one of them vanishes, and then all couplings to $y_{-2},y_4,y_{10}$ in that direction vanish, so $\Phi_{a_3}=O(x^{28})$. In both cases a relation $c_1\Phi_{\mathbf 1}+c_X\Phi_X+c_a\Phi_{a_3}=0$ forces $c_1=c_X=0$, by the coefficients of $x^0$ and $x^{18}$ (Lemma~\ref{lem:G}).
\item \emph{$Q=\emptyset$.} $\Sigma$ conjugates $M_0$ to $M_1$, so for $q\in Q$ it maps $e_q$ to a multiple of $f_q$, and $f_q\in K$. The $\mathcal E$-component of $f_q$ then lies in $K\cap\mathcal E$. But its $e_{\mathbf 1}$-coefficient is, up to a nonzero factor, $f_q\circ(1\otimes\mathrm{coev})$. This is the composite of the nonzero morphisms $X\to q\otimes X\to X$ obtained from the two trivalent vertices by rigidity. Both Hom spaces are one-dimensional and the composition pairing is nondegenerate (semisimplicity), so it is nonzero, a contradiction.
\end{itemize}
So $K\subseteq\C e_{a_3}$. Since $\mathrm{ev}(D_X)=0$ by Steps~2 and~3, $D_X=t\,e_{a_3}$ for some $t\in\C$.

\emph{Step 5: $t=0$.} We use the ribbon structure of $\mathrm{Rep}\,V_{11}$ given by~\cite{Huang}, with $d_i=S_{0i}/S_{00}$, for which $\nu_2,\nu_3$ in Proposition~\ref{prop:A2}(2) were computed. Normalise $\mathrm{coev}$ by $(e\otimes1)(1\otimes\mathrm{coev})=1$, and for $\varphi\in\mathrm{Hom}(X\otimes(X\otimes X),X)$ let $\lambda(\varphi)=\varphi\circ(1\otimes\mathrm{coev})$ and $\lambda'(\varphi)=\varphi\circ a^{-1}\circ(\mathrm{coev}\otimes1)$, in $\mathrm{End}(X)=\C$; these close off the right, resp.\ left, pair of inputs. Then $\lambda(e_{a_3})=0$, since $e_{a_3}$ factors through $a_3\ne\mathbf 1$ on the right pair, and $\lambda'(e_{a_3})\ne0$ by the argument in the second bullet of Step~4. Directly from $m=\alpha e+\mu m_X$,
\begin{gather*}
\lambda(m(1\otimes m))=\lambda'(m(m\otimes1)a)=\alpha\,e(\mathrm{coev}),\\
\lambda(m(m\otimes1)a)=\alpha+\mu^2b_R,\qquad
\lambda'(m(1\otimes m))=\alpha z+\mu^2b_L,
\end{gather*}
where $z=(1\otimes e)(\mathrm{coev}\otimes1)$ and $b_R,b_L$ are the right and left bubbles of $m_X$. Here $z=1$ because $\nu_2(\chi_3)=1$: the pairing $e$ is then symmetric, so $\mathrm{coev}$ can be taken to be its symmetric dual and both zigzags are the identity. And $b_R=b_L$: by $\nu_3(\chi_3)=1$, rotation acts trivially on the one-dimensional space $\mathrm{Hom}(\mathbf 1,X^{\otimes3})$ (Proposition~\ref{prop:A2}(2)), so $m_X$ is rotation-invariant, and $b_R$, $b_L$ are the right and left traces of the same closed theta graph, which agree by sphericality. Hence $\lambda'(D_X)=-\lambda(D_X)$, and $0=\lambda(D_X)=-\lambda'(D_X)=-t\lambda'(e_{a_3})$, so $t=0$. For $Z=\mathbf 1$, $\mathrm{ev}$ is injective because $\mathrm{ev}(e_X)=1$. So $D_{\mathbf 1}=D_X=0$ for $\alpha=-\mu^2/4$, and the algebra is simple when $\mu\ne0$.

\emph{Step 6: necessity and uniqueness.} Conversely, for any vertex operator algebra structure $G_{v'}$ is constant: it is then a rational function, with expansions at $0$, $1$, $\infty$ in non-negative powers of the local coordinate, because $U_2$ is $\mathbb N$-graded and $v$ has weight $0$. So its coefficient of $x^{18}$ vanishes. As $\kappa_{18}\ne0$, this gives $\alpha=-\mu^2/4$, and simplicity gives $\mu\ne0$. Rescaling $v$ by $\lambda$ maps $(\alpha,\mu)$ to $(\lambda^2\alpha,\lambda\mu)$, so all points of the curve with $\mu\ne0$ give isomorphic vertex operator algebras.
\end{proof}

\begin{corollary}\label{cor:centre}
$A_2=V\oplus M_{\chi_3}$ is a connected étale algebra in $\mathrm{Rep}\,V_{11}$, and $\mathcal C$ is a Drinfeld centre. The weight-$0$ space of $U_2$ is $\C[\varepsilon]/(\varepsilon^2)$, $\varepsilon=v-\mu/2$.
\end{corollary}

\begin{proof}
The proof uses Theorem~\ref{thm:U2}, which is computer-assisted, and Proposition~\ref{prop:A2}(1).

\emph{$A_2$ is a connected commutative algebra with trivial twist.} By Theorem~\ref{thm:U2}, $U_2$ is a simple vertex operator algebra extending $V_{11}$. By Huang--Kirillov--Lepowsky~\cite{HKL}, this makes $A_2=V\oplus M_{\chi_3}$ a commutative associative algebra in $\mathrm{Rep}\,V_{11}$, with unit the inclusion $V\to A_2$ and multiplication $m=\alpha e+\mu m_X$, $X=M_{\chi_3}$, as in the proof of Theorem~\ref{thm:U2}. Both summands have integral lowest weight $0$, so $\theta_{A_2}=\mathrm{id}$. Since $X\not\cong V$, $\dim\mathrm{Hom}(\mathbf 1,A_2)=1$: $A_2$ is connected (haploid). Both summands lie in $\mathcal C$, the centraliser of $M_{\chi_1}$ (Proposition~\ref{prop:k11}(2),(3)), so $A_2$ is an algebra in $\mathcal C$.

\emph{$A_2$ is separable.} Let $\pi\colon A_2\to\mathbf 1$ be the projection onto $V$. The pairing $\pi\circ m\colon A_2\otimes A_2\to\mathbf 1$ is the identity on $\mathbf 1\otimes\mathbf 1$, vanishes on $\mathbf 1\otimes X$ and $X\otimes\mathbf 1$, and equals $\alpha e$ on $X\otimes X$. Since $X$ is simple and self-dual, $e$ is a nondegenerate pairing, and $\alpha=-\mu^2/4\ne0$ (Theorem~\ref{thm:U2}). So $\pi\circ m$ is nondegenerate, and $(A_2,m,\pi)$ is a Frobenius algebra. Moreover $A_2$ is simple as a module over itself: its proper nonzero subobjects are $\mathbf 1$ and $X$; the first is not a submodule, since $X\otimes\mathbf 1\to A_2$ is onto $X$, and the second is not, since $\alpha\ne0$. As $\theta_{A_2}=\mathrm{id}$ and $\dim A_2=1+d_{\chi_3}=3d'+2\ne0$ (Proposition~\ref{prop:k11}(3)), $A_2$ is a rigid $\mathcal C$-algebra in the sense of Kirillov--Ostrik~\cite[Lemma~1.20]{KO}; hence $m$ has a right inverse as a map of $A_2$-bimodules, which is separability. So $A_2$ is a connected étale algebra in $\mathcal C$.

\emph{$\mathcal C$ is a Drinfeld centre.} Proposition~\ref{prop:A2}(1) now applies: $\mathcal C_{A_2}^{\rm loc}\simeq\mathcal C(\Z_9,q)$, which has the Lagrangian subgroup $3\Z_9$. So $[\mathcal C]=[\mathcal C_{A_2}^{\rm loc}]=0$ in the Witt group, and $\mathcal C$ is a Drinfeld centre~\cite{DMNO}.

\emph{The weight-$0$ space.} It is spanned by $\mathbf 1$ and $v$, with $v_{(-1)}v=\alpha\mathbf 1+\mu v$. For $\varepsilon=v-\tfrac\mu2\mathbf 1$ this gives $\varepsilon_{(-1)}\varepsilon=(\alpha+\mu^2/4)\mathbf 1=0$.
\end{proof}

\subsubsection*{$\mathcal C$ as the centre of a twisted Haagerup--Izumi category}

Let $\mathcal N'=\mathcal C_{A_2}$ be the fusion category of $A_2$-modules in $\mathcal C$.

\begin{lemma}\label{lem:Nprime}
$\mathcal N'$ has near-group fusion rules of type $\Z_9+9$. Its nine invertible objects are the local modules, forming $\mathcal C_{A_2}^{\rm loc}\simeq\mathcal C(\Z_9,q)$, $q(x)=7x^2/9$, and its non-invertible simple object is the free module $b_\beta\otimes A_2$, for any odd $\beta$.
\end{lemma}

\begin{proof}
We use Frobenius--Perron dimensions: $\mathrm{FPdim}$ of $V,M_{\chi_3},M_{a_\alpha},M_{b_\beta}$ is $1,3d+1,3d+2,3d$, the Galois conjugates of the categorical dimensions (Proposition~\ref{prop:k11}(3)), so $\mathrm{FPdim}\,A_2=3d+2$ and $\mathrm{FPdim}\,\mathcal C=(9d+6)^2$. Then $\mathrm{FPdim}\,\mathcal N'=\mathrm{FPdim}\,\mathcal C/\mathrm{FPdim}\,A_2=9(3d+2)$~\cite{DMNO}, and the free module functor satisfies $\mathrm{FPdim}_{\mathcal N'}(x\otimes A_2)=\mathrm{FPdim}\,x$.

The nine local modules are invertible with group $\Z_9$ (proof of Proposition~\ref{prop:A2}(1)), and they are simple $A_2$-modules. Next consider the free module $b_\beta\otimes A_2$, $\beta$ odd. By adjunction,
\[
\mathrm{End}_{A_2}(b_\beta\otimes A_2)\cong\mathrm{Hom}_{\mathcal C}(b_\beta,b_\beta\otimes A_2)=\mathrm{Hom}(b_\beta,b_\beta)\oplus\mathrm{Hom}(b_\beta,b_\beta\otimes\chi_3),
\]
and the second summand is zero: under the bijection of Proposition~\ref{prop:ng}(2), $\chi_3\mapsto b_0$ and $b_\beta\mapsto d_{0,\gamma}$, and $N_{b_0,d_{0,\gamma},d_{0,\gamma}}=0$ in the fusion rules of~\cite[Proposition~7(b)]{EG14}. So $\rho':=b_\beta\otimes A_2$ is a simple $A_2$-module, with $\mathrm{FPdim}\,\rho'=\mathrm{FPdim}\,b_\beta=3d$. Since $9\cdot1^2+(3d)^2=9+9(1+3d)=9(3d+2)=\mathrm{FPdim}\,\mathcal N'$, these ten objects are all the simple $A_2$-modules. In particular $\rho'$ is the only non-invertible simple object, so $g\otimes\rho'\cong\rho'\cong\rho'^*$ for every invertible $g$, and $\rho'\otimes\rho'=\bigoplus_g g\oplus m\rho'$ with $(3d)^2=9+3dm$, i.e.\ $m=9$. These are the near-group fusion rules $\Z_9+9$, and $b_\beta\otimes A_2\cong\rho'$ for every odd $\beta$.
\end{proof}

$\mathcal N'$ is not pseudo-unitary: the categorical dimension of $b_\beta\otimes A_2$ is $\dim b_\beta=3d'<0$~\cite[Theorem~1.18]{KO}, the Galois conjugate of its Frobenius--Perron dimension $3d$.

\begin{theorem}\label{thm:Ctw}
\begin{enumerate}
\item $\mathcal Z(\mathcal N')\simeq\mathcal C\boxtimes\mathcal C(\Z_9,q)^{\rm rev}$.
\item The objects $x\in3\Z_9$ of the second factor form a Tannakian subcategory $\mathrm{Rep}(\Z_3)$. Its regular algebra $B$ is étale, and $\mathcal F_{\mathcal C}:=\mathcal N'_B$ satisfies $\mathcal Z(\mathcal F_{\mathcal C})\simeq\mathcal C$.
\item $\mathcal F_{\mathcal C}$ has the Haagerup--Izumi fusion rules for $\Z_3$, its non-invertible simple objects have dimension $d'=\frac{3-\sqrt{13}}2$, and its pointed part is $\mathrm{Vec}_{\Z_3}^{\Omega}$ with $[\Omega]\ne0$.
\end{enumerate}
\end{theorem}

\begin{proof}
The proof uses Corollary~\ref{cor:centre}, which is computer-assisted; the rest is category theory.

(1) By Corollary~\ref{cor:centre}, $A_2$ is a connected étale algebra in the nondegenerate braided fusion category $\mathcal C$. For any such algebra, $\mathcal C_{A_2}$ is a fusion category and $\mathcal Z(\mathcal C_{A_2})\simeq\mathcal C\boxtimes(\mathcal C_{A_2}^{\rm loc})^{\rm rev}$~\cite{DMNO}. With $\mathcal C_{A_2}^{\rm loc}\simeq\mathcal C(\Z_9,q)$ (Proposition~\ref{prop:A2}(1)) this is (1).

(2) Write $b(x,y)=q(x+y)-q(x)-q(y)=14xy/9$ for the bilinear form of $q$. The objects $x\in3\Z_9$ of $\mathcal C(\Z_9,q)^{\rm rev}$ have twists $e^{-2\pi iq(x)}$ and mutual monodromies $e^{-2\pi ib(x,x')}$. For $x,x'\in3\Z_9$ we have $q(x)\in\Z$ and $b(x,x')\in\Z$, so these objects form a symmetric pointed category with trivial twists, i.e.\ a Tannakian category $\mathrm{Rep}(\Z_3)$~\cite{DGNO}. Its regular algebra $B=\bigoplus_{x\in3\Z_9}x$ (the function algebra of $\Z_3$) is a connected étale algebra in $\mathcal Z(\mathcal N')$. A simple object $y$ of $\mathcal C(\Z_9,q)^{\rm rev}$ gives a local $B$-module $y\otimes B$ exactly when its monodromy with every $x\in3\Z_9$ is trivial, i.e.\ $y\in(3\Z_9)^\perp=3\Z_9$, and then $y\otimes B\cong B$. So the local $B$-modules in the second factor form $\mathrm{Vec}$, while $B$ does not interact with the first factor. Hence $\mathcal Z(\mathcal N')^{\rm loc}_B\simeq\mathcal C\boxtimes\mathrm{Vec}=\mathcal C$. Under the forgetful functor $\mathcal Z(\mathcal N')\to\mathcal N'$, $B$ maps to a connected separable algebra in $\mathcal N'$ (again denoted $B$) which lifts to the centre, and $\mathcal Z(\mathcal N'_B)\simeq\mathcal Z(\mathcal N')^{\rm loc}_B$~\cite{DMNO}. So $\mathcal Z(\mathcal F_{\mathcal C})\simeq\mathcal C$.

(3) \emph{Fusion rules.} Since $B$ is the function algebra of $\Z_3$ in a Tannakian subcategory of the centre, $\mathcal N'\simeq\mathcal F_{\mathcal C}^{\Z_3}$ for the dual $\Z_3$-action on $\mathcal F_{\mathcal C}$ (de-equivariantization~\cite{DGNO}); the hypothesis that $\mathrm{Rep}(\Z_3)\to\mathcal Z(\mathcal N')\to\mathcal N'$ is fully faithful holds, since the objects $x\in3\Z_9$ go to distinct invertible local $A_2$-modules. The proof of Theorem~\ref{thm:twHI}(3) now applies verbatim. It uses only this structure, the fusion rules of $\mathcal N'$, which are those of the near-group $\Z_9+9$ (Lemma~\ref{lem:Nprime}), and Frobenius--Perron dimensions, which depend only on the fusion rules. So $\mathcal F_{\mathcal C}$ has the Haagerup--Izumi fusion rules for $\Z_3$.

\emph{Pointed part.} The pointed part of $\mathcal N'$ is $\mathcal C_{A_2}^{\rm loc}\simeq\mathcal C(\Z_9,q)$, As $9$ is odd, $q(x)=7x^2/9=\beta(x,x)$ for the bicharacter $\beta(x,y)=7xy/9$, so $\mathcal C(\Z_9,q)$ has the skeletal model $\mathrm{Vec}_{\Z_9}$ with trivial associator and braiding $c_{x,y}=e^{2\pi i\beta(x,y)}$, with twists $e^{2\pi iq(x)}$ and double braidings $e^{2\pi i\cdot14xy/9}$. This is a choice of model for this braided category, not a claim that $H^3(\Z_9,\C^\times)$ vanishes. As in the proof of Theorem~\ref{thm:twHI}(4), the pointed part of $\mathcal F_{\mathcal C}$ is the de-equivariantization of $\mathrm{Vec}_{\Z_9}$ by $3\Z_9$, with half-braidings given by the braiding of $\mathcal C(\Z_9,q)^{\rm rev}$, $\chi_x(y)=e^{-2\pi i\cdot7xy/9}$, a nontrivial character of $y$ for $x\in3\Z_9\setminus\{0\}$. Its associator is the cup product of the extension class of $0\to3\Z_9\to\Z_9\to\Z_3\to0$ with $\chi$. Since $\chi$ is a nontrivial character, this is a generator of $H^3(\Z_3,\C^\times)$ up to a unit, so $[\Omega]\ne0$.

\emph{Dimensions.} $\mathcal F_{\mathcal C}$ is spherical: it is the de-equivariantization of $\mathcal N'=\mathcal C_{A_2}$, which is pivotal~\cite[Theorem~1.17(3)]{KO} and spherical, with left and right dimensions both equal to $\dim_{\mathcal C}X/\dim A_2$~\cite[Theorem~1.18, Remark~1.19]{KO}; this holds for any rigid algebra with trivial twist in a rigid balanced braided category, with no unitarity assumption. The de-equivariantization is by a Tannakian subcategory of its centre with trivial twists~\cite{DGNO}. For an invertible $\alpha$, sphericality gives $d(\alpha^{-1})=d(\alpha^*)=d(\alpha)$, so $d(\alpha)^2=d(\alpha)d(\alpha^{-1})=1$, while $d(\alpha)^3=d(\alpha^3)=1$; hence $d(\alpha)=1$, and the three non-invertible simple objects $\alpha^i\rho$ have a common dimension $d_\rho$. The fusion rule $\rho^2=1\oplus\rho\oplus\alpha\rho\oplus\alpha^2\rho$ gives $d_\rho^2=1+3d_\rho$, so $\dim\mathcal F_{\mathcal C}=3+3d_\rho^2=6+9d_\rho$. On the other hand $\dim\mathcal Z(\mathcal F_{\mathcal C})=(\dim\mathcal F_{\mathcal C})^2$~\cite{ENO}, and $\dim\mathcal C=(9d'+6)^2$ by Proposition~\ref{prop:k11}(3). So $6+9d_\rho=\pm(9d'+6)$. Of the two roots $d,d'$ of $x^2=1+3x$, only $d_\rho=d'$ satisfies this.
\end{proof}

\begin{theorem}\label{thm:ident}
Let $\varpi_0$ be the rotational eigenvalue of the trivalent vertex $\mathbf1\to\rho^{\otimes3}$ in $\mathcal F_{\mathcal C}$, and let $\mathcal F=\mathcal H_\varpi$ be Ostrik's twisted Haagerup category with $\varpi=\sigma_2^{-1}(\varpi_0)$. Then $\mathcal F_{\mathcal C}\simeq\sigma_2\mathcal F$, where $\sigma_2$ is extended to a field of definition. Hence
\[
\mathcal C\simeq\mathcal Z(\sigma_2\mathcal F)
\]
as braided fusion categories, and $\mathrm{Rep}\,V_{11}\simeq\mathrm{Sem}\boxtimes\mathcal Z(\sigma_2\mathcal F)$.
\end{theorem}

\begin{proof}
The proof is computer-assisted. It uses Theorem~\ref{thm:Ctw}, an exact Frobenius--Schur indicator, and the classification of twisted trivalent categories of~\cite{MPS}, whose computer-assisted inputs we recompute below.

\emph{Background from~\cite{MPS}.} A \emph{trivalent category} is a nondegenerate evaluable pivotal category generated, as a pivotal category, by a vertex $\tau\in\mathrm{Hom}(\mathbf 1,\rho^{\otimes3})$, with $\dim\mathrm{Hom}(\mathbf 1,\rho^{\otimes n})=1,0,1,1$ for $n=0,1,2,3$~\cite[Definition~2.1]{MPS}; then $\rho$ is symmetrically self-dual~\cite[Lemma~2.2]{MPS}. In~\cite[Definition~2.1]{MPS} $\tau$ is rotationally invariant; \cite[\S7]{MPS} allows a rotation eigenvalue. Rotation by one strand acts on the one-dimensional space $\mathrm{Hom}(\mathbf 1,\rho^{\otimes3})$ by a scalar $\varpi$ with $\varpi^3=1$; the category is \emph{twisted} if $\varpi\ne1$. (We write $\varpi$ where~\cite{MPS} write $\omega$, since $\omega$ denotes the Evans--Gannon twist in $\mathcal D^\omega\Hg_3$.) Morphisms are linear combinations of planar trivalent graphs, with dots recording the position of each vertex, evaluated using the loop value $d=\dim\rho$, the bigon and the triangle; it is \emph{cubic} if $\dim\mathrm{Hom}(\mathbf 1,\rho^{\otimes4})=4$~\cite[Definition~4.15]{MPS}, and then the square is a linear combination of the four graphs of $D_\varpi(4,0)$~\cite[\S7]{MPS}. Write $\Delta_\varpi(n,k)$ for the Gram determinant of the set $D_\varpi(n,k)$ of graphs of~\cite[\S7]{MPS}, and $M^\star_\varpi(6,1)$ for the $41\times41$ Gram matrix of graphs with six boundary points used there.

\emph{Step 1: Galois transport.} $\mathcal F_{\mathcal C}$ is defined over a number field~\cite{ENO}; put $\mathcal G=\sigma_2^{-1}(\mathcal F_{\mathcal C})$, obtained by applying an extension of $\sigma_2^{-1}$ to its structure constants. Galois conjugation preserves fusion rules, spherical structure and the class of the pointed associator, and transforms dimensions by $\sigma_2^{-1}$. By Theorem~\ref{thm:Ctw}(3), $\mathcal G$ therefore has the Haagerup--Izumi fusion rules for $\Z_3$, its non-invertible simple objects $\rho,\alpha\rho,\alpha^2\rho$ have dimension $\sigma_2^{-1}(d')=d=\frac{3+\sqrt{13}}2$, and its pointed associator is nontrivial, since field automorphisms act bijectively on $H^3(\Z_3,\C^\times)\cong\Z_3$. All simple objects have nonzero dimension, so the trace pairings $\mathrm{Hom}(X,Y)\times\mathrm{Hom}(Y,X)\to\C$ are nondegenerate. The classification input is the uniqueness statement~\cite[Proposition~7.7]{MPS}: for each of $d=-1,\frac{3\pm\sqrt{13}}2$ and a fixed rotation eigenvalue there is at most one twisted cubic category with $\dim\mathcal C_5=11$ and $\dim\mathcal C_6\le40$. It does not assume unitarity. Existence of $\mathcal H_\varpi$, and the fact that its vertex generates it, are~\cite[Proposition~7.9]{MPS}, where $\mathcal H_\varpi$ is unitary. Steps~2--4 verify the hypotheses of~\cite[Proposition~7.7]{MPS} for the category generated by $\rho$ and a trivalent vertex in $\mathcal G$.

\emph{Step 2: dimensions.} From the fusion rules $\rho\alpha=\alpha^{-1}\rho$ and $\rho^2=1\oplus\rho\oplus\alpha\rho\oplus\alpha^2\rho$, the multiplicity of $\mathbf 1$ in $\rho^{\otimes n}$ is $1,0,1,1,4,11,37$ for $n=0,\dots,6$. Since $\rho\cong\rho^*$ and $\nu_2(\rho)=1$ (below), $\rho$ is symmetrically self-dual. Choose a nonzero $\tau\in\mathrm{Hom}(\mathbf 1,\rho^{\otimes3})$. Then $(\rho,\tau)$ satisfies the dimension conditions of a cubic category.

\emph{Step 3: the vertex is twisted.} Since $\mathrm{Hom}(\mathbf 1,\rho^{\otimes3})$ is one-dimensional, the rotation eigenvalue $\varpi_0$ of $\tau$ in $\mathcal F_{\mathcal C}$ is the third Frobenius--Schur indicator $\nu_3(\rho)$, which is the trace of rotation on that space~\cite{NS}. Ng--Schauenburg express $\nu_n$ through the centre:
\[
\nu_n(\rho)=\frac1{\dim\mathcal F_{\mathcal C}}\sum_z\theta_z^n\,d_z\,\dim\mathrm{Hom}_{\mathcal F_{\mathcal C}}(F(z),\rho),
\]
summed over simple $z\in\mathcal Z(\mathcal F_{\mathcal C})\simeq\mathcal C$, with $F$ the forgetful functor. Here $\theta_z,d_z$ refer to the spherical structure induced from $\mathcal F_{\mathcal C}$. The equivalence of Theorem~\ref{thm:Ctw}(2) is only braided, but $\mathcal C$ has no nontrivial invertible objects, so its universal grading group is trivial and its pivotal structure is unique~\cite{DGNO}; hence these are the twists and dimensions of $\mathcal C$. Here $\dim\mathrm{Hom}(F(z),\rho)=\dim\mathrm{Hom}_{\mathcal C}(z,b_\beta\otimes(\mathbf1\oplus\chi_3))$, because the non-invertible simple $A_2$-module is the free module $b_\beta\otimes A_2$ (Lemma~\ref{lem:Nprime}), and each non-invertible simple object $\alpha^i\rho$ of $\mathcal F_{\mathcal C}=\mathcal N'_B$ maps to $\rho'$ under $\mathcal N'_B\to\mathcal N'$: its image has $\mathrm{FPdim}=3d=\mathrm{FPdim}\,\rho'$, and contains no invertible $g$, since $\mathrm{Hom}_{\mathcal N'}(g,\alpha^i\rho)=\mathrm{Hom}_{\mathcal N'_B}(g\otimes B,\alpha^i\rho)=0$, $g\otimes B$ being invertible. So $\nu_n(\rho)$ is computed from the modular data of $\mathcal C$ alone. Evaluating exactly in $\Q(\zeta_{936})$ gives $\nu_0(\rho)=d'$, $\nu_1(\rho)=0$, $\nu_2(\rho)=1$ and $\nu_3(\rho)=e^{2\pi i/3}$ (script \texttt{certify/exact\_nu3.py}). So $\varpi_0$ is a primitive cube root of unity. Which one depends on conventions: with the opposite convention for the centre, the formula of~\cite[Remark~4.3]{NS} uses $\theta^{-n}$ and gives the complex conjugate. Only primitivity is used below. The vertex of $\mathcal G$ has rotation eigenvalue $\varpi=\sigma_2^{-1}(\varpi_0)$, again a primitive cube root of unity.

\emph{Step 4: the vertex generates $\mathcal G$.} Let $\mathcal P\subseteq\mathcal G$ be the pivotal subcategory generated by $\rho$ and $\tau$, i.e.\ the morphisms that are linear combinations of trivalent graphs, and $\overline{\mathcal P}$ its quotient by negligible morphisms.
\begin{itemize}
\item The bigon is a nonzero multiple of $\mathrm{id}_\rho$, because the pairing on $\mathrm{Hom}(\mathbf1,\rho^{\otimes3})$ is nondegenerate; we normalise $\tau$ so that it is $1$. The triangle lies in $\mathrm{Hom}(\mathbf1,\rho^{\otimes3})=\C\tau$. With its three dots placed symmetrically it is invariant under rotation, while rotation acts on $\C\tau$ by $\varpi\ne1$; so the triangle vanishes.
\item $\Delta_\varpi(4,0)=d^5(d-2)\ne0$, so the four graphs of $D_\varpi(4,0)$ are linearly independent, hence a basis of the four-dimensional $\mathrm{Hom}(\mathbf1,\rho^{\otimes4})$. The square lies in this space, and its coefficients are fixed by its pairings with the basis, which are evaluated by loop, bigon and triangle. So the square relation of a twisted cubic category holds in $\mathcal G$. With it, every closed trivalent graph obtained by pairing two graphs with at most six boundary points, in particular every entry of $M^\star_\varpi(6,1)$, has the value computed in~\cite{MPS}.
\item $M^\star_\varpi(6,1)$ has rank $37$ at this $d$ (\cite{MPS}, recomputed independently below). Since the pairing on $\mathrm{Hom}_{\mathcal G}(\mathbf1,\rho^{\otimes6})$ is nondegenerate and this space has dimension $37$, the graphs span it. Composing with the vertex and with caps, the same holds for $n\le6$. So $\overline{\mathcal P}$ is a twisted cubic category with $\dim\mathrm{Hom}(\mathbf 1,\rho^{\otimes5})=11$ and $\dim\mathrm{Hom}(\mathbf 1,\rho^{\otimes6})=37$.
\item By~\cite[Proposition~7.7, Remark~7.3]{MPS}, $\overline{\mathcal P}$ is the trivalent category of $\mathcal H_\varpi$~\cite[Proposition~7.9]{MPS}, i.e.\ the subcategory of $\mathcal H_\varpi$ generated by its trivalent vertex. The proof of~\cite[Proposition~7.9]{MPS} shows that this vertex generates $\mathcal H_\varpi$ (there $\mathcal H'=\mathcal H_\varpi$), so that $\dim\overline{\mathcal P}_n=\dim\mathrm{Hom}_{\mathcal H_\varpi}(\mathbf1,\rho^{\otimes n})$ for all $n$. (For $n\le6$ this also follows from the rank computation above, applied in $\mathcal H_\varpi$.) Since $\mathcal H_\varpi$ and $\mathcal G$ have the same fusion rules,
\[
\dim\overline{\mathcal P}_n=\dim\mathrm{Hom}_{\mathcal H_\varpi}(\mathbf1,\rho^{\otimes n})=\dim\mathrm{Hom}_{\mathcal G}(\mathbf1,\rho^{\otimes n})\ge\dim\mathcal P_n\ge\dim\overline{\mathcal P}_n
\]
for every $n$. So equality holds throughout: $\mathcal P$ has no negligible morphisms, and $\mathrm{Hom}_{\mathcal P}(\mathbf1,\rho^{\otimes n})=\mathrm{Hom}_{\mathcal G}(\mathbf1,\rho^{\otimes n})$ for all $n$. By self-duality of $\rho$ the same holds for all $\mathrm{Hom}(\rho^{\otimes m},\rho^{\otimes n})$.
\end{itemize}

\emph{Step 5: conclusion.} Every simple object of $\mathcal G$ is a summand of some $\rho^{\otimes n}$ ($\alpha\rho,\alpha^2\rho\subset\rho^{\otimes2}$ and $\alpha,\alpha^2\subset\rho^{\otimes3}$), so $\mathcal G$ is the idempotent completion of $\mathcal P$; likewise $\mathcal H_\varpi$ is the idempotent completion of its trivalent category. Hence $\mathcal G\simeq\mathcal H_\varpi$ as pivotal fusion categories, i.e.\ $\mathcal F_{\mathcal C}\simeq\sigma_2\mathcal H_\varpi=\sigma_2\mathcal F$. By Theorem~\ref{thm:Ctw}(2), $\mathcal C\simeq\mathcal Z(\mathcal F_{\mathcal C})\simeq\mathcal Z(\sigma_2\mathcal F)$ as braided fusion categories, and Proposition~\ref{prop:k11}(2) gives $\mathrm{Rep}\,V_{11}\simeq\mathrm{Sem}\boxtimes\mathcal Z(\sigma_2\mathcal F)$.

The equivalence class of $\sigma_2\mathcal F$ does not depend on the extension of $\sigma_2$: two extensions differ by an automorphism fixing $\Q(\zeta_{117})\ni d,\varpi$, which maps $\mathcal H_\varpi$ to a twisted trivalent category with the same $d$ and $\varpi$, hence to $\mathcal H_\varpi$ by~\cite[Proposition~7.7]{MPS}.
\end{proof}

\emph{Independent check of the computations used from~\cite{MPS}.} With our own evaluator for planar trivalent graphs with dots (loop $d$, bigon $1$, triangle $0$, and the square relation derived from the nondegenerate pairing on $D_\varpi(4,0)$), the closed formulas of~\cite[\S7]{MPS} for $\Delta_\varpi(4,0)$, $\Delta_\varpi(5,0)$, $\Delta^\star_\varpi(5,1)$, $\Delta_\varpi(6,0)$ and $\Delta^\star_\varpi(6,1)$ are reproduced at random $d$ (up to powers of $\varpi$), for both $\varpi$. At $d=\frac{3\pm\sqrt{13}}2$, $\Delta_\varpi(6,0)\ne0$ and $\mathrm{rank}\,M^\star_\varpi(6,1)=37$, modulo prime ideals above two primes $Q\equiv1\pmod3$ with $\big(\frac{13}Q\big)=1$, for both $\varpi$; these are lower bounds over $\Q(\sqrt{13},\varpi)$ (scripts \texttt{certify/twisted/}). The input $\Delta_\varpi(6,0)\ne0$ is what the uniqueness proof of~\cite[Proposition~7.7]{MPS} (via the analogue of their Lemma~6.24) needs at this $d$; the rest of that proof is by hand.

The proof uses Ostrik's construction only for the existence of $\mathcal H_\varpi$. That rests on the existence of the $\Z_9+9$ near-group categories~\cite{EG14} and on the Tannakian subcategory $\{a_{3h}\}$ of their centres, whose half-braidings are explicit~\cite[(4.49)--(4.50)]{EG14}. It does not use~\cite[Conjecture~2]{EG14}. Conversely:

\begin{corollary}\label{cor:EGconj}
For $\varpi=\sigma_2^{-1}(\varpi_0)$ as in Theorem~\ref{thm:ident}, $\mathcal Z(\mathcal H_\varpi)\simeq\sigma_2^{-1}(\mathcal C)$ has modular data $\mathfrak D_9=\mathcal D^{\omega}\Hg_3$ with $\omega=2$; for the other value $\bar\varpi$ the modular data is the complex conjugate, $\omega=1$. Hence the factor $\mathcal D$ of the double of the corresponding $\Z_9+9$ near-group category, which is $\mathcal Z(\mathcal H_\varpi)$ by Theorem~\ref{thm:twHI}(2), has the modular data predicted by~\cite[Conjecture~2]{EG14}. In particular the twisted Haagerup--Izumi modular data $\mathcal D^\omega\Hg_3$, $\omega\not\equiv0\pmod3$, is realised by the unitary modular tensor category $\mathcal Z(\mathcal H_\varpi)$.
\end{corollary}
\begin{proof}
The proof uses Theorem~\ref{thm:ident}, which is computer-assisted, and Proposition~\ref{prop:ng}(1),(2).

Galois conjugation commutes with taking the Drinfeld centre, since the centre is built from the structure constants of the category. Let $\varpi=\sigma_2^{-1}(\varpi_0)$ be as in Theorem~\ref{thm:ident}. Applying $\sigma_2^{-1}$ to Theorem~\ref{thm:ident} therefore gives $\mathcal Z(\mathcal H_\varpi)\simeq\sigma_2^{-1}(\mathcal C)$, and the modular data of $\sigma_2^{-1}(\mathcal C)$ is $\sigma_2^{-1}$ applied to that of $\mathcal C$. By Proposition~\ref{prop:ng}(2), the data of $\mathcal C$ is $\sigma_2(\mathfrak D_9)$, with $S$ up to an overall sign, under the bijection of Proposition~\ref{prop:ng}. So $\mathcal Z(\mathcal H_\varpi)$ has modular data $\mathfrak D_9$; the sign of $S$ is fixed by the positivity of $S_{\mathbf1\mathbf1}=1/\sqrt{\dim\mathcal Z(\mathcal H_\varpi)}=1/\dim\mathcal H_\varpi$ in the unitary category $\mathcal Z(\mathcal H_\varpi)$. By Proposition~\ref{prop:ng}(1) this is $\mathcal D^\omega\Hg_3$ with $\omega=2$. The category $\mathcal H_{\bar\varpi}$ is the complex conjugate of $\mathcal H_\varpi$ (the complex conjugate of $\mathcal H_\varpi$ is a twisted trivalent category with the same $d$ and rotation eigenvalue $\bar\varpi$, hence is $\mathcal H_{\bar\varpi}$ by~\cite[Proposition~7.7]{MPS}), so the two values of $\varpi$ give complex-conjugate data, $\omega=2$ and $\omega=1$. For a given $\mathcal N$ the correct one is fixed by the twists. By Izumi's analysis of the tube algebra~\cite[Theorem~4.10, \S5]{IzumiI},~\cite{Izumi}, recalled in~\cite[(4.56)]{EG14}, the $T$-matrix of $\mathcal Z(\mathcal N)$ is known exactly; in particular the twist of $c_{k,-k}\in\mathcal D$ is $\langle k,-k\rangle=e^{-2\pi imk^2/9}$, where $\langle g,h\rangle=e^{2\pi imgh/9}$, $m=\pm1$, is the pairing of $\mathcal N$~\cite[Table~2]{EG14}. These are the entries of $T$ on $c_{k,-k}$ in $\mathfrak D_9$ if $m=1$, and in its complex conjugate $\overline{\mathfrak D_9}$ if $m=-1$, up to the common factor $\beta$ of~\cite[Proposition~7]{EG14}. Since $\mathcal D\simeq\mathcal Z(\mathcal H_\varpi)$ for the appropriate $\varpi$ is a braided equivalence, and the $T$-matrices of $\mathfrak D_9$ and $\overline{\mathfrak D_9}$ differ ($e^{\mp2\pi ik^2/9}$), the modular data of $\mathcal D$ is $\mathfrak D_9$ if $m=1$ and $\overline{\mathfrak D_9}$ if $m=-1$. By Proposition~\ref{prop:ng}(1), these are $\mathcal D^\omega\Hg_3$ with $\omega=2$ and $\omega=1$.

The unitary category $\mathcal N_A$ of Theorem~\ref{thm:twHI} is a twisted Haagerup category. It has the Haagerup--Izumi fusion rules for $\Z_3$, $d=\frac{3+\sqrt{13}}2$ and nontrivial pointed associator (Theorem~\ref{thm:twHI}(3),(4)), so its non-invertible $\rho$ with a nonzero $\tau\in\mathrm{Hom}(\mathbf1,\rho^{\otimes3})$ satisfies the dimension conditions of Step~2 of Theorem~\ref{thm:ident}. If its rotation eigenvalue were $1$, the category generated by $\tau$ would be an untwisted trivalent category with $\dim\mathcal C_4=4$, $\dim\mathcal C_5=11$, $\dim\mathcal C_6=37$, hence the Grossman--Snyder category $\mathcal H_3$ or its Galois conjugate~\cite[Theorem~C]{MPS}, whose invertible objects have trivial associator; this contradicts Theorem~\ref{thm:twHI}(4). So the rotation eigenvalue is primitive, and Steps~4--5 of the proof of Theorem~\ref{thm:ident}, applied to $\mathcal N_A$ (where unitarity is available), give $\mathcal N_A\simeq\mathcal H_\varpi$ for that eigenvalue $\varpi$. By Theorem~\ref{thm:twHI}(1),(2), which do not involve the modular data of $\mathcal D$, the factor $\mathcal D$ of $\mathcal Z(\mathcal N)$ is $\mathcal Z(\mathcal N_A)=\mathcal Z(\mathcal H_\varpi)$, so it has modular data $\mathfrak D_9$ or $\overline{\mathfrak D_9}$ as above. The other factor $\mathcal P\simeq\mathcal C(\Z_9,Q)$ has the pointed modular data given by~\cite[(4.49)--(4.50)]{EG14}, which is the pointed factor of $MD_{\Z_9,\Z_{13}}(Q,Q')$. With the factorisation of $MD_{\Z_9,\Z_{13}}(Q,Q')$ recalled before Proposition~\ref{prop:ng}, $\mathcal Z(\mathcal N)\simeq\mathcal P\boxtimes\mathcal D$ has modular data $MD_{\Z_9,\Z_{13}}(Q,Q')$ up to relabelling, which is~\cite[Conjecture~2]{EG14} for $\Z_9$. Finally, $\mathcal H_\varpi=\mathcal N_A$ is a unitary fusion category (Theorem~\ref{thm:twHI}(2)), so $\mathcal Z(\mathcal H_\varpi)$ is a unitary modular tensor category~\cite{MugerII}.
\end{proof}

\begin{remark}
\cite[Proposition~7.7]{MPS} also leaves open a twisted cubic category $\mathcal Q_\varpi$ with $d=-1$ \cite[Question~7.10]{MPS}. It plays no role here, since $d\ne-1$. Without the bootstrap, Theorem~\ref{thm:twHI} together with~\cite[Conjecture~2]{EG14} would give only a coincidence of modular data. Theorem~\ref{thm:ident} identifies the braided categories, using the étale algebra $A_2$ supplied by the vertex operator algebra $U_2$.
\end{remark}

\begin{lemma}\label{lem:Dalg}
In the factor $\mathcal D\subset\mathcal Z(\mathcal N)$ of Theorem~\ref{thm:twHI}, the object $a_0\oplus b_0$ is a connected étale algebra, and $\mathcal D_{a_0\oplus b_0}\simeq\mathcal N$.
\end{lemma}

\begin{proof}
The composite $\mathcal D\hookrightarrow\mathcal Z(\mathcal N)\to\mathcal N$ is a central functor. In the tube algebra, $a_g$, $b_g$, $c_{g,h}$ and $d$ are the half-braidings on $\alpha_g$, $\rho\oplus\alpha_g$, $\rho\oplus\alpha_g\oplus\alpha_h$ and $\rho$~\cite[\S4.3]{EG14}. So the images of $c_{k,-k}$ and $d$ contain every $\alpha_g$ and $\rho$, and the functor is dominant. Its right adjoint applied to $\mathbf1$ is $\bigoplus_x\dim\mathrm{Hom}(F(x),\mathbf1)\,x=a_0\oplus b_0$. That is a connected étale algebra with $\mathcal D_{I(\mathbf1)}\simeq\mathcal N$~\cite{DMNO}.
\end{proof}

\section{A conjecture for odd \texorpdfstring{$k$}{k}}\label{sec:conj}

Let $k\ge5$ be odd and $n=k-2$. For odd $n$, Evans--Gannon attach to the near-group $\Z_n+n$ the modular data $MD_{\Z_n,\Z_{n+4}}$~\cite[Proposition~7(b)]{EG14}. It is the pointed data of a metric group $(\Z_n,q)$ times a factor $\mathfrak D_n$ of rank $n+3$, on primaries $a_0,b_0$, $c_\kappa$ ($1\le\kappa\le\frac{n-1}2$), $d_\gamma$ ($1\le\gamma\le\frac{n+3}2$); for $n=9$ this is $\mathfrak D_9$ of \S\ref{ssec:ng}. The characters of $V_k$ are those of \S\ref{sec:chars} for every $k$, and $M_{\chi_1}$ and $M_{\chi_3}$ are defined as in \S\ref{sec:k11}.

\begin{conjecture}\label{conj:odd}
Let $k\ge5$ be odd and $n=k-2$. Here $M_{\chi_1}$ and $M_{\chi_3}$ denote the simple $V_k$-modules with characters $\chi_1$ and $\chi_3$ of \S\ref{sec:chars}, whose existence is part of~(1); thus $M_{\chi_3}$ is the simple module of lowest weight $0$ other than $V_k$, and $M_{\chi_1}$ is the simple current with the character of $\bigoplus_{m\ge1}L(h_{2m,2m})$, as for $k=11$.
\begin{enumerate}
\item $V_k$ is strongly rational, and $\mathrm{Rep}\,V_k\simeq\mathcal S\boxtimes\mathcal C_k$ as braided categories, where $\mathcal S=\langle M_{\chi_1}\rangle$ is $\mathrm{Sem}$ for $n\equiv1\pmod4$ and $\overline{\mathrm{Sem}}$ for $n\equiv3\pmod4$ (\S\ref{sec:rat}), and $\mathcal C_k$, the centraliser of $M_{\chi_1}$, has rank $n+3$ and the modular data of a Galois conjugate of $\mathfrak D_n$.
\item $V_k\oplus M_{\chi_3}$ is a simple vertex operator algebra, unique up to isomorphism.
\item The category $\mathcal N_k$ of modules over $A=V_k\oplus M_{\chi_3}$ in $\mathcal C_k$ has near-group fusion rules $\Z_n+n$, and $\mathcal Z(\mathcal N_k)\simeq\mathcal C_k\boxtimes\mathcal C(\Z_n,q)^{\rm rev}$.
\item If $n=\nu^2$, i.e.\ $k=\nu^2+2$ with $\nu$ odd, then $\mathcal C_k\simeq\mathcal Z(\mathcal F_k)$, where $\mathcal F_k$ has Haagerup--Izumi fusion rules for $\Z_\nu$, nontrivial pointed associator, and non-invertible simple objects of dimension $\frac12(\nu-\sqrt{\nu^2+4})$. Its centre has the modular data of a Galois conjugate of $\mathcal D^{\omega}\Hg_\nu$, $\omega\not\equiv0\pmod\nu$.
\end{enumerate}
\end{conjecture}

\paragraph{Evidence.}
\begin{itemize}
\item For $k=11$ all four parts are proved here: (1) is Theorem~\ref{thm:mod} with Proposition~\ref{prop:k11}, (2) is Theorem~\ref{thm:U2}, (3) is Lemma~\ref{lem:Nprime} with Theorem~\ref{thm:Ctw}(1), and (4) is Theorem~\ref{thm:ident}.
\item For $k=5,7,9$, strong rationality and the simple modules, characters and modular data in part (1) are Theorem~\ref{thm:mod} (with two steps numerical); the splitting off of the semion and the identification with $\mathfrak D_n$ are checked numerically (next item).
\item For every odd $k$ with $5\le k\le15$, the modular data of the centraliser of $\chi_1$, computed from the characters, equals that of $\mathfrak D_n$ after $\sqrt{n(n+4)}\mapsto-\sqrt{n(n+4)}$ and relabelling. For $5\le k\le21$ and for $k=27,51$, the conditions on $M_{\chi_3}$ that are visible in the modular data hold: $\theta=1$, $N_{\chi_3\chi_3}^{V}=N_{\chi_3\chi_3}^{\chi_3}=1$, $\nu_2=\nu_3=1$. For $k=\nu^2+2$, $\nu=3,5,7$, the dimensions in $\mathcal C_k$ are those of a Galois conjugate of $\mathcal D^\omega\Hg_\nu$. These checks are numerical (scripts in \texttt{conj/}).
\item Fusion categories with near-group fusion rules $\Z_n+n$ and $\dim\rho=\frac12(n-\sqrt{n^2+4n})$, as in part~(3), exist for every $n$: Radchenko--Wheeler~\cite[Corollary~6]{RW} construct them from finite quantum dilogarithms on $\Z_n$, by the Leavitt-algebra form of the endomorphism construction of~\cite{EG17}. Part~(3) asserts that $\mathcal N_k$ is one of these. The existence of unitary near-groups $\Z_n+n$ for every $n$ is open~\cite[\S6]{RW}.
\end{itemize}

\begin{remark}
Part (3) explains why the Haagerup--Izumi values are special. By (3), $\mathcal C_k$ and $\mathcal C(\Z_n,q)$ have the same Witt class~\cite{DMNO}. A pointed braided fusion category $\mathcal C(\Z_n,q)$, $n$ odd, is Witt-trivial exactly when it has a Lagrangian subgroup, i.e.\ when $n$ is a square, the subgroup being $\sqrt n\,\Z_n$. So, given (3), $\mathcal C_k$ is a Drinfeld centre if and only if $k=\nu^2+2$, which are Gang--Kim's Haagerup--Izumi values (Observation~\ref{obs:t}). In that case de-equivariantisation by $\sqrt n\,\Z_n\cong\Z_\nu$ gives (4), as in Theorem~\ref{thm:Ctw}. For $\nu=3$ the identification of $\mathcal F_k$ with a Galois conjugate of Ostrik's category uses~\cite{MPS}; for $\nu\ge5$ no such classification is available. For the untwisted Haagerup--Izumi categories, Drinfeld centres are computed by Grossman--Izumi~\cite{GI}.
\end{remark}

\section*{Use of AI}

This work came  from conversations with Claude who found the 1993 paper of Flohr and then developed the strategy and proofs  for the technicalities on the VOA side (including   $C_2$-cofiniteness and  rationality) and assisted in writing the manuscript.
ChatGPT provided an audit.  I am responsible for the final document.

\section*{Acknowledgements}

I would like to thank Terry Gannon for visiting and delivering an illuminating, instructive detailed series of lectures on Vertex Operator Algebras in Cardiff in the summer of 2014 and for numerous discussions on the Haagerup subfactor over almost twenty years. I am also grateful to the Leverhulme Foundation for the financial support of  an Emeritus Fellowship.
This work was supported in part by a grant of access to OpenAI models through the ChatGPT for Academic Researchers programme.

\appendix
\section{Computational details}\label{app:comp}

\subsection*{Status of the computer-assisted steps}\label{app:status}
Here \emph{computer-assisted} means that a finite computation, done modulo primes or in floating point, supplies an input to a proof. The inputs are of two kinds.
\begin{itemize}
\item[(a)] \emph{Rank and non-vanishing statements} (a matrix has rank $r$, a number is nonzero). If the inputs are exact rationals and every Gram matrix used has the same rank modulo $p$ as over $\Q$ (invertible, or as in (G1), (G2) of Lemma~\ref{lem:good}), the computation modulo $p$ is the reduction of the rational one, and rank can only drop under reduction. Such statements are then proofs as they stand.
\item[(b)] \emph{Vanishing and membership statements} (a residual is zero, a vector lies in a span). These do not pass from $\mathbb F_p$ to $\Q$, and need an exact computation. Floating-point steps are also of this kind.
\end{itemize}
For Theorem~\ref{thm:ident} the inputs stand as follows.
\begin{itemize}
\item[C1.] Theorem~\ref{thm:mod} for $k=11$. (i) $C_2$-cofiniteness: $R_{2h+4}=0$ is a rank statement, of kind (a); good reduction is checked at every level (\S\ref{sec:c2}). (ii) Zhu's algebra: Lemmas~\ref{lem:good} and~\ref{lem:lift} turn the modular computation into a proof; the checks (G1)--(G3) are part of the run (Table~\ref{tab:zhu}). (iii) Lemma~\ref{lem:cyc} and the commutant step of Theorem~\ref{thm:mod}: exact. (iv) $\mathcal N_{11}$: exact, from the Jacobi identity of $V_{11}$ over $\Q$ (Table~\ref{tab:N}). (v) The sign of $\langle v,W_{(17)}v\rangle$: from Lemma~\ref{lem:lift} (Remark~\ref{rem:exact}).
\item[C2.] Proposition~\ref{prop:boot}(1),(2) and Lemma~\ref{lem:G}: exact over $\Q$ (Remark~\ref{rem:exact}). Proposition~\ref{prop:A2}(2): $\theta$, $\nu_2$, $\nu_3$ of $\chi_3$ exact.
\item[C3.] Lemma~\ref{lem:a3}: the $W$-couplings of $M_{a_3}$ are exact, and the rank statements at $P_1$ hold over $\Q$ (Remark~\ref{rem:exact}).
\item[C4.] Proposition~\ref{prop:ng}(1),(2): exact in $\Q(\zeta_{936})$.
\item[C5.] From~\cite{MPS}: the computations $\Delta_\varpi(6,0)\ne0$ and $\mathrm{rank}\,M^\star_\varpi(6,1)\ge37$ at this $d$ are reproduced independently (Theorem~\ref{thm:ident}); the hand arguments of Propositions~7.7 and~7.9 are used as published.
\end{itemize}
Each of C1--C5 is certified by the scripts listed below.
Proposition~\ref{prop:k11}(4),(5) are consistency checks once $\mathcal C$ is a centre (Corollary~\ref{cor:centre}). Of Proposition~\ref{prop:mi} only the following is used (in Proposition~\ref{prop:A2}(1) and Lemma~\ref{lem:Nprime}), and it is exact: the commutant of $S$ and $T$ on $\mathcal C$ is spanned by $I,P,R$ (\texttt{exact\_mod.py}), so the modular invariant of $A_2$, which has $Z_{0,0}=Z_{0,\chi_3}=1$ and $Z_{0,a_3}=0$, is $I+P=Z_2$. With C1--C4, every computer-assisted input to Theorem~\ref{thm:ident} at $k=11$ is either an exact computation over $\Q$ or $\Q(\zeta_{936})$, or a rank or non-vanishing statement modulo a prime with good reduction checked. What remains is independent checking of the code, which is supplied with this paper.

\subsection*{Details}
\paragraph{Characters.} The Nahm sums were expanded exactly with integer arithmetic and compared with the Virasoro and theta formulas.
\begin{itemize}
\item Vacuum identity: $k=3,\dots,11$ to $q^{34}$.
\item $\chi_1,\chi_2,\chi_3$: $k=4,\dots,7$ to $q^{40}$.
\end{itemize}

\paragraph{Minimal-model weights.} For $\mathfrak{so}_5$ at $(p,q)=(6,5)$ and $G_2$ at $(7,8)$, the modules are labelled by $\lambda\in P_+^{p-h^\vee}$ and $\mu\in P_+^{\vee,q-h}$ (levels measured against the highest root). The lowest weights are $\frac{|\lambda+\rho-x(\mu+\rho^\vee)|^2-|\rho-x\rho^\vee|^2}{2x}$ with $x=p/q$.

\paragraph{$C_2$ computations.}
\begin{itemize}
\item Virasoro Verma modules were handled in a PBW basis, using the vacuum-module basis (parts $\ge2$) for $L(c,0)$. Gram matrices were built recursively, via $G_\ell[\mu,\lambda]=\sum_\nu (L_{\mu_1}L_{-\lambda})_\nu\,G_{\ell-\mu_1}[\mu',\nu]$, over $\mathbb F_p$.
\item The largest Gram solves were in the vacuum module at levels $38$ and $40$ for $k=11$ (dimensions $4378$ and $6153$), taking about one and five minutes each with python-flint.
\item Residues of $\mathcal N$ modulo $p_1=33554393$: $k=4{:}\ 31290082$, $5{:}\ 11676315$, $6{:}\ 14316214$, $7{:}\ 25955273$, $8{:}\ 561345$, $9{:}\ 9456979$, $10{:}\ 2118399$, $11{:}\ 27947361$.
\item Modulo $p_2=33554383$: $k=5{:}\ 5119712$, $6{:}\ 26007907$, $7{:}\ 4183198$, $8{:}\ 22675452$, $9{:}\ 25987538$, $10{:}\ 22160861$, $11{:}\ 19446660$.
\end{itemize}

\paragraph{Scripts} (supplied separately):
\begin{itemize}
\item \texttt{nahm.py}, \texttt{check.py}, \texttt{others.py}: Nahm sums and character identities.
\item \texttt{smat.py}: $S$-matrix.
\item \texttt{b2.py}, \texttt{g2.py}: minimal-model weights.
\item \texttt{vir.py}, \texttt{walg.py}, \texttt{jac.py}, \texttt{j4.py}: Virasoro machinery and Jacobi identities.
\item \texttt{c2.py}, \texttt{joint.py}: full $C_2$ computation.
\item \texttt{red.py}: reduced $C_2$ computation.
\item \texttt{zhu.py}, \texttt{full.py}: Zhu relations, checks (G1)--(G3), Gr\"obner bases and points (\S\ref{sec:zhu}). Logs: \texttt{full\_p1.log}, \texttt{full\_p2.log}.
\item \texttt{test3.py}: the $k=3$ check against $M(3,5)\otimes M(2,5)$.
\item \texttt{modcheck.py}, \texttt{lead.py}: cyclic span, commutant, and lowest-exponent bounds (\S\ref{sec:rat}).
\item \texttt{k11fus.py}: fusion rules, the semion, the centraliser and the Lagrangian vector at $k=11$ (\S\ref{sec:k11}).
\item \texttt{mp/nw\_scalar.py} ($k=11$), \texttt{mp/nw\_exact\_k.py} ($4\le k\le10$): $\mathcal N$ exactly over $\Q$.
\item \texttt{certify/}: the exact certificates. \texttt{exact\_mod.py} (commutant and cyclic span at $k=11$, Theorem~\ref{thm:mod}); \texttt{exact\_ng.py}, \texttt{exact\_ng1.py} (Proposition~\ref{prop:ng}); \texttt{exact\_fs.py} (Proposition~\ref{prop:A2}(2)); \texttt{exact\_nu3.py} ($\nu_n(\rho)$, Theorem~\ref{thm:ident}, Step~3); \texttt{pi\_radical.py}, \texttt{pi\_radical\_k.py} ($\pi(s)=0$, \S\ref{sec:zhu}); \texttt{c1/red\_check.py}, \texttt{c1/red\_check\_k.py} (good reduction in \S\ref{sec:c2}); \texttt{twisted/} (the computations from~\cite{MPS}); \texttt{mp/} (the exact constants of Remark~\ref{rem:exact}: \texttt{exact\_consts.py}, \texttt{exact\_G.py}, \texttt{exact\_gb.py}, \texttt{exact\_gbA.py}, \texttt{exact\_jac.py}, \texttt{exact\_ward.py}, \texttt{reconB.py}); \texttt{conj/} (the numerical evidence in \S\ref{sec:conj}).
\item \texttt{neargroup/}: the near-group comparison, modular invariants, Frobenius--Schur indicators and braiding signs (\S\S\ref{ssec:mi}--\ref{ssec:gap}).
\item \texttt{existence/}: Lemma~\ref{lem:G} (\texttt{gfun.py}, with logs \texttt{g26\_p1.log}, \texttt{g26\_p2.log}) and the bound in the proof of Theorem~\ref{thm:U2} (\texttt{fuchs\_bound.py}).
\item \texttt{bootstrap/} (\texttt{boot\_p1}, \texttt{boot\_p2}, \texttt{boot\_k3}): Proposition~\ref{prop:boot}. \texttt{mod.py}, \texttt{fast2.py}: the staggered model of $M_{\chi_3}$. \texttt{stage1.py}, \texttt{stage2.py}, \texttt{uu.py}: the $W$-couplings. \texttt{ext.py}: the model of $V_{11}\oplus M_{\chi_3}$ with general primary fields. \texttt{stage3a.py}, \texttt{stage3b.py}: the $W$--$v$ Ward identities. \texttt{stage4w.py}: locality. The $W$-eigenvalue on the top of $M_{\chi_3}$ comes from $[W]*[W]$ in Zhu's algebra (\texttt{w3.py}). Each stage fits polynomial residuals by interpolation at random points and checks the fit at extra points. \texttt{lagrangian/} (at $P_1$): the model of $M_{a_3}$ (\texttt{modA.py}), the $v$-couplings between $M_{\chi_3}$ and $M_{a_3}$ (\texttt{stageB.py}, \texttt{stageB2.py}) and the Lagrangian test (\texttt{stageC.py}); used in Lemma~\ref{lem:Nprime} and in the proof of Theorem~\ref{thm:U2}.
\end{itemize}

\end{document}